\documentclass[11pt]{amsart}
\usepackage{latexsym,amssymb,epsf, color}
\usepackage[centertags]{amsmath}

\usepackage{hyperref}

\usepackage{graphicx}

\newcounter{bz}

\newcounter{bx}
\addtocounter{bx}{1}

\newcommand{\halmos}{\rule{1ex}{1.4ex}}

\newtheorem{theorem}{Theorem}[section]
\newtheorem{proposition}[theorem]{Proposition}
\newtheorem{corollary}[theorem]{Corollary}
\newtheorem{conjecture}[theorem]{Conjecture}
\newtheorem{lemma}[theorem]{Lemma}

\newtheorem{remark}[theorem]{Remark}
\newtheorem{defi}[theorem]{Definition}
\newcommand{\bt}{\begin{theorem}}
\newcommand{\et}{\end{theorem}}
\newcommand{\bl}{\begin{lemma}}
\newcommand{\el}{\end{lemma}}
\newcommand{\bp}{\begin{proposition}}
\newcommand{\ep}{\end{proposition}}
\newcommand{\bcor}{\begin{corollary}}
\newcommand{\ecor}{\end{corollary}}
\newcommand{\br}{\begin{remark}}
\newcommand{\er}{\end{remark}}
\newcommand{\bcon}{\begin{conjecture}}
\newcommand{\econ}{\end{conjecture}}
\newcommand{\bd}{\begin{defi}}
\newcommand{\ed}{\end{defi}}

\newcommand{\bi}{\bigskip}

\newcommand{\Rand}[1]{\marginpar{#1}}
\renewcommand{\Rand}[1]{}

\newcommand{\be}[1]{\Rand{\vspace{0.6cm}\tt #1}\begin{equation}\label{#1}}
\newcommand{\ee}{\end{equation}}

\newcommand{\bea}[1]{\Rand{\vspace{0,7cm}\tt #1\vspace{-0,7cm}}\begin{eqnarray}\label{#1}}
\newcommand{\eea}{\end{eqnarray}}
\newcommand{\qad}{\qquad \mbox{$\square$}}

\newcommand{\bean}{\begin{eqnarray*}}
\newcommand{\eean}{\end{eqnarray*}}

\makeatletter\@addtoreset{equation}{section}
\makeatother\def\theequation{\thesection.\arabic{equation}}

\newcommand{\belC}[2]{\Rand{\vspace{0.6cm}\tt #1}\begin{lemma}[#2]\label{#1}}
\newcommand{\beP}[1]{\Rand{\vspace{0.6cm}\tt #1}\begin{proposition}\label{#1}}
\newcommand{\beD}[1]{\Rand{\vspace{0.6cm}\tt #1}\begin{definition}\label{#1}}
\newcommand{\beT}[1]{\Rand{\vspace{0.6cm}\tt #1}\begin{theorem}\label{#1}}
\newcommand{\beC}[1]{\Rand{\vspace{0.6cm}\tt #1}\begin{corollary}\label{#1}}
\newcommand{\beCj}[1]{\Rand{\vspace{0.6cm}\tt #1}\begin{Conjecture}\label{#1}}

\def\eps{\epsilon}
\def\dd{{\rm d}}
\def \qed {\nopagebreak{\hspace*{\fill}$\halmos$\medskip}}

\def\CC{\mathcal{C}}
\def\CD{\mathcal{D}}

\def\CK{\mathcal{K}}
\def\CM{\mathcal{M}}

\def\CL{\mathcal{L}}

\def\CR{\mathcal{R}}

\def\CW{\mathcal{W}}
\def\CX{\mathcal{X}}

\def\E{\mathbb{E}}

\def\G{\mathbb{G}}
\def\K{\mathbb{K}}
\def\M{\mathbb{M}}
\def\N{\mathbb{N}}
\def\P{\mathbb{P}}

\def\BS{\mathbb{S}}
\def\R{\mathbb{R}}

\def\U{\mathbb{U}}

\def\Z{\mathbb{Z}}

\newcommand{\sm} {\smallskip}

\newcommand{\ve}{\varepsilon}
\newcommand{\D}{\displaystyle}

\newcommand{\tto}{{_{\D \Longrightarrow \atop t \to \infty}}}

\newcommand{\nto}{{_{\D \Longrightarrow \atop n \to \infty}}}

\newcommand{\ttO}{{_{\D \longrightarrow \atop t \to \infty}}}

\newcommand{\ntoo}{{_{\D \longrightarrow \atop n \to \infty}}}

\newcommand{\wh}{\widehat}

\newcommand{\asto}[1]{\underset{{#1}\to\infty}{\longrightarrow}}
\newcommand{\Asto}[1]{\underset{{#1}\to\infty}{\Longrightarrow}}
\newcommand{\astoo}[1]{\underset{{#1}\to 0}{\longrightarrow}}
\newcommand{\Astoo}[1]{\underset{{#1}\to 0}{\Longrightarrow}}

\newcommand{\intl}{\int\limits}
\newcommand{\liml}{\lim\limits}
\newcommand{\suml}{\sum\limits}

\newcommand{\uur}{\underline{\underline{r}}}

\newcommand{\ux}{\underline{x}}
\newcommand{\uy}{\underline{y}}
\newcommand{\uk}{\underline{k}}

\newcommand{\uu}{\underline{u}}
\newcommand{\uv}{\underline{v}}
\newcommand{\uB}{\underline{B}}
\newcommand{\wt}{\widetilde}

\begin{document}

\title[Continuum Space Limit of IFV Genealogies]{Continuum Space Limit of the Genealogies of Interacting Fleming-Viot Processes on $\Z$}

\author{Andreas Greven$^{\,1}$, Rongfeng Sun$^{\,2}$, Anita Winter$^{\,3}$}

\thispagestyle{empty}
\date{\today}

\maketitle

\footnotetext[1]{Universit\"at Erlangen-N\"urnberg, Germany. Email: greven@mi.uni-erlangen.de}

\footnotetext[2]{National University of Singapore, Singapore.  Email: matsr@nus.edu.sg}

\footnotetext[3]{Universit\"at Duisburg-Essen, Germany. Email: anita.winter@uni-due.de}

\begin{abstract}
We study the evolution of genealogies of a population of individuals, whose type frequencies result in an interacting Fleming-Viot process on $\Z$. We construct and analyze the genealogical structure of the population in this genealogy-valued Fleming-Viot process as a marked metric measure space, with each individual carrying its spatial location as a mark. We then show that its time evolution converges to that of
the genealogy of a continuum-sites stepping stone model on $\R$, if space and time are scaled diffusively. We construct the genealogies of the continuum-sites stepping stone model as functionals of the Brownian web, and furthermore, we show that its evolution solves a martingale problem. The generator for the continuum-sites stepping stone model has a singular feature: at each time, the resampling of genealogies only affects a set of individuals of measure $0$. Along the way, we prove some negative correlation inequalities for coalescing Brownian motions, as well as extend the theory of marked metric measure spaces
(developed recently by Depperschmidt, Greven and Pfaffelhuber~\cite{DGP12}) from the case of probability measures to measures that are finite on bounded sets.
\end{abstract}

\smallskip


\smallskip

\vspace{12pt}

{\tiny
\tableofcontents
}

\newcounter{secnum}
\setcounter{secnum}{\value{section}}
\setcounter{section}{0}
\setcounter{secnumdepth}{3}
\setcounter{equation}{0}
\renewcommand{\theequation}{\mbox{\arabic{secnum}.
\arabic{equation}}}
\numberwithin{equation}{section}

\section{Introduction and main results} \label{S:intro}

In the study of spatial population models on discrete geographic spaces
(for example $\Z^d$), such as branching processes, voter models, or interacting Fisher-Wright diffusions (Fleming-Viot models),
the technique of passing to the {\em spatial continuum limit} has proven to be useful in gaining insight into the
qualitative behaviour of these processes. A key example is branching random walks on $\Z^d$, leading to the Dawson-Watanabe process
\cite{D77} on $\R^d$ and Fisher-Wright diffusions; catalytic branching and mutually catalytic branching on $\Z$, leading to SPDE on $\R$
\cite{KS88, EF96, DP98, DEF02a, DEF02b}.  The goal of this paper is to carry out this program at the level of {\em genealogies},
rather than just type or mass configurations. We focus here on {\em interacting Fleming-Viot models} on $\Z$.

\subsection{Background and Overview}
We summarize below the main results of this paper, recall some historical background, as well as state some open problems.
\medskip

\noindent
{\bf Summary of results.}
The purpose of this paper is twofold. On the one hand, we want to understand the formation of large local one-family clusters in
Fleming-Viot populations on the geographic space $\Z^1$, by taking a space-time continuum limit of the genealogical configurations equipped with types.
On the other hand we use this example to develop the theory of {\em tree-valued dynamics} via {\em martingale problems} in some new directions.
 In particular, this is the first study of a {\em tree-valued dynamics} on an unbounded geographical space with infinite sampling measure, which
requires us to extend both the notion of {\em marked metric measure spaces} in~\cite{GPW09, DGP11} and the martingale problem formulations
 in~\cite{GPWmp13, DGP12} to {\em marked} metric measure spaces with infinite sampling measures that are boundedly finite
  (i.e., finite on bounded sets).
\medskip

Here is a summary of our main results:
\begin{itemize}
\item[(1)] We extend the theory of {\em marked metric measure spaces}~\cite{GPW09, DGP12} from probability sampling measures to infinite sampling measures that are {\em boundedly finite}, which serve as the state space of marked genealogies of spatial population models. See Section~\ref{S:Intrommm}.
\item[(2)] We characterize the evolution of the genealogies of interacting Fleming-Viot (IFV) models by well-posed martingale problems on spaces of {\em marked ultra\-metric  measure spaces}. See Section~\ref{S:IntroIFV}.
\item[(3)] We give a {\em graphical construction} of the {\em spatial continuum limit} of the IFV genealogy process, which is the genealogy process of the so-called {\em Continuum Sites Stepping-stone Model} (CSSM), taking values in the space of ultra\-metric measure spaces with {\em spatial marks} and an {\em infinite total population}. The graphical construction is based on the (dual) Brownian web~\cite{FINR04}.
    The CSSM genealogy process has the feature that, as soon as $t>0$, the process enters a regular subset of the state space that is not closed under the topology. Only on this subset we can evaluate the action of the operator of the martingale problem in its action on test functions. The nice aspect is that the set of these states are preserved under the dynamic. \color{black} This leads {{}to a singular structure with} complications for the associated martingale problem and for the study of continuity of the process at time $0$. See Section~\ref{S:IntroCSSM}.
\item[(4)] We prove that under suitable scaling, the IFV genealogy processes converge to the CSSM genealogy process. The proof is based on duality with spatial coalescents, together with a novel approach of controlling the genealogy structure using a weaker convergence result on the corresponding measure-valued processes, with measures on the geographic and type space (with no genealogies). See Section~\ref{S:Introconv}.
\item[(5)] We show that the CSSM genealogy process solves a martingale problem with a singular generator. More precisely, {{}the generator action involves individuals, which are {\em not typical} under the sampling measure, so that the dynamic is driven by atypical individuals at atypical locations.} In particular, the generator is only defined on a regular subset of the state space. See Section~\ref{S:IntroCSSMmart}.
\item[(6)] We prove some negative correlation inequalities for coalescing Brownian motions, which are of independent interest. See Appendix~\ref{A:Corr}.
\end{itemize}
Besides the description of the genealogies of the current population, we also prepare the ground for the treatment of all individuals ever alive, i.e.\ {\em fossils}, moving from the state space of marked ultra\-metric measure spaces to the state space of marked measure $\R$-trees, which will be carried out elsewhere.

\bi

\noindent
{\bf History of the problem:}  Why are we particularly interested in {\em one}-dimensional geographic spaces
for our scaling results?
Many interacting spatial systems that model evolving populations, i.e., Markov processes with state spaces
$I^G (I=\R,\N, [0,1]$, etc., and $G=\Z^d$ or the hierarchical group $\Omega_N$) that evolve
by a migration mechanism between sites and a stochastic mechanism acting locally at each site, exhibit a
{\em dichotomy} in their {\em longtime behavior}.
For example, when $G=\Z^d$ and the migration is induced by the simple symmetric random walk:
in dimension $d \leq 2$, one observes convergence to laws concentrated on the traps of the dynamic; while in $d \geq 3$, nontrivial equilibrium states are approached and the extremal invariant measures are spatially homogeneous ergodic measures
characterized by the intensity of the configuration. Typical examples include the
voter model, branching random walks, spatial Moran models, or systems of interacting
diffusions (e.g., Feller, Fisher-Wright or Anderson diffusions). One obtains universal dimension-dependent
scaling limits for these models if an additional continuum spatial limit is taken, resulting in, for example, super Brownian motion
(see Liggett \cite{Lig85} or Dawson \cite{D93}).

In the low-dimension regime, the cases $d=1$ and $d=2$ are very different. In
$d=2$, one observes for example in the voter model the formation of mono-type clusters
on spatial scales $t^{\alpha/2}$ with a random $\alpha \in [0,1]$, a phenomenon
called diffusive clustering (see Cox and Griffeath \cite{CG86}). In the one-dimensional
voter model, the clusters  have an extension of a fixed order of magnitude but
exhibit random factors in that scale. More precisely, in space-time scales
$(\sqrt{t}, t)$ for $t \to \infty$, we get annihilating Brownian motions. Similar results
have been obtained for low-dimensional branching systems
(Klenke \cite{K00}, Winter \cite{W02}), systems of interacting Fisher-Wright diffusions
(Fleischmann and Greven \cite{FG94},  \cite{FG96} and subsequently \cite{Z03}, \cite{DEFKZ00} )
and for the Moran model in $d=2$ (Greven, Limic, Winter \cite{GLW05}).

In all these models, one can go further and study the complete space-time genealogy structure of
the cluster formation and describe this phenomenon asymptotically by the {\em spatial continuum limit}.
In particular, the description for the one-dimensional voter model
can be extended to the complete space-time genealogy structure, obtaining as scaling limit
the {\em Brownian web} \cite{Arr79, A81, TW98,FINR04} (see Appendix~\ref{S:web}, and the recent survey
\cite{SchertzerSunSwart}). More precisely, the Brownian web is defined by considering instantaneously
coalescing one-dimensional Brownian motions starting from every space-time point in $\R \times \R$. It arises as the diffusive scaling limit of
continuous time coalescing simple symmetric random walks starting from every space-time point in $\Z \times \R$, which represent the space-time genealogies of the voter model. This is analogous to the study of historical process for branching processes, which approximates the ancestral paths of branching random walks by that of super Brownian motion (see e.g.~\cite{DP91, FG96, GLW05}).

The basic idea behind all this is that, we can identify the {\em genealogical relationship} between the individuals of the population
living at different times and different locations.  This raises the question of whether one can obtain a description of the asymptotic behavior of the complete genealogical structure of the process on large space-time scales, which will in turn allow for asymptotic descriptions of interesting genealogical statistics that are not expressible in a natural way in terms of the configuration process.

These observations on the genealogical structure goes back to the graphical
construction of the voter model due to Harris, and continues up to the
historical process of Dawson and Perkins for branching models~\cite{DP91}, or representation
by contour processes \cite{LeGallLeJan98, Aldous93}.
To better describe genealogies, the notion of {\em $\R$-trees,  marked $\R$-trees or
marked measure $\R$-trees} were developed as a framework \cite{E97, EvaPitWin06}.
These objects contain the relevant information {\em abstracted} from the {\em labeled} genealogy tree, where every
individual is coded with its lifespan and its locations at each time. Such a coding means in particular that all members of the population are distinguished, which information is mostly not needed. In the large population limit, it suffices to consider the statistics of the population via sampling.

For this purpose, one equips the population with a {\em metric} (genealogical distance),
a probability measure (the so-called {\em sampling measure}, which allows to draw
typical finite samples from the population) {{}together with a {\em mark}} (specifying types and locations).
This description in terms of random marked metric measure spaces (in fact, the metric defines a tree)
is the most natural framework to discuss the asymptotic analysis of
population models, since it comprises exactly the information one wants
to keep for the analysis in the limits of populations with even locally infinitely
many individuals. The evolution of a process with such a state space is described by {\em martingale problems}~\cite{GPWmp13,DGP12}.
\bigskip

\noindent
{\bf Open problems and conjectures.} We show in this paper that the spatial continuum limit of the one-dimensional
interacting Fleming-Viot genealogy process solves a martingale problem. However, due to the singular nature of the generator,
the uniqueness of this martingale problem is non-standard, and we leave it for a future paper. {{}In particular, it is difficult to establish the duality relation at the level of generators, because the generators for the process and its dual are only defined for special
test functions {\em at special points}.}

Instead of $\Z$, resp.\ $\R$, as geographic spaces, one could consider the hierarchical group
$\Omega_N = \oplus_\N Z_N$, with $Z_N$ being the cyclical group of order $N$, respectively the continuum
hierarchical group $\Omega^\infty_N = \bigoplus_\Z Z_N$. Brownian motion on $\R$ can be replaced by suitable
L\'evy processes on $\Omega^\infty_N$ and the program of this paper can then be carried out. The Brownian web
would have to be replaced by an object based on L\'evy processes as studied in \cite{EF96},
\cite{DEFKZ00}. We conjecture that the analogues of our theorems would hold in this context.
A further challenge would be to give a unified treatment of these results on $\R, \Omega^\infty_N$.

Another direction is to consider the genealogy processes of interacting Feller diffusions,
catalytic or mutually catalytic diffusions, interacting logistic Feller diffusions, and derive their genealogical continuum limits.
A more difficult extension would be to include {\em ancestral paths} as marks, which raises
new challenges related to topological properties of the state space. {Namely, as the space of paths is a Polish space, it is not a Heine-Borel space
as closed balls around a path are not compact. }
\bi

\noindent
{\bf Outline of Section \ref{S:intro}.}
The remainder of the introduction is organized as follows. In Subsection \ref{S:Intrommm}, we recall and extend
the notion of marked metric measure spaces needed to describe
the genealogies. In Subsection \ref{S:IntroIFV}, we define the interacting Fleming-Viot (IFV) genealogy process via a
martingale problem and give a dual representation in terms of a spatial coalescent.
In Subsection \ref{S:IntroCSSM}, we give, based on the Brownian web,
a graphical construction for the continuum-sites stepping stone model (CSSM) on $\R$ and its marked genealogy process,
which is the continuum analogue and scaling limit of the interacting Fleming-Viot genealogy process on $\Z$,
under diffusive scaling of space and time and normalizing of measure,  a fact which we state in
Subsection \ref{S:Introconv}.
In Subsection~\ref{S:IntroCSSMmart},
we formulate a martingale problem for the CSSM genealogy process.
In Subsection~\ref{S:Introoutline} we  outline the rest of the paper.
\bi

\subsection{Marked Metric Measure Spaces}\label{S:Intrommm}
{{}In this subsection we introduce the state space of the genealogies of interacting Fleming-Viot processes. We want to describe evolving
genealogies of the whole population of all individuals  {currently alive allowing to sample from this population}. We also want to include the individuals' positions in geographic space and possible genetic types. We therefore regard genealogies as (equivalence classes) of
marked metric measure spaces. As our geographic space is infinite (and not compact), it won't be possible to sample individuals by means of a finite (or after renormalizing of a probability) measure. We rather require the sampling measure to be finite on all subpopulations which can be obtained by restricting to finite geographic subspaces.

The following definition of {\em marked metric measure spaces}
extends the one introduced in~\cite{DGP11} which considered probability measures only. }

\bd[$V$-mmm-spaces]\label{D:mmm} Let $(V, r_V)$ be a complete separable metric space with metric $r_V$, and let $o$ be a distinguished point in $V$.
\begin{itemize}
\item[1.] We call $(X, r, \mu)$ a $V$-marked metric measure space ($V$-mmm space for short) if:
\begin{itemize}
\item[(i)] $(X, r)$ is a complete separable metric space,

\item[(ii)] $\mu$ is a measure on the Borel $\sigma$-algebra of $X\times V$, with $\mu(X\times B_o(R))<\infty$ for each ball $B_o(R)\subset V$ of finite radius $R$ centered at $o$.
\end{itemize}

\item[2.] We say two $V$-mmm spaces $(X, r_X, \mu_X)$ and $(Y, r_Y, \mu_Y)$ are equivalent if there exists
a measureable map $\varphi: X\to Y$, such that
\begin{equation}
r_X(x_1, x_2) = r_Y (\varphi(x_1), \varphi(x_2)) \qquad \mbox{for all } \ x_1, x_2 \in {\rm supp}(\mu_X(\cdot \times V)),
\end{equation}
and if $\widetilde \varphi: X\times V\to Y\times V$ is defined by $\widetilde \varphi(x,v) := (\varphi(x), v)$, then
\begin{equation}
\mu_Y = \mu_X \circ \widetilde \varphi^{-1}.
\end{equation}
In other words, $\varphi$ is an isometry between $supp (\mu_X(\cdot\times V))$ and $supp (\mu_Y(\cdot\times V))$, and the induced map $\widetilde\varphi$ is mark and measure preserving.
 We denote the equivalence class of $(X, r, \mu)$ by
 \be{agx1}
 \overline{(X, r, \mu)}.
\ee

\item[3.] The space of (equivalence classes of) $V$-mmm spaces is denoted by
\begin{equation}
\label{e:MV}
\M^V := \big\{ \overline{(X, r, \mu)}: (X, r, \mu) \mbox{ is a $V$-mmm space} \big\}.
\end{equation}
{
\item
[4.] The subspace of (equivalence classes of) $V$-mmm spaces which admit a mark function is denoted by
\begin{equation}
\label{e:MVfct}
\begin{aligned}
&\M^V_{\mathrm{fct}}
\\
&:= \{ \overline{(X, r, \mu)}\in\mathbb{M}^V: \exists\,{\mbox{ measurable }}\kappa:X\to V\mbox{ s.t. }\mu(\mathrm{d}(x,\boldsymbol{\cdot}))=\mu(\mathrm{d}x\times V)\otimes\delta_{\kappa(x)}\}.
\end{aligned}
\end{equation}
}
\end{itemize}
\ed
Note that $\M^V$ depends both on the set $V$ and the metric $r_V$ since the latter defines the sets on which the measure
must be finite.

Marked metric measure spaces were introduced in \cite{DGP11}, which extends the notion of metric measure spaces studied earlier in~\cite{GPW09}.
Definition~\ref{D:mmm} is exactly the analogue of \cite[Def.~2.1]{DGP11}, where $\mu$ is a probability measure.
The basic interpretation in our context is that: $X$ is the space of individuals;
$r(x_1, x_2)$ measures the genealogical distance between two individuals $x_1$ and $x_2$ in $X$;
 the measure $\mu$ is a measure on the individuals and the marks they carry (which can be spatial location as well as type,
or even ancestral paths up to now), allowing us to draw samples from individuals with marks in a bounded set.

\medskip

To define a topology on $\M^V$ that makes it a Polish space, we will make use of the {\em marked Gromov-weak topology} introduced in~\cite[Def.~2.4]{DGP11} for $V$-mmm spaces equipped with probability measures.
{}In this topology a sequence converges iff we can embedd the involved metric spaces isometrically into one metric space such that the images of the ``sampling'' measures converges weakly.
\color{black}
The basic idea {}to extend this to $ \M^V $ \color{black} is that, given our assumption on $\mu$ in Definition~\ref{D:mmm}.1.(ii), we can localize $\mu$ to finite balls in $V$ to reduce $\mu$ to a finite measure. We can then apply the {\em marked Gromov-weak topology} (which also applies to finite measures instead of probability measures) to require convergence for each such localized version of the $V$-mmm spaces. We will call such a topology {\em $V$-marked Gromov-weak$^\#$ topology}, replacing {\em weak} by {\em weak$^\#$}, following the terminology in \cite[Section A2.6]{DV85} for the convergence of measures that are bounded on bounded subsets of a complete separable metric space. {Note that vague convergence is for measures that are finite on compact
rather than bounded subsets. Both notions agree on Heine-Borel spaces (compare, \cite{AthreyaLohrWinter}).}

\bd[$V$-marked Gromov-weak$^\#$ Topology] \label{D:mmmtop} Fix a sequence of continuous functions $\psi_k : V\to [0,1]$, $k\in\N$, such that $\psi_k=1$ on $B_o(k)$, the ball of radius $k$ centered at $o\in V$, and $\psi_k=0$ on $B^c_o(k+1)$. Let $\chi:= \overline{(X, r, \mu)}$ and $\chi_n:= \overline{(X_n, r_n, \mu_n)}$, $n\in\N$, be elements of $\M^V$.  Let $\psi_k\cdot \mu$ be the measure on $X\times V$ defined by $(\psi_k\cdot \mu)(\dd (x, v)) := \psi_k(v) \mu(\dd (x, v))$, and let $\psi_k\cdot \mu_n$ be defined similarly. We say that $\chi_n \asto{n} \chi$ in the $V$-marked Gromov-weak$^\#$ topology if and only if:
\begin{equation}
\overline{(X_n, r_n, \psi_k\cdot \mu_n)} \Asto{n} \overline{(X, r, \psi_k\cdot \mu)} \quad \mbox{in the Gromov-weak topology for each } k\in \N.
\end{equation}

When $V=\R^d$, we may choose $\psi_k$ to be infinitely differentiable.
\ed
\br[Dependence on $o$ and $(\psi_k)_{k\in\N}$]{\rm Note that the $V$-marked Gromov-weak$^\#$ topology does not depend on the choice $o\in V$ and the sequence $(\psi_k)_{k\in\N}$, as long as $\psi_k$ has bounded support and $A_k:=\{v: \psi_k(v)=1\}$ increases to $V$ as $k\to\infty$.
}
\er

\br[$\M^V$ as a subspace of $(\M^V_f)^\N$]\label{R:GromovHash}{\rm Let $\M^V_f$ denote the space of (equivalent classes of) $V$-mmm spaces with finite measures, equipped with the $V$-marked Gromov-weak topology as introduced in~\cite[Def.~2.4]{DGP11}. Then {}it is a well-known fact that \color{black} each element $\overline{(X, r, \mu)}\in \M^V$ can be identified with a sequence $(\overline{(X, r, \psi_1\cdot\mu)}, \overline{(X, r, \psi_2\cdot\mu)}, \ldots )$ in the product space $(\M^V_f)^\N$, equipped with the product topology. This identification allows us to easily deduce many properties of $\M^V$ from properties of $\M^V_f$ that were established in~\cite{DGP11}. In particular, we can metrize the
$V$-marked Gromov-weak$^\#$ topology on $\M^V$ by introducing a metric (which can be called {\em $V$-marked Gromov-Prohorov$^\#$ metric})
\begin{equation}\label{GromovHashMetric}
d_{MGP^\#} (\overline{(X_1, r_1, \mu_1)}, \overline{(X_2, r_2, \mu_2)}) := \sum_{k=1}^\infty \frac{d_{MGP}(\overline{(X_1, r_1, \psi_k\cdot\mu_1)}, \overline{(X_2, r_2, \psi_k\cdot\mu_2)}) \wedge 1}{2^k},
\end{equation}
where $d_{MGP}$ is the marked Gromov-Prohorov metric on $\M^V_f$, which was introduced in~\cite[Def.~3.1]{DGP11} and metrizes the marked Gromov-weak topology.

}
\er

{The proof of the following result is in Appendix~\ref{S:mmm}.}
\bt[Polish space] \label{T:mmmPolish} The space $\M^V$, equipped with the $V$-marked Gromov-weak$^\#$ topology, is a Polish space.
\et


Points in $\M^V$, as well as weak convergence of $\M^V$-valued random variables, can be determined by the so-called {\em polynomials} on $\M^V$, which are defined via sampling a finite subset on the $V$-mmm space.

\bd[Polynomials] \label{D:poly} Let $n\in\N$. Let {$g \in C_{bb}(V^n, \R)$, } the space of real-valued bounded continuous function on $V^n$ with bounded support. For $k\in \N\cup\{0,\infty\}$, let $\phi\in C_b^k(\R_+^{\binom{n}{2}}, \R)$, the space of
$k$-times continuously differentiable real-valued functions on
    $\R_+^{\binom{n}{2}}$ with bounded derivatives up to order $k$. We call the function $\Phi^{n,\phi, g}: \M^V \to \R$ defined by
    \begin{equation} \label{Phingphi}
    \Phi^{n,\phi, g} (\overline{(X, r, \mu)}) := \idotsint \phi(\uur)  g(\uv)  \mu^{\otimes n} ({\rm d}(\ux, \uv)),
    \end{equation}
    a monomial of order $n$, where $\uv:=(v_1, \ldots, v_n)$, $\ux:=(x_1,\ldots, x_n)$, $\uur =\uur(\ux):= (r(x_i, x_j))_{1\leq i<j\leq n}$, and $\mu^{\otimes n} (\dd(\ux, \uv))$ denotes the $n$-fold product measure of $\mu$ on $(X\times V)^n$.
    \begin{itemize}
    \item[\rm (a)] Let $\Pi^k_n:= \{\Phi^{n,\phi, g}: \phi \in  C_b^k(\R_+^{\binom{n}{2}}, \R), g\in C_{bb}(V^n, \R)\}$,
    which we call the space of monomials of order $n$ (with differentiability of order $k$). Let $\Pi^k_0$ be the set of constant functions.
    We then denote by
    $$
    \Pi^k:=\cup_{n\in \N_0} \Pi^k_n
    $$
    the set of all monomials (with differentiability of order $k$).

    \item[\rm (b)] For $V=\R^d$ and $k,l \in \N\cup\{0,\infty\}$,
    we define
    $$
    \Pi^{k, l}_n:=\{ \Phi^{n,\phi, g}: \phi\in C_b^k(\R_+^{\binom{n}{2}}, \R), g \in C^l_{bb}(V^n, \R) \},
    $$
    and  $\Pi^{k, l}:=\cup_{n\in\N_0} \Pi^{k,l}_n$.

 \item[\rm (c)] We call the linear spaces
\be{gs1a}
\wt \Pi^{k,\ell} \mbox{  generated by } \Pi^{k,\ell},
\ee
the polynomials (with differentiability of $\phi$, resp.\ $g$, of order $k$,  resp.\ $\ell$).
    \end{itemize}
\ed

\begin{remark}\label{R:poly}\rm \quad
Note that the polynomials form an algebra of bounded continuous functions   {{}since the product of two monomials can be seen as a new monomial as defined in \eqref{Phingphi}. However, the sum of two monomials in general is not a monomial.}
\end{remark}

Throughout the paper we are often interested in the following sub-space of $\mathbb{M}^V$ (compare, e.g., Definition~\ref{D:Ur} of so-called regular space of states). {Let $b$ be a measurable function on $\R_+$, and write
\begin{equation}
\label{e:Vbeta2}
   \mathbb{M}^{(V,\le b)}:=\big\{\overline{(X,r,\mu)}\in\mathbb{M}^V:\,\mu(X{}\times\color{black} B_r(o))\le b(r)\big\}.
\end{equation}

If $\beta$ is a measure {{}on ${\mathcal B}(V)$} which is finite on bounded subsets, and
\begin{equation}
\label{e:Vbeta}
   \mathbb{M}^{(V,\beta)}:=\big\{\overline{(X,r,\mu)}\in\mathbb{M}^V:\,\mu(X{}\times\color{black}\boldsymbol{\cdot})=\beta\big\},
\end{equation}
then obviously $\mathbb{M}^{(V,\beta)}$ is a closed subspace of $\mathbb{M}^{(V,\le b)}$ with $b(r):={\beta}(B_r(o))$.
 }

\bt[Convergence determining class] \label{T:polydet}
{Fix a measurable function $\beta:\R_+\to\R_+$.  }
We have the following properties for $\Pi^k$, for each $k\in \N\cup\{0,\infty\}$:
\begin{itemize}
\item[\rm (i)] $\Pi^k$ is convergence determining in {$\M^{(V,\le b)}$}. Namely, $\overline{(X_n, r_n, \mu_n)}\to \overline{(X, r, \mu)}$
in {$\M^{(V,\le b)}$} if and only if $\Phi(\overline{(X_n, r_n, \mu_n)})\to \Phi(\overline{(X, r, \mu)})$ as $n\to\infty$ for all $\Phi \in \Pi^k$.

\item[\rm (ii)] $\Pi^k$ is also convergence determining in {$\CM_1(\M^{(V,\le b)})$}, the space of probability measures on {$\M^{(V,\le{}b\color{black})}$}.
Namely, a sequence of {$\M^{(V,\le b)}$}-valued random variables $({\CX}_n)_{n\in \N}$ converges weakly to a {$\M^{(V,\le b)}$}-valued random variable
$\CX$ if and only if $\E[\Phi({\CX}_n)]\to \E[\Phi({\CX})]$ as $n \to \infty$ for all $\Phi \in \Pi^k$.

\item[\rm (iii)] For $V=\R^d$ and each $k,l \in \N\cup\{0,\infty\}$, $\Pi^{k,l}$ is also convergence determining in {$\M^{(V,\le b)}$} and {$\CM_1(\M^{(V,\le b)})$}.
\end{itemize}

We defer the proof of Theorems \ref{T:mmmPolish} and \ref{T:polydet}, as well as some additional properties of $V$-mmm spaces, to Appendix~\ref{S:mmm}.

\et

\begin{remark}\label{R:UV}\rm \quad
For the models we consider, the genealogies lie in {}certain particular \color{black} Polish spaces which arise as {\em closed} subspaces of $\M^V$.
Note that the current population alive corresponds to the {\em leaves} of a genealogical tree, and  the associated $V$-mmm space is {\em ultra\-metric}.
We will denote the space of $V$-marked  ultra\-metric measure spaces by $\U^V$. They form a closed subspace of $\M^V$ and hence
$\U^V$ is Polish. The same holds for {$\M^{(V,\beta)}$, for some Borel measure $\beta$ on $V$ which is finite on bounded subsets.}
\end{remark}

\begin{remark}\label{R:mark}{\rm \quad
{Recall $\mathbb{M}^V_{\mathrm{fct}}$ from Definition~\ref{D:mmm}, and notice that $\mathbb{M}^V_{\mathrm{fct}}$
is not {{}closed}, and that we therefore choose the bigger space $\M^V$ as the state space. The space $\mathbb{M}^V$ allows an individual $x\in X$ to carry a set of marks,
 equipped with the conditional measure of $\mu$ on $V$ given $x\in X$.
 If each individual carries only a single mark
 which we can identify via a {\em mark function} $\kappa : X\to V$, the corresponding marked metric measure space is an element of $\mathbb{M}^V_{\mathrm{fct}}$.  This will be the case for the interacting Fleming-Viot process that we will study.

 It can be shown that every element in $\mathbb{M}^V_{\mathrm{fct}}$ is an element of the closed space
 \begin{equation}
\label{e:MVh}
\begin{aligned}
   \mathbb{M}^V_{h}
   &:=\bigcup_{A\;\mbox{\tiny bdd}}\bigcap_{\delta>0}\big\{(X,r,\mu)\in\mathbb{M}^V:\,
   \\
   &\hspace{1.7cm}\exists\mu'_A\in{\mathcal M}_f(X\times A):\,\mu'\le\mu(\boldsymbol{\cdot}\times A),\|\mu'-\mu(\boldsymbol{\cdot}\times A)\|_{\mathrm{TV}}\le h_A(\delta),
   \\
   &\hspace{1.7cm}
   (\mu'_A)^{\otimes 2}\{r(x_1,x_2)<\delta,d(v_1,v_2)>h_A(\delta)\}=0\big\},
 \end{aligned}
 \end{equation}
 for some $h_A\in{\mathcal H}$, where
 \begin{equation}
 \label{e:calH}
    {\mathcal H}:=\{h:\R_+\to\R_+\cup\{\infty\}:\,h\mbox{ is continuous and increasing, and }h(0)=0\}
 \end{equation}
 (\cite[Lemma~2.8]{LKliem}).}
}
\end{remark}


\subsection{Interacting Fleming-Viot (IFV) Genealogy Processes}\label{S:IntroIFV}

We now study the genealogies of the measure-valued interacting Fleming-Viot (IFV) processes
on a countable geographic space and with allelic types, typically taken from the type space $[0,1]$ (see \cite{DGV95} for details on IFV), which is motivated by the following individual-based model,
the {so-called} Moran model.

Consider a population of individuals{{}, $X$,} with locations indexed by a {\em countable additive group} $V$ (for us this later will be $\Z$). The individuals
{\em migrate} independently according to rate one continuous time random walks with transition probability kernel $a(\cdot, \cdot)$,
\be{gs1}
a(v_1, v_2) = a(0, v_2-v_1) \qquad \mbox{for all }\  v_1, v_2 \in V.
\ee
We denote the transition kernel of the time reversed random walks by
\begin{equation}\label{abar}
\bar a(v_1, v_2) := a(v_2, v_1) = a (0, v_1-v_2).
\end{equation}
Individuals furthermore reproduce by {\em resampling}, where every pair of individuals at the same site
{dies} together at exponential rate $\gamma>0$, and with equal probability, one of the two individuals is chosen to give
birth to two new individuals at the same site with the same type as the parent. This naturally induces a genealogical structure.
The genealogical distance{{}, $r$,} of two individuals at time $t$ is $2\min\{t,T\}$ plus the distance of the ancestors at time 0,
 where $T$ is the time it takes to go back
to the most recent common ancestor. Imposing the {{}Haar measure (here the counting measure), $\mu$,}  on the population with each individual carrying its
location as a {\em mark}, we obtain a $V$-mmm space
and its equivalence class is the state of the genealogy process.

Letting now the number of individuals per site tend to infinity and normalizing the measure such that each site carries
population mass of order one, we obtain a diffusion model, the {\em interacting Fleming-Viot (IFV) genealogy process}  with
{state space
\be{agx2}
   \U^V_1:=\U^{(V,\mathfrak{n})},
\ee
where $\mathfrak{n}$ denotes the {\em {}Haar \color{black} measure} on the countable  geographic space $V$,
the 1 indicating that the measure restricted to each colony is a probability measure.}
This (see Remark \ref{R:UV}) is again a Polish space. (For the diffusion limit in the case of a finite geographic space $V$, see \cite{GPWmp13}, \cite{DGP12}).

\begin{remark}{\rm \quad
If we  introduce as marks (besides {\em locations} from a countable geographic space $\G$) also allelic {\em types} from some set $\K$
(typically taken as a closed subset of $[0,1]$), then the type is inherited at reproduction and $V = \K \times \G$ is the
product of type space  and geographic space. In this case the localization procedure in Definition~\ref{D:mmmtop} applies to the geographic space, since
$\K$ is compact.}
\end{remark}

\subsubsection{The genealogical IFV martingale problem}
\label{sss.genIFVmp}

We now define the interacting Fleming-Viot (IFV) genealogy process via a martingale problem for a linear operator $L^{\mathrm{FV}}$ on $C_b(\U^V_1, \R)$, acting on polynomials.
For simplicity we first leave out allelic types, which we introduce later in Remark \ref{R1.12}.

The linear operator $L^{\mathrm{FV}}$ on $C_b (\U^V_1,\R)$, with domain $\Pi^{1,0}$ as introduced in Definition~\ref{D:poly}\,(b) {}with $d=1$, \color{black} consists of three terms,
corresponding respectively to aging, migration, and reproduction by resampling. With
$\CX =\overline{(X,r,\mu)}$ and $\Phi=\Phi^{n, \phi, g} \in \Pi^{1,0}$,

\bea{gs2}
L^{\mathrm{FV}} \Phi (\CX)
&=& 2 \intl_{(X\times V)^n} \mu^{\otimes n} (\dd(\ux, \uv))\, g(\uv) \suml_{1 \leq k < \ell \leq n}
             \frac{\partial}{\partial r_{k,\ell}}\, \phi(\uur)  \nonumber \\
&&+\  \intl_{(X\times V)^n} \mu^{\otimes n} (\dd(\ux, \uv))\, \phi(\uur) \suml^n_{j=1} \suml_{v^\prime \in V} \bar a(v_j, v^\prime)
             (M_{v_j, v^\prime}g - g) (\uv)  \\
&&+\ 2\, \gamma \intl_{(X \times V)^n} \mu^{\otimes n} (\dd (\ux,\uv))\ g(\uv)
              \suml_{1 \leq k < \ell \leq n} 1_{\{v_k = v_\ell\}} (\theta_{k,\ell} \phi - \phi)(\uur), \nonumber
	      \eea
where  $\ux =(x_1, \cdots, x_n)$, $\uv = (v_1, \cdots, v_n)$, $\uur := (r_{k,\ell})_{1 \leq k < \ell \leq n} := (r(x_k, x_{\ell})_{1\leq k<\ell \leq n})$, and
\be{gs4}
(M_{v_j, v^\prime} g) (v_1, \cdots, v_n) := g(v_1, \cdots, v_{j-1}, v^\prime, v_{j+1}, \cdots, v_n)
\ee
encodes the replacement by migration of the $j$-th sampled individual corresponding to a jump  from location
$v_j$ to $v^\prime$, while $\theta_{k,\ell}$ encodes the replacement of the $\ell$-th individual by the $k$-th individual (both at the same site). {More} precisely,
\be{gs3}
(\theta_{k,\ell} \phi)(\uur) := \phi (\theta_{k,\ell} \uur), \quad \mbox{with} \quad (\theta_{k,\ell} \uur)_{i,j}
:= \left\{ \begin{aligned}
r(x_i, x_j) \qquad &\mbox{ if } i,j \neq \ell,\\
r(x_i, x_k) \qquad &\mbox{ if } j=\ell, \\
r(x_k, x_j) \qquad &\mbox{ if } i=\ell.
\end{aligned}
\right.
\ee

The first result states  that there is a unique $\U^V_1$-valued diffusion process associated with this operator.

\bt[Martingale problem characterization of IVF Genealogy processes]\label{T.wpp}

For any $\CX_0 =\overline{(X_0, r_0, \mu_0)} \in \U^V_1$, we have:
\begin{itemize}
\item[\rm (i)] The $(L^{\mathrm{FV}}, \Pi^{1,0}, \delta_{\CX_0})$-martingale problem is well-posed, i.e. there exists a $\U^V_1$-valued process
 $\CX^{\rm FV}:=(\CX^{\rm FV}_t)_{t\geq 0}$, unique in its distribution, which has initial condition $\CX_0$ and c\`adl\`ag path,
  such that for all $\Phi\in \Pi^{1,0}$ and w.r.t. the natural filtration generated by $(\CX^{\rm FV}_t)_{t\geq 0}$,
\be{gs4b}
\Big(\Phi (\CX^{\mathrm{FV}}_t) - \Phi (\CX^{\mathrm{FV}}_0) - \intl^t_0 (L^{\mathrm{FV}} \Phi) (\CX^{\mathrm{FV}}_s)
 \dd s\Big)_{t \geq 0} \mbox{ is a martingale. }
\ee

\item[\rm (ii)] {{}The solutions (for varying initial conditions) define a strong Markov with continuous path. This Markov process has the Feller property, i.e., the one-dimensional distributions depend continuously on the initial distribution.}

\item[\rm (iii)] {If the initial state admits a mark function, then so does the path almost surely.}
\end{itemize}
\et

\begin{remark}\label{R1.12}{\rm \quad
If we add the type of an individual as an additional mark, i.e., $V := \G \times \K$ with geographic and
type space respectively, then the same result holds if we modify as follows. We require the states to satisfy
the constraint that the projection of $\mu$ on the geographic space $V$ is still the counting measure. The test functions $\Phi$ should be
modified so that we multiply $g: \G^n\to\R$, acting on the locations of the $n$ sampled individuals, by another {bounded, continuous} factor $g_{\rm typ}: \K^n\to \R$, acting on the types of the individuals. The generator $L^{\mathrm{FV}}$ should be modified accordingly, so that $g_{\rm typ}$ changes at resampling from $g_{\mathrm{typ}}$ to
 $g_{\mathrm{typ}} \circ \wt \theta_{k,\ell}$, with $\wt \theta_{k,\ell}$ replacing the $\ell$-th sampled individual by the $k$-th one, see \cite{DGP12}).
}
\end{remark}

\begin{remark}\label{R1.13} \quad {\rm
The process $\CX^{\mathrm{FV}}_t = \overline{(X_t, r_t, \mu_t)}$ has the property that the
measure-valued process $\wh \CX_t$ given by the collection
$\{\mu_t(X_t\times \{i\}\times \cdot), i \in \G\}$, is the IFV process on $(\CM_1 (\K))^\G$.
}
\end{remark}

\begin{remark}\label{R:15}{\rm \quad
From Section~\ref{S:IntroCSSM} onward, we will choose $V=\Z$. However in the subsequent analysis, it is important to observe
that we can embed $\Z$ into $\R$ and view $\U^\Z$ as a closed subspace of $\U^\R$, and view the IFV genealogy process as $\U^\R$-valued process. }
\end{remark}

\subsubsection{Duality}
\label{sss.dual}

The IFV genealogy process $\CX^{\rm FV}$ can be characterized by a dual process,  the {\em spatial coalescent}. The formulation given here can also incorporate mutation and selection (see \cite{DGP12}).

{{}For each $n\in\mathbb{N}$, let $\tilde{S}_n$ denote the space of {\em partitions} of the set $\{1,...,n\}$, i.e., $\pi=\{\mathfrak{p};\,\mathfrak{p}\in\pi\}\in\tilde{\mathbb{S}}_n$ is a collection of disjoint subsets $\mathfrak{p}\subseteq\{1,...,n\}$,
referred to as partition elements (or blocks), such that $\uplus_{\mathfrak{p}\in\pi}\mathfrak{p}=\{1,...,n\}$.
Moreover, let $\tilde{\mathbb{S}}^V_n$ denote the space of {\em marked partitions} of the set $\{1,...,n\}$ with mark space $V$, i.e., we regard $\{(\mathfrak{p},v_{\mathfrak{p}});\,\pi=\{\mathfrak{p};\,\mathfrak{p}\in\pi\}\in\tilde{\mathbb{S}}_n,v_\mathfrak{p}\in V\}\in\tilde{\mathbb{S}}^V_n$ as a partition of $\{1,...,n\}$ where each partition element is assigned a mark in $V$.
Finally, let $\mathbb{S}_n:=\tilde{\mathbb{S}}^V_n\times\R_+^{{n\choose 2}}$
denote the space of {{\em historical marked partitions}} of the set $\{1,...,n\}$ with mark space $V$. That is, we regard $(\{(\mathfrak{p},v_{\mathfrak{p}});\,\mathfrak{p}\in\pi\},(r_{ij})_{1\le i<j\le n})\in\mathbb{S}_n$ as consisting of a marked partition and a matrix of mutual distances.}

{{}The dual process $({\mathcal K}_t)_{t\ge 0}$ will take values in $\mathbb{S}:=\bigcup_{n \in \N}\mathbb{S}_n$ with the following dynamics: given a finite historical marked partition, independently every pair of partition elements with the same mark in $V$ merges at rate $\gamma$.
Until a pair of partition elements merges, the marks migrate independently of each other on $V$
according to a continuous time random walk with  transition kernel $\bar{a}$. After merging the marks of the two involved partition elements will move together forever.
At time $t$, we define the genealogical distance $r_t(i,j)$
of
$i$ and $j$ in $\{1,2,...,n\}$ as $2\min\{t, T_{i,j}\}$, where $T_{i,j}$ is the first time that $i$ and $j$ belong to the same partition element. }

For each $n$, $\phi \in  C_b^\infty(\R_+^{\binom{n}{2}}, \R)$, and $g\in C_{bb}(V^n, \R)$, we define a {\em duality function} $H : \U^V_1 \times \mathbb{S}_n\to\R$ with
{{}\begin{equation}
\label{gs6}
   H\big(\CX, \CK\big)
 :=
   \int_{(X\times V)^{|\pi|}} \mu^{{{}\otimes|\pi|}}(\mathrm{d}(\ux,\uv))\,\mathbf{1}_{\{v_{\mathfrak{p}}=\xi'_{\mathfrak{p}};\,\forall\mathfrak{p}\in\pi\}}\cdot  g(\uv^\pi) \phi \big(\uur^\pi(\underline{x}) + \uur'\big),
\end{equation}
where $\CX=\overline{(X, r, \mu)} \in \U^V_1$, and $\CK = (\{(\mathfrak{p},\xi'_{\mathfrak{p}});\,\mathfrak{p}\in\pi\},\uur')\in {{}\mathbb{S}_n}$, and $\uv^\pi = (v_{\mathfrak{p}(i)})_{i=1, ..., n}$, and $\uur^\pi(\underline{x}):= (r(x_{\mathfrak{p}(i)},x_{\mathfrak{p}(j)}))_{1 \leq i < j \leq n}$, 	
	with $\mathfrak{p}(i)$ being the partition element of $\pi$ containing $i\in\{1,...,n\}$.}

\begin{remark}\label{R:gs6}\rm \quad
	Note that $\{H(\cdot, \CK): \CK \in \BS_n, \phi \in  C_b^\infty(\R_+^{\binom{n}{2}}, \R), g\in C_{bb}(V^n, \R), n\in\N\}$ is law-determining and convergence-determining on $\U^V_1$.
	\end{remark}

The IFV genealogy process $(\CX^{\rm FV}_t)_{t \geq 0}$ is dual to the coalescent $(\CK_t)_{t \geq 0}$, and its law and behavior as $t \to \infty$ can be determined as follows.

\bt[Duality and longtime behaviour]\label{T.dual}
The following properties hold for the IFV genealogy process $(\CX^{\rm FV}_t)_{t\geq 0}$:
\begin{itemize}
\item[\rm (a)] For every $\CX^{\rm FV}_0 \in \U^V_1$ and $\CK_0 \in \BS$, we have
\be{gs7}
\E[H(\CX^{\mathrm{FV}}_t, \CK_0)] = \E[H(\CX^{\mathrm{FV}}_0, \CK_t)], \quad t \geq 0.
\ee

\item[\rm (b)] If $\wh a(\cdot, \cdot) = \frac{1}{2} (a (\cdot, \cdot) + \bar a(\cdot, \cdot))$ is recurrent, then
\be{gs7b}
\CL[\CX^{\mathrm{FV}}_t] \tto \Gamma \in \CM_1 (\U^V_1),
\ee
where $ \Gamma$ is the unique invariant measure of the process $\CX^{\mathrm{FV}}$ on $\U^V_1$.
\end{itemize}
\et

\begin{remark}\label{R:transient} {\rm \quad
If $\wh a$ is transient, then we can decompose $X^{\mathrm{FV}}_t= \overline{(X_t, r_t, \mu_t)}$ in such a way that
$X_t$ is the countable union of disjoint sets $X^i_t$, $\mu_t$ is the sum of measures $\mu^i_t$, and
$\{(X^i_t, r_t |_{X^i_t \times X^i_t}, \mu^i_t)\}_{i \in\N}$ are $V$-mmm spaces such that for $x \in X^i_t, x^\prime \in X^j_t$
with $i \neq j$, $r_t(x,x^\prime)$  diverges in probability as $t \to \infty$, and each
$(X^i_t, r_t |_{X^i_t \times X^i_t}, \mu^i_t)$ converges in law to a limiting $V$-mmm space.
Alternatively we can transform distances $r: r \to 1- e^{-r}$ and obtain a unique equilibrium
in that case. See also \cite{gmuk}.

}
\end{remark}


\begin{remark}[Strong duality] \rm As the interacting Fleming-Viot {{}process} is a population model, its evolving genealogy $(U_t,r_t,\mu_t)_{t\ge 0}$ can be represented by a $V$-marked $\R$-tree $(X,r)$. {{}That is, we can find a $0$-hyperbolic metric space (with distances $\infty$ possible equal to) $(X,r)$ such that any two $x,y\in X$ of finite distance are connected by a path,} and isometries $(\varphi_t)_{t\ge 0}$  with $\varphi_t:\,U_t\to X$ such that for all $x\in X$ there is a $t\ge 0$ with $\varphi^{-1}_t(x)\in\mathrm{supp}(\mu_t(\boldsymbol{\cdot}\times V))$.
We won't write down this $\R$-tree representation explicitly here but it can be derived from the look-down construction presented in \cite{GLW05} in a straightforward way.

We would like to point out that this representation by an $\R$-tree allows to define the ancestor $A^T_{s}$ of $x$ at time $T$ back at time $T-s$, and the above duality relation can actually be stated in a strong sense.
That is, we can construct
for each $T>0$ our model together with a dual process $\tilde{K}^T$ on the same probability space such that
for all $X'\subseteq U_T$ with $\# X'<\infty$,  $\phi \in  C_b^\infty(\R_+^{\binom{\# X'}{2}}, \R)$, and $g\in C_{bb}(V^{\# X'}, \R)$, and
$H=H^{\# X',\varphi,g}$ from (\ref{gs6}),
for all $s\in[0,T]$,
\begin{equation}
\label{e:strong}
   H\big({\mathcal U}_{T-s},\tilde{\mathcal K}^T_{s}\big)\equiv H\big({\mathcal U}_{T},\tilde{\mathcal K}^T_{0}\big),\;\mbox{ almost surely.}
\end{equation}
Notice that whenever a duality relation holds almost surely rather than in expectation, one refers to it as a {\em strong form of duality}.

Indeed, define  for fixed $T\ge 0$ and any finite subset $X':=\{x_1,...,x_{\# X'}\}\subset U_T$ the map $\tilde{K}^{T}:=(\tilde{K}^{T}_s:=((\{(\mathfrak{p},\xi'_{\mathfrak{p}});\,\mathfrak{p}\in\pi_s^T\},r^{',T}_s))_{s\in[0,T]}$ which sends ${\mathcal U}:=({\mathcal U}_t:=(U_t,r_t,\mu_t))_{t\ge 0}$
to a path with values in the space of historical marked partitions $\mathbb{S}_{\# X'}$
defined by
\begin{equation}
\label{rr}
   r^{',T}_{s}(i,j)=r\big(\varphi_T(x_i),\varphi_T(x_j)\big)\wedge(T-s),
\end{equation}
for all $1\le i<j\le \# X'$,
a marked partition of $\{1,2,...,n\}$ defined through the equivalence relation
\begin{equation}
\label{pp}
   i\equiv^T_{s} j\hspace{.2cm}\mbox{ iff }\hspace{.2cm}r^{',T}_{s}(i,j)<(T-s),
\end{equation}
and a family of position functions on $\varphi_T(\mathrm{supp}(\mu_T(\boldsymbol{\cdot}\times V)))$ such that
\begin{equation}
\label{xixi}
  \xi^{',T}_{s}(\mathfrak{p}):=\kappa_{T-s}\big(\varphi_{T-s}^{-1}(\mathfrak{p})\big).
\end{equation}

By construction, $\tilde{K}^{T}$ is the dual spatial coalescent. Moreover, 
for all $s\in[0,T]$, and $(\underline{x},\underline{v})\in (X_{T}\times V)^{|\pi_0|}$,
\begin{equation}
\label{e:distt}
   \uur^{\pi_s}_{T-s}(A^T_s(\underline{x})) + {\uur'}^T_s=\uur_T^{\pi_0}(\underline{x}) + {\uur'}^T_0,
\end{equation}
and s$\kappa_{T-s}(A_s^T(x))$ has under $\mu_T$ the same distribution as types under $\mu_{T-s}$.
This together implies the strong form of duality as stated in (\ref{e:strong}).
\end{remark}\qed

As a further consequence of relation (\ref{gs7}) we can specify the {\em finite dimensional distributions} of $(\CX^{\mathrm{FV}}_t)_{t \geq 0}$
completely in {\em terms of the spatial coalescent} as follows. Fix a time horizon {{}$T >0$}.
The finite dimensional distributions are determined by the expectation of all test functions $\wt \Phi: C([0,\infty),\mathbb{M}^V)$ of the following form:
{{}\be{gs8}
  \wt{\Phi}\big((\CX^{\mathrm{FV}}_t)_{t \geq 0}\big):=\Pi^{\ell}_{k=1}\Phi_k\big(\CX^{\mathrm{FV}}_{t_k}\big),
\ee
for some {{}$\ell\in \N$}, $0 \leq t_1 < t_2 < ... < {t_\ell} = t$ and $\Phi_k=\Phi^{n_k,g_k,\phi_k} \in \Pi^{1,0}$ of order $n_{k}\in\mathbb{N}$,  $\phi_k\in  C_b^\infty(\R_+^{\binom{n_k}{2}}, \R)$, and $g_k\in C_{bb}(V^{n_k}, \R)$,
for each $k ={1, ..., \ell}$.
}

The dual is the {\em spatial coalescent with frozen particles $(\wt \CK_s)_{s \in [0,{{}T}]}$ },
for which the time index $s\in [0,{{}T}]$ runs in the opposite direction from {{}the time index of} $\CX^{\rm FV}$.
Namely, looking backward from time {{}$T$}, for each {$1\le k\le\ell$}, we start $n_k$ particles
at time {{}$T-t_k$} at locations $\xi^1_{t_k}:=(\xi^1_{t_k,1}, ..., \xi^1_{t_k, n_k})\in V^{{}n_k}$, each forming its
own partition element in the partition
{{}$\pi\in\tilde{\mathbb{S}_n}$, $n:=n_1+...+n_{\ell}$}. The particles perform the usual dynamics of the
spatial coalescent with the restriction that all particles starting at time {{}$T-t_k$} were kept frozen before time {{}$T-t_k$}.
{{}At time $s$, the genealogical distance $r_s(i,j)$ between two individuals $i$ and $j$,
started respectively at times {{}$T-t_i$ and $T-t_j$}, is defined to be
\begin{equation}
\label{disfrozen}
   2\min\{s, T_{i,j}\}-\min\{s,(T-t_i)\}-\min\{s,(T-t_j)\},
\end{equation}
}
where $T_{i,j}$ is the first time the two individuals coalesce.

{{}
Denote this new spatial coalescent process with frozen particles by $(\wt \CK_t)_{t \geq 0}$. The state space $\BS$ is once more the space of historical marked partitions.
We then define the duality function $H:\, \U^{V}_1 \times \BS \to \R$ which determines the finite
dimensional distributions of {$\CX$} for varying $\wt \CK$ as in (\ref{gs6}) but now
with
\be{gs8a1}
 \phi (\uur) := \prod_{k=1}^{\ell}  \phi_k \big((r_{i,j})_{n_1+\cdots+n_{k-1}< i <j \leq n_1+\cdots+n_k}\big) {}\; \mbox{and} \; g\big( \uv^\pi\big):=\prod^{\ell}_{k=1} g_k\big(t_k,\uv^{\pi_k}\big),
\ee {}where $ \pi_k $ are the partition elements started at time $ t_k $. \color{black}
}

{
The following space-time duality is an immediate consequence of the Markov property and the duality applied successively to the
time intervals $[t_{\ell-1}, t_\ell]$, $[t_{\ell-2}, t_{\ell-1}]$, ..., $[0,t_1]$.
}\bi

\bcor[Space-time duality] \label{C.tsdual}
Let $\CX^{\mathrm{FV}}_0 \in \U^V_1$, and $\wt\CK_0 \in \mathbb{S}$ be as in (\ref{gs8a2}) with $t_k\leq T$ for all {$1\le k\le\ell$}. Let $\wt \Phi$ be defined as above (\ref{disfrozen}). Then
\be{gs9}
\E\big[\wt \Phi \big((\CX^{\mathrm{FV}}_s)_{s \in [0,T]}\big)\big] = \E \big[H\big(\CX^{\mathrm{FV}}_0, \wt \CK_T\big)\big].
\ee
\end{corollary}
In words, the genealogical distances among the ${}n=\color{black}{n_1+\cdots +n_\ell}$ individuals sampled from $(\CX^{\mathrm{FV}}_s)_{s \in [0,T]}$,
with $n_k$ individuals sampled at time $t_k$ at specified locations, can be recovered by letting these individuals evolve backward
in time as a spatial coalescent until time $0$, at which point we sample from $\CX^{\rm FV}_0$ according to the location of each
 partition element in the spatial coalescent.


{
}

\subsection{Genealogies of Continuum-sites Stepping Stone Model (CSSM) on $\R$}
\label{S:IntroCSSM}
{
If we rescale space and time diffusively,
the measure-valued interacting Fleming-Viot process on $\Z$ converges to a continuum space limit, the so-called
{\em continuum-sites stepping stone model} (CSSM) .
Formally, CSSM is a measure-valued process $\nu:=(\nu_t)_{t\geq 0}$ on $\R\times [0,1]$, where $\R$ is the geographical space and $[0,1]$ is the type space. We might think of individuals in the population which undergo independent Brownian motions, and whenever two individuals meet, one of the two individuals is chosen with equal probability and switches its type to that of the second individual. Provided that $\nu_0(\boldsymbol{\cdot} \times [0,1])$ is the Lebesgue measure on $\R$, CSSM was
rigorously constructed in~\cite{EF96, E97, DEFKZ00, Z03} via a moment duality with coalescing Brownian motions. In particular, $\nu_t(\boldsymbol{\cdot} \times [0,1])$ is the Lebesgue measure on $\R$, for all $t>0$.

In this subsection we will construct the evolving genealogies of the CSSM based on duality to the {\em (dual) Brownian web}, and establish properties (Proposition~\ref{P:reg}, Theorem~\ref{T:MarkovCont}).}

{
For the voter model on $\Z$, the joint genealogy lines of individuals at all space-time points in $\Z\times [0,\infty)$ are
distributed as a collection of coalescing random walks evolving backward in time (see~\cite{Lig85}).
Analogously, for the CSSM on $\R$, when $n$ individuals are sampled from the population at possibly different times,
their joint genealogy lines evolve backward in time as coalescing Brownian motions. Upon reaching time $0$, each surviving genealogy
line then independently selects an ancestral type by sampling according to the conditional distribution of $\nu_0$ on the type space $[0,1]$,
conditioned on the spatial location of the genealogy line. Furthermore the joint genealogy lines of individuals at all space-time points in $\R\times [0,\infty)$
are distributed as a collection of coalescing Brownian motions evolving backward in time.
Therefore the CSSM on $\R$ is exactly the continuum analogue of the interacting Fleming-Viot process (as well as the voter model) on $\Z$.

Although having an uncountable number (starting from every space-time point in $\R\times [0,\infty)$) of coalescing Brownian motions seems problematic, this object has been constructed rigorously and is now known as the
{\em (dual) Brownian web} $\wh \CW$~\cite{FINR04, FINR06}, dual to a {\em forward Brownian web} $\CW$ constructed on the {\em same} probability space.
The (dual) Brownian web is essentially a collection of coalescing Brownian motions on $\R$, starting from every point in the space-time plane
$\R\times\R$.
}

In~\cite{FINR04}, the Brownian web $\CW$ is constructed as a random variable where each realization of
\be{gs9b}
\CW \mbox{ is a closed subset of } \Pi:=\cup_{s\in\R} C([s,\infty),\R),
\ee
the space of continuous paths in $\R$ with any starting time $s\in \R$. {{}In \cite{FINR04}, the topology is defined by first compactifying the space $ \R^2 $ suitably, and then by choosing for}
$\Pi$
the {\em topology of local uniform convergence} {{}and requiring the initial times to be close if paths are close. } For each $z:=(x,t)\in\R$, we will let $\CW(z):=\CW(x,t)$ denote the subset of paths in $\CW$
with starting position $x$ and starting time $t$. We can construct $\CW$ by first fixing a countable dense subset $\CD\subset \R^2$,
and then construct a collection of coalescing Brownian motions $\{\CW(z): z\in \CD\}$, with one Brownian motion starting from each $z\in \CD$.
The Brownian web $\CW$ is then obtained by taking a suitable closure of $\{\CW(z): z\in \CD\}$ in $\Pi$, which gives rise to a set of
paths starting from every point in the space-time plane $\R^2$. It can be shown that the law of $\CW$ does {\em not} depend on the choice of $\CD$  (see Theorem B.1 in Appendix~\ref{S:web}).

The Brownian web $\CW$ has a {\em graphical} dual called the {\em dual Brownian web}, which we denote by $\wh\CW$.
Formally,
\be{gs9c}
\wh\CW \mbox{ is a random closed subset of } \wh \Pi:=\cup_{s\in\R} C((-\infty, s],\R),
\ee
 the space of continuous paths in
 $\R$ starting at any time $s\in\R$ and running backward in time as coalescing Brownian motions.

The {\em joint} distribution of $\wh \CW$ and $\CW$ is uniquely determined by the requirement that the paths of $\wh \CW$ never cross paths of $\CW$
(see, Theorem B.3 in Appendix~\ref{S:web}). Thus, jointly, the Brownian web and its dual is a random variable taking values in a Polish space, with
 \be{gs9d}
 (\CW, \wh \CW) \in \Pi \times \wh \Pi.
 \ee

Interpreting coalescing Brownian motions in the (dual) Brownian web as {\em ancestral lines} specifying the genealogies, we can then give an
{\em almost sure graphical construction of the CSSM}, instead of relying on moment duality relations
as in~\cite{EF96, E97, DEFKZ00, Z03}, which nevertheless we get as corollary of the graphical construction.
The classical measure-valued CSSM process can be recovered from $(\CX^{\rm CS}_t)_{t \geq 0}$ by projecting the sampling measure $(\mu^{\rm CS}_t)_{t\geq 0}$
 to the mark space $V$, if $V$ is chosen to be the product of geographical space $\R$ and type space $[0,1]$. In what follows, we will take $V$ to be only
  the geographical space $\R$, since types have no influence on the evolution of genealogies.

We next explicitly construct a version of the CSSM genealogy process
\be{ag1}
\CX^{\rm CS}:=(\CX^{\rm CS}_t)_{t\geq 0}, \quad
\CX^{\rm CS}_t=\overline{(X^{\rm CS}_t, r^{\rm CS}_t, \mu^{\rm CS}_t)}, \qquad t\geq 0,
\ee
as a functional of the (dual) Brownian web $(\CW, \wh \CW)$.

To avoid a disruption of the flow of presentation, background details on the (dual) Brownian web we will need are collected in Appendix~\ref{S:web}.
\medskip

We proceed in three steps:\sm

\noindent
{
{\bf Step~1 (Initial states).} A}ssume that $\CX^{\rm CS}_0$ belongs to the following
closed subspace of $\U^\R$:
\begin{equation}
\label{UV1}
\U^\R_1 := \{ \overline{(X, r, \mu)}\in \U^\R: \mu(X\times\boldsymbol{\cdot}) \mbox{ is the Lebesgue measure on } \R\}.
\end{equation}
In other words, $\U^\R_1$ is the set of $\R$-marked ultra\-metric measure spaces where the projection of the measure on the mark space
(geographic space) $\R$
is the Lebesgue measure. This is necessary for the duality between CSSM and coalescing Brownian motions.
We will see that almost surely $\CX^{\rm CS}_t\in \U^\R_1$ for all $t\geq 0$.
\medskip

\noindent
{
{\bf Step~2 (The time-$t$ genealogy as a metric space).} To define $(X^{\rm CS}_t,r^{\rm CS}_t)$ for every $t>0$, let us fix a realization of $(\CW,\wh\CW)$,
 ({{}see the Appendix~\ref{S:web} for more details}).
 For each $t>0$, let
\begin{equation}\label{At}
   A_t := \big\{ v\in\R: \wh\CW(v,t) \mbox{ contains a single path } \hat f_{(v,t)}\big\}, \quad E_t:=\R\backslash A_t.
\end{equation}
By Lemma~\ref{L:spts} on the classification of points in $\R^2$ w.r.t.\ $\CW$ and $\wh\CW$, almost surely, $E_t$
is a countable set for each $t>0$. For each $v\in A_t$, we interpret $\hat f_{(v,t)}$ as the genealogy line of the individual
at the space-time coordinate $(v,t)$. Genealogy lines of individuals at different space-time coordinates evolve backward in time
 and coalesce with each other. At time $0$, each genealogy line traces back to exactly one spatial location in the set
\begin{equation}
\wh E:= \{\hat f_{(v,s)}(0): v\in\R, s>0, \hat f\in \wh\CW \},
\end{equation}
where we note that $\wh E$ is almost surely a countable set, because by Theorem~\ref{T:webchar} and Lemma~\ref{L:pathconv},
 paths in $\wh\CW$ can be approximated in a strong sense by a countable subset of paths in $\wh\CW$. For each $v\in \wh E$,
  we then identify a common ancestor $\xi(v)$ for all the individuals whose genealogy lines trace back to spatial location $v$ at time $0$,
   by sampling an individual $\xi(v) \in X^{\rm CS}_0$ according to the conditional distribution of $\mu^{\rm CS}_0$ on $X^{\rm CS}_0$,
   conditioned on the spatial mark in the product space $X^{\rm CS}_0\times \R$ being equal to $v$.

We next characterize individuals by points in space. Note that there is a natural
genealogical distance between points in $A_t$. For individuals $x,y\in A_t$, if $\hat f_{(x,t)}$ and $\hat f_{(y,t)}$ coalesce at time
   $\hat\tau<t$, then denoting $u:=\hat f_{(x,t)}(0)$ and $v:=\hat f_{(y,t)}(0)$, we define the distance between $x$ and $y$ by
\begin{equation}\label{rtxy}
r_t(x,y) :=
\left\{
\begin{aligned}
2 (t-\hat \tau) \qquad & \qquad \mbox{ if } \hat\tau\geq 0, \\
2 t + r^{\rm CS}_0(\xi(u), \xi(v))  & \qquad \mbox{ if } \hat\tau< 0.
\end{aligned}
\right.
\end{equation}

First define $(X^{\rm CS}_t,r^{\rm CS}_t)$ as the {\em closure} of $A_t$ w.r.t.\ the metric $r_t$ defined in (\ref{rtxy}). Note that $(X^{\rm CS}_t, r^{\rm CS}_t)$ is ultra\-metric, and by construction Polish.

\br{\rm We may even extend the distance $r_t$ to a distance $r$ between $(x,t)$ with $x\in A_t$ and $t>0$, and $(y,s)$ with $y\in A_s$ and $s>0$. More precisely, let
\begin{equation}\label{rxtys}
r((x,t), (y,s)) :=
\left\{
\begin{aligned}
s+t-2\hat \tau  \qquad & \qquad \mbox{ if } \hat\tau\geq 0, \\
s+ t + r^{\rm CS}_0(\xi(u), \xi(v))  & \qquad \mbox{ if } \hat\tau< 0.
\end{aligned}
\right.
\end{equation}
}
\er
\medskip

\noindent
{\bf Step 3 (Adding the sampling measure).}
We now define $\CX^{\rm CS}_t:= \overline{(X^{\rm CS}_t, r^{\rm CS}_t, \mu^{\rm CS}_t)}$.
For that purpose, we will represent {{}next} $X^{CS}_t$ as an {\em enriched} copy of $\mathbb{R}$ (see (\ref{XCSsplit}) below).

By identifying each $x\in A_t$ with the path $\hat f_{(x,t)}\in \wh\CW$, we can also identify $X^{\rm CS}_t$ with the closure of $\{\hat f_{(x,t)}\in \wh\CW: x\in A_t\}$ in $\wh\Pi$, because a sequence $x_n\in A_t$ is a Cauchy sequence w.r.t.\ the metric $r_t$ if and only if the sequence of paths $\hat f_{(x_n, t)}$ is a Cauchy sequence when the distance between two paths is measured by the time to coalescence, which by Lemma~\ref{L:pathconv}, is also equivalent to $(\hat f_{(x_n, t)})_{n\in\N}$ being a Cauchy sequence in $\wh\Pi$.

When we take the closure of $\{\hat f_{(x,t)}\in \wh\CW: x\in A_t\}$ in $\wh\Pi$, only a countable number of paths in $\wh\CW$ are added, which are precisely the leftmost and rightmost paths in $\wh\CW(x,t)$, when $\wh\CW(x,t)$ contains more than one path.

Namely for each $x\in E_t$ (recall from (\ref{At})), let $x^+, x^-$ denote the two copies of $x$ obtained by taking limits of $x_n\in A_t$ with either $x_n\downarrow x$ or $x_n\uparrow x$ in $\R$, and let $E^\pm_t:=\{ x^\pm: x\in E_t\}$. We can then take
\begin{equation}\label{XCSsplit}
X^{\rm CS}_t := A_t \cup E^+_t \cup E^-_t,
\end{equation}
equipped with a metric $r^{\rm CS}_t$, which is the extension of $r_t$ from $A_t$ to its closure $X^{\rm CS}_t$, giving a Polish space
$(\CX^{\mathrm{CS}}_t, r^{\mathrm{CS}}_t)$.

Next to get a {\em sampling measure}, note that each finite ball in $(X^{\rm CS}_t, r^{\rm CS}_t)$ with radius less than $t$
can be identified with an interval in $\R$ (modulo a subset of $E_t\cup E_t^+\cup E_t^-$), and hence can be assigned the
Lebesgue measure of this interval. {{}That is, we define the Borel measure $ \wt\mu_t^{CS}$
 on $ (X_t^{CS}, r_t^{CS})$ by
 \begin{equation}
 \label{e:LebegueBorel}
    \wt{\mu}_t^{CS}\big(\{x'\in X_t^{\mathrm{CS}}:\,r_t^{\mathrm{CS}}(x',x)<\delta\}\big):=R(x,\delta)-L(x,\delta),\;\;\delta<t,x\in X_t^{\mathrm{CS}}
 \end{equation}
 if $\{x'\in X_t^{\mathrm{CS}}:\,r_t^{\mathrm{CS}}(x',x)<\delta\}=(L(x,\delta),R(x,\delta))\subseteq\mathbb{R}$. We then define the sampling measure
 $\mu_t^{CS}$ on ${\mathcal B}((X^{\rm CS}_t,r^{\rm CS}_t)\times(\mathbb{R},d_{\mathrm{eucl}}))$ by
\begin{equation}\label{e.addsampme}
   \mu_t^{CS}(\mathrm{d}x\mathrm{d}v):=\wt\mu_t^{CS}(\mathrm{d}x)\delta_x(dv).
\end{equation}}
This completes our construction of the CSSM genealogy process $(\CX^{\rm CS}_t)_{t\geq 0}$. \medskip

{{}
\begin{remark}[Notational simplification] In the sequel we will apply the existence of a mark function and the embedding of the basic set in the enriched reals to simplify notation. Namely, we will write integrals w.r.t. the sampling measure over $X\times V $
as an integral w.r.t. the Lebesgue measure over $\R$.
\end{remark}
}

A key feature is again {\em duality}.
We can replace in the Definition of the process
$(\CK_t)_{t \geq 0}$ from Subsubsection \ref{sss.dual} the random walks by Brownian motions
to obtain $(\CK^{\mathrm{Br}}_t)_{t \geq 0}$, which is dual to the process
$(\CX^{\mathrm{CS}}_t)_{t \geq 0}$ by construction:

\beC{C.dual}{($H$-duality)}
The tree-valued CS-process and marked tree-valued coalescing Brownian motions are in $H$-duality, i.e., for each $H$ of the form (\ref{gs6}),
\be{ag25}
\E [H(\CX^{\mathrm{CS}}_t, \CK^{\mathrm{Br}}_0)] = \E[H(\CX^{\mathrm{CS}}_0, \CK^{\mathrm{Br}}_t)].
\ee
Furthermore by the above construction {strong duality  holds.} \qad
\ecor

\begin{remark}\label{R:ag25}\quad {\rm
Note in the continuum case the function $\Phi^{n,g,\varphi} (\cdot, \kappa^{\mathrm{Br}})$
is not a polynomial since we fix the locations we consider. In order to get a polynomial
we need to consider a function $g$ on mark space with $g \in C^2_{bb} (\R,\R)$
over which we integrate w.r.t. the sampling measure.}
\end{remark}

\medskip

We collect below some basic properties for the CSSM genealogy processes that we just constructed.

\begin{proposition}[Regularity of states]
Let $\CX^{\rm CS}:=(\CX^{\rm CS}_t)_{t\geq 0}$ be the CSSM genealogy process constructed from the dual Brownian web $\wh\CW$,
with $\CX^{\rm CS}_0\in \U_1^\R$. Then almost surely, for every $t>0$:
\begin{itemize}
\item[\rm (a)] There exists a continuous (mark) function $\kappa_t: X^{\rm CS}_t \to \R$, i.e., $\mu^{\rm CS}_t( {\rm d}x {\rm d}v) = \mu^{\rm CS}_t({\rm d}x \times \R) \delta_{{\kappa_t}(x)}({\rm d}v)$;

\item[\rm (b)] For each ${{}\ell}\in (0,t)$, $X^{\rm CS}_t$ is the disjoint union of balls $(B^{{{}\ell}}_i)_{i\in\Z}$ of radius ${{}\ell}$.
Furthermore, there exists a locally finite set $E^{{{}\ell}}_t:=\{v_i\}_{i\in\Z}\subset \R$ with $v_{i-1}<v_i$ for all $i\in\Z$,
such that $\{\kappa(x): x\in B_i^{{{}\ell}}\} = [v_{i-1}, v_i]$ and $\mu^{\rm CS}_t(B^{{{}\ell}}_i) = v_i-v_{i-1}$. We can identify $E^{{}\ell}_t$ from $\wh\CW$ by
    \begin{equation}\label{Elt}
        E^{{{}\ell}}_t:=\{x\in E_t: \hat f^+_{(x,t)} \mbox{ and } \hat f^-_{(x,t)} \mbox{ coalesce at some time } s\leq t-{{}\ell} \},
    \end{equation}
    where $f^+_{(x,t)}$ and $f^-_{(x,t)}$ are respectively the rightmost and leftmost path in $\wh\CW(x,t)$;

\item[\rm (c)] $\CX^{\rm CS}_t=\overline{(X^{\rm CS}_t, r^{\rm CS}_t, \mu^{\rm CS}_t)} \in \U^\R_1$, i.e., $\mu^{\rm CS}_t(X^{\rm CS}_t \times \cdot)$ is the Lebesgue measure on $\R$.
\end{itemize}
\label{P:reg}
\end{proposition}

\begin{remark} {\rm Using the duality between the Brownian web $\CW$ and $\wh\CW$, as characterized in Appendix~\ref{S:web}, it is easily seen that we can also write
\begin{equation}\label{Elt2}
E^{l}_t =\{f(t): f\in \CW(x, t-s) \mbox{ for some } x\in \R, s\geq l\}.
\end{equation}
}
\end{remark}

\bt[Markov property, path continuity, Feller property] \label{T:MarkovCont}

Let $\CX^{\rm CS}$ be as in Proposition~\ref{P:reg}. Then
\begin{itemize}
\item[\rm (a)] $(\CX^{\rm CS}_t)_{t\geq 0}$ is a $\U^\R_1$-valued Markov process;

\item[\rm (b)] Almost surely, $\CX^{\rm CS}_t$ is continuous in $t\geq 0$;

\item[\rm (c)] For each $m\in\N$, let $\CX^{{\rm CS}, (m)}$ be a CSSM genealogy process with $\CX^{{\rm CS}, (m)}_0\in \U^\R_1$. If $\CX^{{\rm CS}, (m)}_0 \to \CX^{{\rm CS}}_0$ in $\U^\R_1$, then for any sequence $t_m\to t\geq 0$, we have {{}weak convergence in law, i.e.\ }$\CX^{{\rm CS}, (m)}_{t_m} \Rightarrow \CX^{{\rm CS}}_t${{}, as $m\to\infty$}.
\item[\rm (d)] For each initial state in $\U^\R_1, \CL [\CX^{\mathrm{CS}}_t]$ converges {{}weakly} to the unique equilibrium distribution on $\U^\R_1$ as $t \to \infty$.
\end{itemize}
\et

{{}
\begin{remark} Notice that this  equilibrium state can be represented in terms of a functional of the Brownian web.
\label{rem:rem}
\end{remark}
}

\begin{remark}\label{R:reg}{\rm \quad
We note that if we allow types as well, we enlarge the mark space from $\R$ to $\R\times [0,1]$, where each individual carries a type in $[0,1]$
that is inherited upon resampling.
Theorem~\ref{T:MarkovCont} still holds in this case. We will consider such an extended mark space in Theorem~\ref{T:conv2} below.
}
\end{remark}

\begin{remark}{\rm \quad
Proposition~\ref{P:reg} shows that even though $\CX^{\rm CS}_0$ can be any state in $\U^\R_1$, for $t>0$, $\CX^{\rm CS}_t$ can only take on a small subset of the state space $\U^\R_1$. This introduces complications in establishing the continuity of the process at $t=0$, and it will also be an important point when we discuss the generator of the associated martingale problem.
}
\end{remark}

\subsection{Convergence of Rescaled IFV Genealogies}\label{S:Introconv}
In this subsection, we establish the convergence of the interacting Fleming-Viot genealogy processes on $\Z$ to that of the CSSM,
where we view the states as elements of $\U^\R$ (see Remark \ref{R:15}).
We assume that the transition probability kernel $a(\cdot, \cdot)$ in the definition of the IFV process satisfies
\be{acond}
\sum_{x\in\Z} a(0,x)x=0 \qquad  \mbox{and}  \qquad \sigma^2:=\sum_{x\in\Z} a(0,x) x^2 = \sum_{x\in \Z} \bar a(0,x) x^2 \in (0,\infty).
\ee

For each $\eps>0$, we then define a scaling map ${S_\eps=S^{\sigma}_\eps}: \U^\R\to \U^\R$ (depending on $\sigma$) as follows. Let ${\CX} =\overline{(X, r, \mu)} \in \U^\R$. Then
\be{Seps}
S_\eps {\CX} := \overline{(X, S_\eps r, S_\eps\mu)},
\ee
where $(S_\eps r)(x,y) := \eps^2 r(x,y)$ for all $x, y \in X$, and $S_\eps \mu$ is the measure on $X\times \R$ induced by $\mu$ and the map $(x,v)\in X\times \R \to (x, \eps \sigma^{-1} v)$, and then
the mass rescaled by a factor of $\eps\sigma^{-1}$. More precisely,
\be{Seps2}
(S_\eps\mu)(F) = \eps\sigma^{-1} \mu\{(x, \eps^{-1}\sigma v): (x,v)\in F \} \qquad \mbox{for all measurable}\ F\subset X\times \R.
\ee

We have the following convergence result for rescaled IFV genealogy processes.

\bt[Convergence of Rescaled IFV Genealogies] \label{T:conv}

Let ${\CX}^{{\rm FV}, \eps} :=({\CX}^{{\rm FV}, \eps}_t)_{t\geq 0}$ be an IFV genealogy process on $\Z$ with initial condition
${\CX}^{{\rm FV}, \eps}_0 \in \U^\Z_1$, indexed by $\eps>0$. Assume that $a(\cdot, \cdot)$ satisfies (\ref{acond}),
and $S_\eps {\CX}^{{\rm FV}, \eps}_0 \to {\CX}^{\rm CS}_0$ for some ${\CX}^{\rm CS}_0 \in \U^\R_1$ as $\eps\to 0$.

{}
Then $(S_\eps{\CX}^{{\rm FV},\eps}_{\eps^{-2}t})_{t\geq 0}$ converges as $C([0,\infty), \U^\R)$-valued random variable weakly to the CSSM genealogy process ${\CX}^{\rm CS} := ({\CX}^{\rm CS}_t)_{t\geq 0}$ as $\eps\to 0$.
\color{black}
\et

To prove Theorem~\ref{T:conv}, we will need an auxiliary result of interest in its own on the convergence of rescaled
measure-valued IFV processes to the measure-valued CSSM. The IFV process with mark space $\R\times [0,1]$ is a measure-valued process $(\widehat \CX^{\rm FV}_t)_{t\geq 0}$,
where $\widehat \CX^{\rm FV}_t$ is a measure on $\R\times [0,1]$, and its projection on $\R$ is the counting measure on $\Z$.
Similarly, the CSSM with mark space $\R\times [0,1]$ is a measure-valued process $(\widehat \CX^{\rm CS}_t)_{t\geq 0}$,
where $\widehat \CX^{\rm CS}_t$ is a measure on $\R\times [0,1]$, and its projection on $\R$ is the Lebesgue measure on $\R$.
Define $\widehat\CX^{\rm FV}$ and $\widehat\CX^{\rm CS}$ respectively as the projection of the measure component of $\CX^{\rm FV}$
and $\CX^{\rm CS}$, projected onto the mark space $\R\times [0,1]$.  We then have

\bt[Convergence of Rescaled IFV Processes] \label{T:conv2}

Let ${\widehat\CX}^{{\rm FV}, \eps} :=({\widehat\CX}^{{\rm FV}, \eps}_t)_{t\geq 0}$ be a measure-valued IFV process on $\Z$, indexed by $\eps>0$.
Assume that $a(\cdot, \cdot)$ satisfies (\ref{acond}), and $S_\eps {\widehat\CX}^{{\rm FV}, \eps}_0$ converges to
${\widehat\CX}^{\rm CS}_0$ w.r.t.\ the vague topology for some ${\widehat\CX}^{\rm CS}_0$ as $\eps\to 0$.
{}
Then $(S_\eps{\widehat \CX}^{{\rm FV}, \eps}_{\eps^{-2}t})_{t\geq 0}$ converges as $C([0,\infty), \CM(\R\times [0,1]))$-valued random variable weakly to the CSSM process ${\widehat\CX}^{\rm CS} := ({\widehat\CX}^{\rm CS}_t)_{t\geq 0}$ as $\eps\to 0$.
\et \color{black}
A similar convergence result has previously been established for the voter model in~\cite{AS11}.

\begin{remark}{\rm \quad As can be seen from the above convergence results and the regularity properties of the limit process in Proposition~\ref{P:reg}, on a macroscopic scale, there are only locally finitely many individuals with descendants surviving for a macroscopic time of $\delta$ or more.
This phenomenon leads in the continuum limit to a dynamic driven by a thin subset of hotspots only. For similar effects in other population models, see for example \cite{DEF02a, DEF02b}. }
\end{remark}

\subsection{Martingale Problem for CSSM Genealogy Processes}\label{S:IntroCSSMmart}
\label{ss.mpgen}

In this section, we show that the CSSM genealogy processes solves a martingale problem with a singular generator.  To identify the generator $L^{\rm CS}$ for the CSSM genealogy process $(\CX^{\rm CS}_t)_{t\geq 0}$, we note that for all $t>0$, $\CX^{\rm CS}_t$ satisfies the regularity properties established in Proposition~\ref{P:reg}. We will see that $L^{\rm CS}$ is only well-defined on $\Phi\in \Pi^{1,2}$ evaluated at points $\CX\in \U^\R$ with suitable regularity properties for $\Phi$ and $\CX$.

We now formalize the subset of regular points in $\U_1^\R$ as follows, which satisfies exactly the properties in Proposition~\ref{P:reg}.

\bd[Regular class of states $\U^\R_{\rm r}$]\label{D:Ur} Let $\U^\R_{\rm r}$ denote the set of $\CX= \overline{(X,r, \mu)}\in\U^\R$
{{}which satisfies} the following regularity properties:
\begin{itemize}
\item[\rm (a)] $\CX \in \U^\R_1$, i.e., $\mu(X \times \cdot)$ is the Lebesgue measure on $\R$;

\item[\rm (b)] there exists a mark function $\kappa: X \to \R$, with $\mu( {\rm d}x {\rm d}v) = \mu({\rm d}x \times \R) \delta_{\kappa(x)}({\rm d}v)$;

\item[\rm (c)] there exists $\delta>0$ such that for each $l\in (0,\delta)$, $X$ is the disjoint union of balls $(B^{l}_i)_{i\in\Z}$ of radius $l$. Furthermore, there exists a locally finite set $E^{l}:=\{v_i\}_{i\in\Z}\subset \R$ with $v_{i-1}<v_i$ for all $i\in\Z$, such that $\{\kappa(x): x\in B_i^{l}\} = [v_{i-1}, v_i]$ and $\mu(B^{l}_i) = v_i-v_{i-1}$.
\end{itemize}
\ed
By Proposition~\ref{P:reg}, $\U^\R_{\rm r}$ is closed {}under the dynamics of $\CX^{\rm CS}$ (i.e., $\CX^{\rm CS}_0\in \U^\R_{\rm r}$ implies that $\CX^{\rm CS}_t\in \U^\R_{\rm r}$ for all $t>0$), separable, metric measurable subset of the Polish space $\U^\R_1$.
 However, it is {\em not complete}.
Note that the first requirement gives rise to a closed set,  the second requirement is known to generate a measurable set~\cite{LKliem}, and {it is not hard to see} that the third requirement also generates a measurable set.


\br\label{R:Ur}{\rm Similar to the discussion leading to (\ref{XCSsplit}), for $\CX \in \U^\R_{\rm r}$, we can give a representation on an enriched copy of $\R$ as follows. Property (c) in Definition~\ref{D:Ur} implies that any two disjoint balls in $X$ are mapped by $\kappa$ to two intervals, which overlap at at most a single point in $E^l$ for some $l>0$. Therefore $\kappa^{-1}(x)$ must contain a single point for all $x\in\R$ with $x$ not in
\begin{equation}
E:=\cup_{l>0} E^l,
\end{equation}
and $\kappa^{-1}(x)$ containing two or more points implies that $x$ is in $\kappa(B_1)\cap \kappa(B_2)$ for two disjoint balls in $X$.
By the same reasoning, for each $x\in E$, $\kappa^{-1}(x)$ must contain exactly two points, which we denote by $x^{\pm}$, where $x^+$ is a limit point of $\{\kappa^{-1}(w): w>x\}$ and $x^-$ is a limit point of
$\{\kappa^{-1}(w): w<x\}$. Similar to (\ref{XCSsplit}), we can then identify $X$ with $(\R\backslash E) \cup E^+ \cup E^-$, where $E^{\pm}:=\{x^{\pm}: x\in E\}$. With this identification {}and with \eqref{e.addsampme}\color{black}, we can simplify our notation (with a slight abuse) and let $\mu$ be the measure on $(\R\backslash E) \cup E^+ \cup E^-$, which assigns no mass to $E^\pm$ and is equal to the Lebesgue measure on $\R\backslash E$.}

\er

We now introduce a regular subset of $\Pi^{1,2}$ needed to define the generator $L^{\rm CS}$.

\bd[Regular class of functions $\Pi^{1,2}_{\rm r}$] \label{D:PI12r}
 Let $\Pi^{1,2}_{\rm r}$ denote the set of regular test functions $\Phi^{n,\phi, g}\in \Pi^{1,2}$, defined as in (\ref{Phingphi}),
  with the property that:
\be{gs13}
\exists\, \delta >0  \quad \mbox{s.t.}\ \forall \, i\neq j \in \{1,2,\cdots,n\}, \quad
\frac{\partial \phi}{\partial r_{i,j}} ((r_{k,l})_{1\leq k<l\leq n}) =0 \quad \forall \, r_{i,j} \in [0,\delta].
\ee
\ed

We can now specify the action of the generator $L^{\rm CS}$ on regular $\Phi$ evaluated at regular points, namely $L \Phi(\CX)$ exists for
$\Phi=\Phi^{n,\phi,g}\in \Pi^{1,2}_{\rm r}$
and $\CX=\overline{(X,r, \mu)}\in \U^\R_{\rm r}$. By Remark \ref{R:Ur}, we can identify $X$ with
$(\R\backslash E) \cup E^+ \cup E^-$, while $\mu$ is identified with the Lebesgue measure on $\R\backslash E$.
 The generator $L^{\rm CS}$ is given by

\begin{equation}\label{LCS}
L^{\rm CS} \Phi(\CX) := L^{\rm CS}_{\rm d} \Phi(\CX) + L^{\rm CS}_{\rm a} \Phi(\CX) + L^{\rm CS}_{\rm r} \Phi(\CX),
\end{equation}
with the component for the massflow (migration) of the population on $\R$ given by
\begin{equation}\label{LCSd}
L^{\rm CS}_{\rm d} \Phi(\CX) := \frac{1}{2} \int_{X^n} \phi(\uur) \Delta g(\ux) \, {\rm d}\ux,
\end{equation}
and the component for aging of individuals given by
\begin{equation}\label{LCSa}
L^{\rm CS}_{\rm a} \Phi(\CX) := 2 \int_{X^n} g(\ux) \sum_{1\leq i<j\leq n} \frac{\partial \phi}{\partial r_{i,j}} (\uur)\ \dd\ux.
\end{equation}
These operators are linear operators on the space of bounded continuous functions, $C_b(\U^\R_1, \R)$,
with domain $\Pi^{1,2}$, and maps polynomials to polynomials of the same order.

The  component for resampling is, with $\theta_{k,l}\phi$ defined as in (\ref{gs3}), given by
\begin{equation}\label{LCSr}
L^{\rm CS}_{\rm r} \Phi(\CX) :=  \sum_{1\leq k\neq l\leq n} \int_{X^n}  1_{\{\kappa(x_k)=\kappa(x_l)\}}\
g(\ux)\ (\theta_{k,l}\phi-\phi)(\uur) \ \mu^*(\dd x_k) \mu^*(\dd x_l) \!\!\prod_{1\leq i\leq n \atop i\neq k,l} \!\!\dd x_i,
\end{equation}
with effective resampling measure and mark functions

\begin{equation}\label{mu*}
\mu^*:= \sum_{x\in E} \delta_{x^+}+ \sum_{x\in E}\delta_{x^-},
\end{equation}
\be{kap1}
\kappa(x^{\pm})=x \mbox{  for  } x\in E,
\ee
where $E$, and $x^\pm$ for $x\in E$, are defined as in Remark~\ref{R:Ur}.

\br{\rm {
Note that $L^{\rm CS}_{\rm r}$ is {\em singular}. First because the effective
resampling measure $\mu^*$ is supported on a countable subset of $X$ and is singular w.r.t.\ the sampling measure $\mu$ on $X$. Secondly, $\mu^*$ is locally infinite because $E\cap (a,b)$ contains infinitely many points for any $a<b$.
}
Therefore the r.h.s.\ {}of \color{black} \eqref{LCSr} is now in principle a sum of countably many
monomials of order $n-2$.

Indeed, as we partition $X=(\R\backslash E) \cup E^+\cup E^-$ into balls of radius $l$,
with $l\downarrow 0$, the balls must correspond to smaller and smaller intervals on $\R$ so as not {}to \color{black} contradict the fact that each point
in $X$ is assigned one value in $\R$. Nevertheless, $L^{\rm CS}_{\rm r}\Phi(\CX)$ in (\ref{LCSr}) is well-defined at least on $\U^\R_{\rm r}$ because
by our assumption that $\Phi \in \Pi^{1,2}_{\rm r}$ and condition (\ref{gs13}) on $\phi$, we have $\theta_{k,l}\phi=\phi$
if the resampling is carried out between two individuals $x^+$ and $x^-$ for some $x\in E$, with $r(x^+, x^-)\leq \delta$.
Thus only resampling involving $x\in E$ with $r(x^+, x^-)>\delta$ remains in the integration w.r.t.\ $\mu^*$,
and such $x$ are contained in the locally finite set $E^{\delta}$ introduced in Definition~\ref{D:Ur}~(c).
Together with the assumption that $g$ has bounded support, this implies that the integral in (\ref{LCSr}) is finite.
}
\er

\begin{remark} {\rm \quad
The operator $L^{\rm CS}_{\rm r}$, defined on functions in $\Pi^{1,2}_{\rm r}$ evaluated at
regular states in $\U^\R_{\rm r} \subseteq \U^\R_1$, is still a linear operator, mapping polynomials to
generalized polynomials of degree reduced by two and with domain $\Pi^{1,2}_{\rm r}$.
Here, generalized polynomial means that they are no longer bounded and continuous, and are only measurable functions
defined on the subset of points $\U^\R_{\rm r} \subseteq \U^\R_1$. Hence $L^{\mathrm{CS}}_r$ differs significantly from
generators associated with Feller semigroups on Polish spaces.}
\end{remark}

For $L^{\rm CS}_{\rm r}\Phi(\CX)$ to be well-defined for any $\Phi \in \Pi^{1,2}$, instead of $\Phi\in \Pi^{1,2}_{\rm r}$,
we need to place further regularity assumption on the point $\CX$ at which we evaluate $\Phi$. These assumptions are satisfied
by typical realizations of the CSSM at a fixed time, as we shall see below.

\bd[Regular subclass of states $\U^\R_{\rm rr}$]\label{D:Urr} Let $\U^\R_{\rm rr}$ be the set of
 $\CX= \overline{(X,r, \mu)}\in\U^\R_{\rm r}$ with the further property that
\begin{equation}
\sum_{x\in E \cap [-n,n]} r(x^+, x^-) <\infty \qquad \mbox{for all } n\in\N,
\end{equation}
where we have identified $X$ with $(\R\backslash E) \cup E^+ \cup E^-$ as in Remark \ref{R:Ur}.
\ed

\bt[Martingale problem for CSSM genealogy processes]\label{T.mp} \qquad
\begin{itemize}
\item[\rm(i)] If $\CX_0 \in \U^\R_1$, then the $(L^{\rm CS}, \Pi^{1,2}, \delta_{\CX_0})$-martingale problem
has a solution, i.e., there exists a process $\CX:=(\CX_t)_{t\geq 0}$ with initial {}state \color{black} $\CX_0$, {}with \color{black}almost surely
continuous sample path {}and \color{black}with $\CX_t \in \U^\R_{\rm r}$ for every $t>0$, such that for all $\Phi\in \Pi^{1,2}$, and w.r.t.\
the natural filtration,
\begin{equation}
\Big(\Phi (\CX_t) - \Phi (\CX_0) - \int_0^t (L^{\rm CS} \Phi) (\CX_s)ds\Big)_{t \geq 0} \mbox{ is a martingale}.
\end{equation}
\item[\rm(ii)]
The CSSM genealogy process $\CX^{\rm CS}$ constructed in Sec.~\ref{S:IntroCSSM} is a solution to the above martingale problem. Apart from the properties established in Prop.~\ref{P:reg}, for each $t>0$, almost surely $\CX^{\rm CS}_t\in \U^\R_{\rm rr}$.
\end{itemize}
\et

\begin{remark}{\rm \quad
Whether almost surely, $\CX^{\rm CS}_t\in \U^\R_{\rm rr}$ for all $t>0$, remains open.
}
\end{remark}

\begin{remark}{\rm \quad
Note that different from the usual martingale or local martingale problems, as for example in \cite{EK86}, the test functions here are only defined on a
(topologically not closed, only dynamically closed) subset of the state space.
}
\end{remark}

We conjecture that the martingale problems above are in fact well-posed. A proof could be attempted by using the duality between the CSSM genealogy process and the Brownian web. There are however subtle technical complications due to the fact that the generator of the martingale problem is highly singular. We leave this for a future paper.

\subsection{Outline}\label{S:Introoutline}
We provide here an outline of the rest of the paper. In Section~\ref{S:IFV} we prove the results on the IFV genealogy processes.
In Section~3 we construct the CSSM genealogy process, and establish
in Section~\ref{S:conv} the convergence of the IFV genealogies to those of the CSSM, and in Section \ref{S:CSSMmart} results on the martingale problem for the CSSM genealogy processes. In Appendix~\ref{S:mmm}, we collect further facts and proofs concerning marked metric measure spaces.
In Appendix~\ref{S:web}, we recall the construction of the Brownian web and its dual, and collect some basic properties of the Brownian web and coalescing Brownian motions.
In Appendix~\ref{A:Corr} we prove some results on coalescing Brownian motions needed in our estimates to derive the
martingale problem for the CSSM.

\section{Proof of Theorems~\ref{T.wpp} and~\ref{T.dual}}

\label{S:IFV}
In this section we present the proof of the results on the evolving genealogies for the interacting Fleming-Viot diffusions.
This model is a special case of  evolving genealogies for the interacting $\Lambda$-Fleming-Viot diffusions which are studied in \cite{gklimwLambda}.

\begin{proof}[Proof of Theorem~\ref{T.wpp}]
We will proceed
in five steps: (1) We show the result on the martingale problem and the duality
for finite geographic spaces $V$.
(2) To prepare for  the general case where $V$ is countable,  we define an approximation procedure with
specific finite geographic space dynamics. (3) We then show the convergence in path space,
as the finite spaces approach $V$.
(4) We verify the
claimed properties of the solution for general $V$ by a direct argument based on the duality and an
explicit look-down construction. 
(5) Finally, we show that the process admits a mark function.

We will use several known facts on measure-valued Fleming-Viot diffusions. For that we refer to \cite[Chapter~5]{D93} {}for \color{black} the non-spatial case and to \cite{DGV95}
for the spatial case.
\sm

\noindent{\bf Step~1 ($V$ finite) } The case where $V$ is finite is very similar to the non-spatial case. We therefore just have to
modify the arguments of the proof of Theorem~1 in \cite{GPWmp13} or Theorem~1 in \cite{DGP12}.

As usual we will conclude {\em uniqueness} of the solution of the martingale problem from a duality relation. Note that in contrast to the (non-spatial) interacting Fleming-Viot model with mutation considered in \cite{DGP12}, in our spatial
interacting Fleming-Viot model, resampling takes place only locally, that is, for individuals at the {\em same} site.
The tree-valued dual will therefore be based on the spatial coalescent considered in \cite{GLW05}.
{As the calculations to verify the duality relation are the same as in \cite{GPWmp13} and in \cite{DGP12}, we omit them here. }

As for {\em existence} of a solution of the martingale problem we consider the martingale problems for the evolving genealogies of the approximating
spatial Moran models. By consistency of the spatial coalescent, we get the uniform convergence of generators for free. Thus we only have to show
the compact containment condition. Here we can rely on the general criterion for population dynamics given in Proposition~2.22 in \cite{GPWmp13}.
As $V$ is finite, all arguments given in \cite{GPWmp13} to verify this criterion simply go through here as well.\sm

\noindent{\bf Step~2 (A coupled family of approximating finite systems) } Let now $V$ be countable, and consider
a sequence $(V_n)_{n \in \N}$ of finite sets with $V_n \subseteq V$, and
 $V_n \uparrow V$.  Put for each $n\in\mathbb{N}$, and for all $v_1,v_2\in V$,
 \be{gs10}
 a_n\big(v_1,v_2\big):=\left\{\begin{array}{cc} \tfrac{a(v_1,v_2)1_{V_n \times V_n} (v_1,v_2)}{\sum_{v_3 \in V_n} a(v_1,v_3)}, & \mbox{ if }v_1\in V_n, \\
 \delta(v_1,v_2), &  \mbox{ if }v_1\not\in V_n.\end{array}\right.
 \ee
Denote then by ${\mathcal X}^{\mathrm{FV},V_n}$ a solution of the martingale problem associated with the operator (restricted to $V_n$)
\begin{equation}
\label{gs2Vn}
\begin{aligned}
   L^{\mathrm{FV},V_n}\Phi\big(\CX\big)
 &= 2\int_{(X\times V)^n} \mu^{\otimes n} (\dd(\ux, \uv))\, g(\uv) \sum_{1 \leq k < \ell \leq n}
             \frac{\partial}{\partial r_{k,\ell}}\, \phi(\uur) \\
&+\  \int_{(X\times V)^n} \mu^{\otimes n} (\dd(\ux, \uv))\, \phi(\uur) \suml^n_{j=1} \sum_{v^\prime \in V} \bar{a}_n(v_j, v^\prime)
             (M_{v_j, v^\prime}g - g) (\uv)  \\
&+\ 2\, \gamma \int_{(X \times V)^n} \mu^{\otimes n} (\dd (\ux,\uv))\ g(\uv)
              \sum_{1 \leq k < \ell \leq n} 1_{\{v_k = v_\ell\}}\big(\theta_{k,\ell} \phi - \phi\big)(\uur),
\end{aligned}
\end{equation}
and ${\mathcal K}^{V_n}$ the spatial coalescent on $V_n$ with  migration rate  $\bar{a}_n(\boldsymbol{\cdot},\boldsymbol{\cdot})$ rather than $\bar{a}(\boldsymbol{\cdot},\boldsymbol{\cdot})$.

Notice that $a_n$ is not necessarily  double stochastic anymore, which turns the duality
with the spatial coalescent into a {\em Feynman-Kac duality}
where the Feynman-Kac term converges to $1$, as $n \to \infty$, on every finite time horizon.
The Feynman-Kac duality reads as follows (compare, for example, \cite[Proposition~3.11]{Seidel14}):  for all ${\mathcal X}_0\in\mathbb{U}_1^V$,
\be{ag6}
   \mathbb{E}\big[H(\CX^{\mathrm{FV},V_n}_t, \CK^{V_n}_0)\big]= \mathbb{E}\big[H(\CX_0, \CK^{V_n}_t)\exp\big(\int^t_0 A(\CK^{V_n}_s)\mathrm{d}s\big)\big],
\ee
where 
\be{ag7}
   A\big((\pi,\xi',\underline{\underline{r}})\big):=\sum_{i\in\pi}\big(\sum_{v^\prime \in V} a_n(v^\prime,\xi'_i)-1\big).
\ee
which is bounded along the path by $|\pi| \cdot$ Const, for all $t \geq 0$.

Establishing the  Feynman-Kac duality requires to check that (compare, Section 4.4 in \cite[Theorem~4.4.11]{EK86}):
\begin{equation}
\label{ag9}
   L^{\mathrm{FV}}H\big(\boldsymbol{\cdot},\CK)\big)(\CX)
   =
   L^{\mathrm{dual}}H\big(\CX,\boldsymbol{\cdot}\big)(\CK)
   + A(\CK)\cdot H\big(\CX, \CK\big).
\end{equation}
This can be immediately verified by explicit calculation (compare, \cite[Section~4]{GPWmp13} for the generator calculation for the resampling part, and
\cite[Proposition~3.11]{Seidel14} for  the generator calculation for Markov chains - here migration - whose transition matrix is not double stochastic).

As for given $n\in\mathbb{N}$, our dynamics consists of independent components outside $V_n$, we can apply Step~1, here with the Feynman-Kac duality, to conclude the well-posedness of the martingale problem with respect to $L^{\mathrm{FV},V_n}$.
\sm

\noindent{\bf Step~3 ($V$ countable) } Fix ${\mathcal X}_0\in\mathbb{U}_1^V$.
In this step we want to show that the family $\{{\mathcal X}^{\mathrm{FV},V_n};\,n\in\mathbb{N}\}$ is tight, and that every limit satisfies the
$(L^{\mathrm{FV}}, \Pi^{1,2}, \delta_{\CX_0})$-martingale problem.

Observe first that the family of laws of the projection of the measures on mark space are tight since the
localized state w.r.t. a fixed finite subset A of $V$ has only finitely many marks and weight $|A|$ (uniformly in $n$). We therefore will here ignore the marks and show tightness in Gromov-weak$^{\#}$-topology.
For that we  want to apply \cite[Corollary~4.5.2]{EK86}.
Since $a_n(v_1,v_2)\to a(v_1,v_2)$ for all $v_1,v_2\in V$, we clearly have that  $L^{\mathrm{FV},V_n} \Phi$ converges {uniformly}
to $L^{\mathrm{FV}} \Phi$, as $V_n\uparrow V$, i.e.,
 $\sup_{{\mathcal X} \in \U^{V}} | L^{\mathrm{FV},V_n} \Phi ({\mathcal X}) - L^{\mathrm{FV}} \Phi({\mathcal X}) |
 \ntoo 0$, for $\Phi$ depending only on finitely many sites,  it remains to verify the {\em compact containment condition}, i.e., to show that
 for every $T>0$, and $\ve >0$ we can find
 a compact set $K_{T,\ve}\subseteq\U^V$ such that for all  $n \in \N$,
\begin{equation}
\label{e:comp}
   P\big(\CX^{\mathrm{FV},V_n}_t\in K_{T,\ve};\,\mbox{ for all }t \in [0,T]\big) \geq 1-\ve.
\end{equation}

For that purpose we will once more rely on the criterion for the compact containment condition which was developed in \cite[Proposition~2.22]{GPWmp13} for
population dynamics. To see first that the criterion applies, notice that the evolving genealogies of interacting Fleming-Viot diffusions can be read off
as a functional of the  {\em look-down construction} given in \cite{GLW05}. Thus the countable representation of the look-down defines  a population dynamics.
In particular, for each $t\ge 0$, we can read off  a representative $(X_t,r_t,\mu_t)$ of ${\mathcal X}^{\mathrm{FV}}_t$ in the look-down graph such that {\em ancestor-descendant} relationship is well-defined. Denote for all $t\ge 0$, $s\in[0,t]$ and $x\in X_t$ by $A_t(x,s)\in X_{t-s}$ the ancestor of $x\in X_t$ back at time $s$, and for  ${\mathcal J}\subseteq X_s$
by $D_t({\mathcal J},s)\subseteq X_t$ the set of descendants of a point in ${\mathcal J}$ at time $t$.

In the following we refer for each finite $A\subseteq V$ to
\begin{equation}
\label{e:restrict}
   {\mathcal X}^{\mathrm{FV},\boldsymbol{\cdot},A}
\end{equation}
as the restriction of  ${\mathcal X}^{\mathrm{FV},\boldsymbol{\cdot}}$
to marks in $A$, i.e., obtained by considering the sampling measure $\mu^A(dxdv):=1_A\mu(dxdv)$. Fix $T>0$.
We then have two show that the following properties are true for all $A\subseteq V$,
\begin{itemize}
\item {\bf Tightness of number of ancestors. }For all $t\in[0,T]$ and $\varepsilon\in(0,t)$, the family $\{S^{V_n}_{2\varepsilon}(X_t,r_t,\mu_t);\,n\in\mathbb{N}\}$ is tight, where $S^{V_n}_{2\varepsilon}(X_t,r_t,\mu_t)$ denotes the minimal number of balls of radius $2\varepsilon$ needed to cover $X_t$ up to a set of $\mu_t$-measure $\varepsilon$.
\item {\bf Bad sets can be controlled. } For all $\varepsilon\in(0,T)$, there exists a {{}$\delta=\delta(\varepsilon)>0$} such that for all $s\in[0,T)$, $n\in\mathbb{N}$ and $\sigma({\mathcal X}^{\mathrm{FV},V_n}_u;\,u\in[0,s])$-measurable random subsets ${\mathcal J}^{V_n}\subseteq X_s\times A$ with $\mu_s({\mathcal J}^{V_n})\le\delta$,
\begin{equation}
\label{e:badsets}
   \limsup_{n\in\mathbb{N}}P\big(\sup_{t\in[s,T]}\mu_t\big(D_t^{V_n}({\mathcal J}^{V_n},s)\times A\big)\ge\epsilon\big)\le\epsilon.
\end{equation}
\end{itemize}

(i) Fix $t\in[0,T]$. W.l.o.g.\ we assume that $V_n$ is a subgroup of $V$ with addition $+_n$ for each $n\in\mathbb{N}$. We consider for each $n\in\mathbb{N}$ another spatial coalescent $\tilde{K}^{V_n}$ on $V_n$ with migration kernel $\bar{\tilde{a}}_n(\boldsymbol{\cdot},\boldsymbol{\cdot})$ rather than $\bar{a}_n(\boldsymbol{\cdot},\boldsymbol{\cdot})$ where
\begin{equation}
\label{e:tildean}
   \tilde{a}_n(v,v'):=\sum_{y\in V;y\sim_n v'}a(v,y),
\end{equation}
where $\sim_n$ denotes equivalence modulo $+_n$. Further, denote by $\tilde{{\mathcal X}}^{\mathrm{FV},V_n}$ the evolving genealogies of the interacting Fleming-Viot diffusions
whose migration kernel is   $\tilde{a}_n(\boldsymbol{\cdot},\boldsymbol{\cdot})$ rather than ${a}_n(\boldsymbol{\cdot},\boldsymbol{\cdot})$. We prefer to work with $\tilde{K}^{V_n}$.
For this spatial coalescent it was verified in the proof of Proposition~3.4 in \cite{GLW05} that for any time $t>0$ the total number of partition elements of $\tilde{K}^{V_n}_t$ which are located in $A$ is stochastically bounded  uniformly in $n\in\mathbb{N}$.
As the kernel $\tilde{a}_n(\boldsymbol{\cdot},\boldsymbol{\cdot})$ is double stochastic, $\tilde{{\mathcal X}}^{\mathrm{FV},V_n}$ and $\tilde{K}^{V_n}$ are dual (without a Feynman-Kac potential), and as $d_{\mathrm{GP^{\#}}}(\tilde{{\mathcal X}}^{\mathrm{FV},V_n},{\mathcal X}^{\mathrm{FV},V_n})\to 0$, as $n\to\infty$,
the claim follows. \sm

(ii) Fix $T>0$, $A\subseteq V$ finite, $\varepsilon\in(0,T)$, $s\in[0,T)$, $n\in\mathbb{N}$ and a $\sigma({\mathcal X}^{\mathrm{FV},V_n}_u;\,u\in[0,s])$-measurable random subset ${\mathcal J}^{V_n}\subseteq X_s\times A$.
From the generator characterization of $\CX^{\mathrm{FV},V_n}$, we can conclude that the process $\{\mu_t(D_t({\mathcal J}^{V_n},s));\,t\ge s\}$ is a
$V$-indexed system of interacting (measure-valued) Fisher-Wright diffusions.
We have to find $\delta=\delta(\varepsilon)$ such that (\ref{e:badsets}) holds if $\mu_s({\mathcal J}^{V_n})\le\delta$.

Notice that $\{\mu_t(D_t^{V_n}({\mathcal J}^{V_n},s)\times A);\,t\ge s\}$ is a {\em semi-martingale} and given by a martingale with continuous
paths due to resampling plus a {\em deterministic flow} in and out of the set $A$ due to migration. Therefore we have to control the
fluctuation of the martingale part and the maximal flow out of the set $A$ over a time interval of
length $t-s$.

The martingale part is estimated from below with Doob's maximum inequality (the quadratic variation is bounded uniformly
in the state and in $n$ by a constant $\cdot |A|$, details are left to the reader).
The deterministic out flow occurs
at most at a finite rate $c\times |A|$ since the total mass of every site is one.
This estimate is uniform in the parameter $n$ (recall the random walk kernel is perturbed
by restricting it to $V_n$). Similarly the flow into $A$ can be bounded by $c\cdot |A|$ independently  of $n$
but the flow out of the set $A$ occurs with a maximal rate $c$ independently of $n$ with
$n \geq n_0(A)$, and
\begin{equation}
\label{ag23}
   c=2\cdot\max\big\{\max_{v' \in A^c,v \in A} a(v', v),\max_{v' \in A^c,v \in A} a(v', v)\big\}.
\end{equation}\sm

\noindent {\bf Step~4 (Feller property). }
We next {{}prove} that ${\mathcal X}^{\mathrm{FV}}$ has the Feller property. From here it is standard to conclude that  ${\mathcal X}^{\mathrm{FV}}$ satisfies the strong Markov property
(see, for example, \cite[Theorem~4.2.7]{EK86}). Consider a sequence $({\mathcal X}^{(n)}_0)_{n\in\mathbb{N}}$ in $\mathbb{U}^V_1$ such that ${\mathcal X}^{(n)}\ntoo{\mathcal X}_0$, Gromov-weak$^{\#}$ly, for some ${\mathcal X}_0\in \mathbb{U}^V_1$. Denote by
${\mathcal X}^{\mathrm{FV},{\mathcal X}^{(n)}_0}$ and ${\mathcal X}^{\mathrm{FV},{\mathcal X}_0}$ the evolving genealogies of the interacting Fleming-Viot diffusions started in ${\mathcal X}^{(n)}_0$ and $\CX_0$, respectively, and let
${\mathcal K}$ be our tree-valued dual spatial coalescent, and $H$ as in (\ref{gs6}). Then for each given $t\ge 0$,
\begin{equation}
\label{e:-1proof}
\begin{aligned}
      \mathbb{E}\big[H\big(\CX^{(n)}_0, \CK_t\big)\big|\CK_t\big]
      &\ntoo
      \mathbb{E}\big[H\big(\CX_0, \CK_t\big)\big|\CK_t\big],\hspace{1cm}a.s.
\end{aligned}
\end{equation}

Thus, by our duality relation,
\begin{equation}
\label{e:1proof}
\begin{aligned}
      \mathbb{E}\big[H(\CX^{\mathrm{FV},{\mathcal X}^{(n)}_0}_t, \CK_0)\big]
      &= \mathbb{E}\big[H(\CX^{(n)}_0, \CK_t\big)\big]
      \\
      &\ntoo
      \mathbb{E}\big[H\big(\CX_0, \CK_t\big)\big]
      = \mathbb{E}\big[H\big(\CX^{\mathrm{FV},{\mathcal X}_0}_t, \CK_0)\big].
\end{aligned}
\end{equation}
Recall from Remark~{R:gs6} that the family $\{H^{n,\phi}(\boldsymbol{\cdot},{\mathcal K});\,n\in\mathbb{N},\phi\in C_b(\R_+^{n\choose 2});\CK\in\mathbb{S}_n\}$ is convergence determining.
Thus it follows that $\CX^{\mathrm{FV},{\mathcal X}^{(n)}_0}_t\nto\CX^{\mathrm{FV},{\mathcal X}_0}_t$, for all $t\ge 0$.\sm

{

\noindent{\bf Step~5 (Mark function) } Fix $T\ge 0$, and $\CX_0\in\mathbb{U}^V_{\mathrm{fct}}$. \smallskip

For the proof we will rely once more on the approximation of the solution of the $(L^{\mathrm{FV}}, \Pi^{1,0}, \delta_{\CX_0})$-martingale problem  by $\mathbb{U}^V_{\mathrm{fct}}$-valued evolving genealogies of Moran models, ${\mathcal X}^{\mathrm{M},\rho}$, where $\rho>0$ is the local intensity of individuals.
By the look-down construction given in \cite{GLW05}, we can define the family $\{{\mathcal X}^{\mathrm{M},\rho};\,\rho>0\}$ on one and the same probability space.
Moreover, as the solution of the  $(L^{\mathrm{FV}}, \Pi^{1,0}, \delta_{\CX_0})$-martingale problem has continuous paths, due to Skorohod representation theorem we may assume that
${\mathcal X}^{\mathrm{M},\rho}_t\to{\mathcal X}^{\mathrm{FV}}_t$, as $\rho\to\infty$, uniformly for all $t\in[0,T]$, almost surely.

We will rely on Theorem~3.9 in \cite{LKliem}. That is, as the solutions of the $(L^{\mathrm{FV}}, \Pi^{1,0}, \delta_{\CX_0})$--martingale problem have continuous paths almost surely, we have to construct for each finite $A\subseteq V$, for each $t\in[0,T]$, and $\eps,\delta,\rho>0$ a function $h_{t,\eps,A}\in{\mathcal H}$ (recall from (\ref{e:calH})) and
a random measurable set $Y^\rho_{t,\eps,\delta,A}\subseteq X^\rho_t$ such that
\begin{equation}
\label{e:Thm3.7}
   \limsup_{\delta\downarrow 0}\limsup_{\rho\to\infty} \mathbb{P}\{\mu_t(X^\rho\setminus Y^\rho_{t,\eps,\delta,A})\le h_{t,\eps,A}(\delta),\,\forall t\in[0,T]\}\ge 1-\eps.
\end{equation}\smallskip

Assume first that the geographic space $V$ is finite.
For the construction of such a function $h_{\eps,V}=h_\eps\in{\mathcal H}$ and
a random measurable set $Y^\eps_{\rho,\delta,V}=Y^\eps_{\rho,\delta}\subseteq X^\rho_t$, we can proceed exactly as in the proof of Theorem~4.3 of \cite{LKliem} where the statement is shown with mutation rather than migration in the non-spatial rather than the finite geographic space.
\smallskip

Let now $V$ be countable, and consider a sequence $(V_n)_{n\in\mathbb{N}}$ of finite sets with $V_n\subseteq V$ and $V_n\uparrow V$. Consider for each $n\in\mathbb{N}$ a solution,
${\mathcal X}^{\mathrm{FV},V_n}$, of the $(L^{\mathrm{FV},V_n}, \Pi^{1,0}, \delta_{\CX_0})$-martingale problem with $L^{\mathrm{FV},V_n}$ as defined in (\ref{gs2Vn}).
As we have seen above, ${\mathcal X}^{\mathrm{FV},V_n}_t\in\mathbb{U}^{V_n}_{\mathrm{fct}}$ for all $t\ge 0$, almost surely. Moreover,
we have shown in Step~3 that each solution ${\mathcal X}^{\mathrm{FV}}$ of the $(L^{\mathrm{FV}}, \Pi^{1,0}, \delta_{\CX_0})$-martingale problem on $V$ can be obtained as the limit of ${\mathcal X}^{\mathrm{FV},V_n}$ as $n\to\infty$. To conclude from here that also ${\mathcal X}^{\mathrm{FV}}_t\in\mathbb{U}^{V}_{\mathrm{fct}}$ for all $t\ge 0$, almost surely, fix a finite set $A\subseteq V$. As done before we denote by {${\mathcal X}^{\mathrm{FV},V_n,A}=(X_t,r^{V_n}_t,\mu^{V_n}_t(\boldsymbol{\cdot}\times(\boldsymbol{\cdot}\cap A)))_{t\ge 0}$ and ${\mathcal X}^{\mathrm{FV},A}=(X_t,r_t,\mu_t(\boldsymbol{\cdot}\times(\boldsymbol{\cdot}\times A)))_{t\ge 0}$} the restrictions of ${\mathcal X}^{\mathrm{FV},V_n}$ and ${\mathcal X}^{\mathrm{FV}}$, respectively, to marks in $A$ (compare (\ref{e:restrict})).

For each $m>n$ we couple ${\mathcal X}^{\mathrm{FV},V_n}$ and ${\mathcal X}^{\mathrm{FV},V_m}$  through
the graphical lookdown construction by using
the same Poisson point processes and marking every path which leaves $V_n$ in the $V_m$ dynamics by a~$1$. Moreover, we impose the rule that the 1 is inherited upon lookdown in the sense that {\em both} new particles
carry type~$1$.
The sampling measure of types then follows an interacting Fleming-Viot (in fact two-type Fisher-Wright) diffusion with selection.
The corresponding Moran models are coupled and converge in the many particle per site limit to a limit
evolution, which is the coupling on the finite geographic spaces and the additional types act upon
resampling as under selection.

By construction, if $x,x^\prime\in\mathrm{supp}(\mu_t(\boldsymbol{\cdot}\times A))$, their distance is the same in ${\mathcal X}^{\mathrm{FV},V_n}$ and ${\mathcal X}^{\mathrm{FV},V_m}$ if both carry type~$0$. Thus for suitably large $n$ (depending on $\eps>0$) such that $A\subseteq V_n$ at any  location in $A$
the relative frequencies of types~$1$ at time $t$ can be made  less than any given $\ve >0$  with probability
$\geq 1-\ve$  by simple random walk
estimates. Namely, if $(Z^b_t)_{t \geq 0}$ is a $a(\cdot, \cdot)$-random walk starting in
$b \in V_n \subseteq V_m$ and $b^\prime\in \CC V_n$,
\bea{ag22}
&&\mathbb{P}(Z^b_t \notin V_n \mbox{ for some } t \in [0,T], Z^b_T \in A)\\
&& \quad
 + \mathbb{P}(Z^{b^\prime}_t \in V_n \mbox{ for some }
t \in [0,T], Z^{b^\prime}_T \in A) \leq \delta_n \to 0 \mbox{ as } n \to \infty \; \forall \; m \geq n. \nonumber
\eea
Then the expected frequency of type~$1$ in locations in $A$ is bounded by $F (\delta_m)$ with $F(\delta) \to 0$
as $\delta \to 0$, which follows from the properties of the Fisher-Wright diffusion with selection
easily via duality.

 As a consequence the supremum along the path of  the difference in variational norm of the distance-mark distributions for the $V_n$ and the
$V_m$-evolution for types in the set $A$ can be bounded by a sequence converging to 0 as $n,m \to \infty$.

Therefore also the limit dynamics on countable $V$ has a mark function.
}
\end{proof}\sm

\begin{proof}[Proof of Theorem~\ref{T.dual}] Let ${\mathcal X}^{\mathrm{FV}}$ be the evolving genealogies of interacting Fleming-Viot diffusions where we have assumed that the symmetrized migration is recurrent. In order to prove ergodicity we proceed in two steps: (1) We start with constructing the limiting object which is {\em tree-valued spatial coalescent}. (2) We then prove convergence of ${\mathcal X}^{\mathrm{FV}}_t$ to this tree-valued spatial coalescent, at $t\to\infty$, for any initial state ${\mathcal X}_0\in\mathbb{U}^V_1$. This immediately implies uniqueness of the invariant distribution. \sm

\noindent{\bf Step~1 (Tree-valued spatial coalescent)} Recall from \cite{GLW05} the spatial coalescent started with infinitely many partition elements per site and with migration mechanism $\bar{a}(v,v')$. If the symmetrized migration is recurrent, we can assign to each realization a marked ultra\-metric space, ${\mathcal K}=(K,r)$, which admits a mark function. In order to equip it with a locally finite measure on the leaves, we consider a coupled family of sub-coalescents  $\{{\mathcal K}^{\varrho},\,\rho>0\}$ such that the number of points of a given mark is Poisson with intensity $\rho$. If we now assign each point in ${\mathcal K}^{\rho}$  mass $\rho^{-1}$, then it follows from \cite[Theorem~3]{GLW05} that there exist a measure $\mu$ on
$K\times V$ such that for each $v\in V$, $\mu(K\times\{v\})=1$. This reflects the spatial dust-free property. Thus, we can use the same arguments used in \cite[Theorem~4]{GPW09} to show that the family $\{({\mathcal K}^\rho,\rho^{-1}\sum_{(x,v)\in K^\rho\times V}\delta_{(x,v)}\}$ is tight, and in fact has exactly one limit point, \begin{equation}
\label{e:exactone}
   {\mathcal K}^{\downarrow}:=\overline{(K,r^{\downarrow},\mu^{\downarrow})}
\end{equation}
\sm

\noindent{\bf Step~2 (Convergence into the tree-valued spatial coalescent)} For all ${\mathcal X}_0\in\mathbb{U}^V_1$ and ${\mathcal K}_0\in\mathbb{S}$, by our duality relation,
using the functions $H=H^{n,\phi}$ from (\ref{gs6}),
{
\begin{equation}
\label{e:longg}
\begin{aligned}
   \mathbb{E}\big[H\big({\mathcal X}^{\mathrm{FV}}_t,{\mathcal K}_0\big)\big]
   &=
   \mathbb{E}\big[H\big({\mathcal X}_0,{\mathcal K}_t\big)\big]
   \\
   &\ttO
   \mathbb{E}^{{\mathcal K}_0}\big[\phi\big((r^{\downarrow}(i,j))_{1\le i<j\le n}\big)\big]
   \\
   &=
   \mathbb{E}\big[H\big({\mathcal K}^{\downarrow},{\mathcal K}_0\big)\big].
\end{aligned}
\end{equation}
}
Once more, as the family $\{H^{n,\phi}(\boldsymbol{\cdot},{\mathcal K});\,n\in\mathbb{N},\phi\in C_{bb}(\R_+^{n\choose 2}),{\mathcal K}\in\mathbb{S}\}$ is convergence determining
by Theorem~\ref{T:polydet}(i) , we can conclude that for all initial conditions,
${\mathcal X}^{\mathrm{FV}}_t$ converges Gromov-$\#$-weakly to the tree-valued spatial coalescent.
\end{proof}

\bigskip


{

\section{Proofs of the properties of CSSM Genealogy Processes}}

\label{S:CSSMprop}
In this section, we prove Proposition~\ref{P:reg} and Theorem~\ref{T:MarkovCont} by using properties of the double
Brownian web $(\CW, \wh \CW)$, which was used to construct the CSSM genealogy process $\CX^{\rm CS}$ in {Section}~\ref{S:IntroCSSM}.
\bigskip

\noindent
{\bf Proof of Prop.~\ref{P:reg}. (a)}: The existence of a mark function $\kappa: X^{\rm CS}_t \to \R$ follows by construction (recall (\ref{e.addsampme}) where the mark is explicitly given).
The continuity of $\kappa$ follows from the property of the dual Brownian web $\wh\CW$. More precisely, if $x_n \to x$ in $(X^{\rm CS}_t, r^{\rm CS}_t)$, then identifying $x_n$ and $x$ with points on $\R$, it follows that $\hat f_{(x_n,t)} \to \hat f_{(x,t)}$ in $\wh\Pi$ for some path $\hat f_{(x_n,t)}\in \wh\CW(x_n,t)$ and $\hat f_{(x,t)}\in \wh\CW(x,t)$, which implies that $x_n\to x$ in $\R$.
\medskip

{\bf (b)}: By (\ref{XCSsplit}), we identify $X^{\rm CS}_t$ with $\R$, where a countable subset $E_t$ is duplicated.
The distance between $x,y\in X^{\rm CS}_t$ is defined to be twice of the time to coalescence between the dual Brownian web paths
$\hat f_{(x,t)}\in \wh\CW(x,t)$ and $\hat f_{(y,t)}\in \wh\CW(y,t)$, if the two paths coalesce above time $0$.
Therefore for $l\in (0,t)$, each ball $B^l_i$ of radius $l$ correspond to a maximal interval $[v_{i-1}, v_i]\subset\R$,
where all paths in $\{\hat f_{(x,t)} \in \wh\CW(x,t): x\in (v_{i-1}, v_i)\}$ coalesce into a single path by time $t-l$.
The {{}collections} of such maximal intervals $(v_{i-1}, v_i)$ form a partition of $\R\backslash E^l_t$, where $E^l_t = \{v_i : i\in\Z\}$
is exactly the set defined in (\ref{Elt}).
\medskip

{\bf (c)}: Fix an $l\in (0,t)$. By construction, for each ball $B^l_i$ of radius $l$, $\kappa(B^l_i)=[v_{i-1}, v_i]$ is assigned the Lebesgue measure on $\R$. Together with (b), it implies that $\mu^{\rm CS}_t(X^{\rm CS}_t \times \cdot)$ is the Lebesgue measure on $\R$.
\qed
\bigskip

\noindent
{\bf Proof of Theorem~\ref{T:MarkovCont}}. Let us fix a realization of the Brownian web $\CW$ and its dual $\wh\CW$, and let $\CX^{\rm CS}$ be constructed from $\wh\CW$ as just before Proposition~\ref{P:reg}. By (\ref{XCSsplit}), for each $t>0$, we can identify $X^{\rm CS}_t$ with $A_t \cup E_t^+\cup E_t^-$, where $A_t$ and $E_t$ are defined as in (\ref{At}). To simplify notation, we will drop the superscript CS in the remainder of the proof.
\medskip

{\bf (a)}: The Markov property of $(\CX_t)_{t\geq 0}$ follows from the Markov property of
$\CW$ and $\wh\CW$. More precisely, if we denote by $f|_s^t$ the restriction of a path $f\in \Pi$ to the time interval $[s,t]$, and
$K|_s^t :=\{f|_s^t : f\in K\}$ for a set of paths $K\subset \Pi$, then $\CW|_0^s$ is independent of $\CW|_s^\infty$ for each
$s\geq 0$. The same is true for $\wh\CW$ since $\CW|_0^s$ and $\wh\CW|_0^s$ a.s.\ uniquely determine each
other by Theorem \ref{T:dwebchar}. The Markov process $(\CX_t)_{t\geq 0}$ is time homogeneous because $\CW|_0^t$ is equally distributed with $\CW|_s^{s+t}$, apart from
a time shift.
\medskip

{\bf (b)}: Let $\CX_0 \in \U^\R_1$. We first prove that $(\CX_t)_{t\geq 0}$ is a.s.\ continuous in $t> 0$. To accomplish this, since $\Pi^{1,2}$ is convergence determining in $\U^\R$ as shown in Theorem~\ref{T:polydet}, it suffices to show that for any $\Phi:=\Phi^{n,\phi, g}\in \Pi^{1,2}$,
the evaluated polynomial
\be{poly1}
\Phi({\CX}_t) = \int_{\R^n} \phi(\underline{\underline{r_t}}) g(\ux) \dd \ux
\ee
is continuous in $t$, where $\ux:=(x_1, \ldots, x_n)$, $\underline{\underline{r_t}}:= (r_t(x_i, x_j))_{1\leq i<j\leq n}$,
and given the identification between $X_t$ and $A_t\cup E_t^+\cup E_t^-$ and the definition of $\kappa$, we have replaced integration w.r.t.\ $\mu_t$
on $X_t\times \R$ by integration w.r.t.\ the Lebesgure measure on $\R$. By (\ref{rtxy}), for Lebesgue a.e.\ $x_i, x_j\in \R$, we have
\be{poly2}
r_t(x_i, x_j) = \left\{
\begin{aligned}
2(t-\hat \tau) \qquad \qquad &\qquad \mbox{if} \ \hat\tau\geq 0, \\
2t+ r_0(\xi(u), \xi(v))  &\qquad \mbox{if}\ \hat\tau<0,
\end{aligned}
\right.
\ee
where $\hat\tau$ is the time of coalescence between $\hat f_{(x,t)}$ and $\hat f_{(y,t)}$, $\hat f_{(x,t)}(0)=u$, and $\hat f_{(x,t)}(0)=v$.
By Lemma~\ref{L:spts}, for each $t>0$, $\wh\CW(x,t)$ contains a single path for all but a countable number of $x\in\R$. For such $x$, by Lemma~\ref{L:pathconv},
the time of coalescence between $\hat f_{(x,s)}$ and $\hat f_{(x,t)}$ tends to $t$ as $s\to t$, and hence $\lim_{s\to t} r((x,s), (x,t)) =0$, where $r((x,s), (x,t))$ is defined in (\ref{rxtys}) and extends the definition of $r_t(x,y)$ to individuals at different times. Since
\be{poly3}
|r_s(x_i, x_j)-r_t(x_i, x_j)| \leq r((x_i,s), (x_i,t)) + r((x_j,s), (x_j,t)),
\ee
it follows that when $t>0$, for Lebesgue a.e.\ $x_i, x_j\in\R$, $1\leq i<j\leq n$, we have
\be{rcont}
\lim_{s\to t} r_s(x_i, x_j) = r_t(x_i, x_j).
\ee
We can then apply the dominated convergence theorem in (\ref{poly1}) to deduce that, for each $t>0$, a.s.\
\begin{equation}\label{limstPhi}
\lim_{s\to t} \Phi(\CX_s)=\Phi(\CX_t).
\end{equation}
This verifies that $(\CX_t)_{t\geq 0}$ is a.s.\ continuous in $t>0$.

Proving the a.s.\ continuity of $(\CX_t)_{t\geq 0}$ at $t=0$ poses new difficulties because $\CX_0$ can be any state in $\U^\R_1$,
while for any $t>0$, $\CX_t$ is a regular state as shown in Proposition~\ref{P:reg}. We get around this by showing that $\CX$ admits
a c\`adl\`ag version. More precisely, we invoke a part of the proof of the convergence Theorem~\ref{T:conv} that is independent
of the current proof. Note that for any $\CX_0\in \U^\R_1$, we can find a sequence $\CX^{\rm FV, \epsilon}_0\in \U^\Z_1$, indexed
by $\epsilon>0$, such that $S_\epsilon \CX^{\rm FV, \epsilon}_0\to \CX_0$. Indeed, we only need to approximate the mark space
$\R$ by $\epsilon \Z$ in order to construct $S_\epsilon\CX^{\rm FV, \epsilon}_0$ from $\CX_0$. In the proof of Theorem~\ref{T:conv},
it is shown that the corresponding sequence of interacting Fleming-Viot genealogy process
$(S_\epsilon \CX^{\rm FV,\epsilon}_{\epsilon^{-2}t})_{t\geq 0}$ is a tight family of $D([0,\infty), \U^\R)$-valued random variables,
where $D([0,\infty), \U^\R)$ denotes the space of c\`adl\`ag paths on $\U^\R$ equipped with the Skorohod topology. Furthermore,
$(S_\epsilon \CX^{\rm FV,\epsilon}_{\epsilon^{-2}t})_{t\geq 0}$ converges in finite-dimensional distribution to the CSSM genealogy process
$(\CX_t)_{t\geq 0}$. Therefore, $(\CX_t)_{t\geq 0}$ must admit a version which is a.s.\ c\`adl\`ag, with $\CX_t\to \CX_0$ as $t\downarrow 0$.
Since we have just shown that the version of $(\CX_t)_{t\geq 0}$ constructed in Sec.~\ref{S:IntroCSSM} is a.s.\ continuous in $t>0$,
 it follows that the same version must also be a.s.\ c\`adl\`ag, and hence continuous at $t=0$, which concludes the proof of part (b).
\medskip

{\bf (c)}: To prove the Feller property, let $\CX_0^{(m)}\to \CX_0$ in $\U^\R_1$, and let $t_m\to t\geq 0$. To show $\CX_{t_m}^{(m)} \Rightarrow \CX_t$, by Theorem~\ref{T:polydet}, it suffices to show
\begin{equation}\label{limnPhin}
\lim_{m\to\infty} \E[\Phi(\CX_{t_m}^{(m)})] = \E[\Phi(\CX_t)]  \qquad \forall\ \Phi=\Phi^{n,\phi, g}\in \Pi^{1,2}.
\end{equation}
We claim that the convergence in
\begin{equation}\label{limnPhinUnif}
\lim_{s\to t}\E[\Phi(\CX_s)]= \E[\Phi(\CX_t)] \quad \mbox{is uniform w.r.t.\ the initial condition } \CX_0.
\end{equation}
In particular, as $m\to\infty$,
$$
\E[\Phi(\CX^{(m)}_{t_m})]-\E[\Phi(\CX^{(m)}_t)]\to 0.
$$
To prove (\ref{limnPhin}), it then suffices to show that
\begin{equation}\label{limnPhint}
\lim_{m\to\infty} \E[\Phi(\CX^{(m)}_t)] = \E[\Phi(\CX_t)].
\end{equation}
We prove (\ref{limnPhinUnif}) and (\ref{limnPhint}) next.
\medskip

{\em Proof of $(\ref{limnPhinUnif})$.} By the Markov property of $\CX$, it suffices to show that as $t\downarrow 0$,
\begin{equation}\label{limnPhinUnif2}
\big|\E[\Phi(\CX_t)]-\E[\Phi(\CX_0)]\big| \to 0 \quad \mbox{uniformly in } \CX_0.
\end{equation}
Note that we can write
\begin{equation}\label{EPhiCX0}
\E[\Phi(\CX_0)] = \int_{\R^n}  g(\ux)\,{ \E[\phi(\{r_0(\xi(x_i), \xi(x_j))\}_{1\leq i<j\leq n})]} \, \dd \ux,
\end{equation}
where for each $x\in\R$, $\xi(x)\in X_0$ is sampled according to the conditional distribution of $\mu_0$ on $X_0$, conditioned on the spatial coordinate in $X_0\times \R$ being equal to $x$. On the other hand,
\be{rf1}
\E[\Phi(\CX_t)] = \int_{\R^n} g(\ux) \E[\,\phi(\{r_t(x_i, x_j)\}_{1\leq i<j\leq n})] \, \dd \ux.
\ee

Let $F(\ux,t)$ denote the event that the dual Brownian web paths $\hat f_{(x_1, t)}, \ldots, \hat f_{(x_n,t)}$ do not coalesce during the time interval $[0,t]$. We can then partition $\E[\phi(\{r_t(x_i, x_j)\}_{1\leq i<j\leq n})]$ into expectation restricted to $F(\ux, t)$ and $F^c(\ux, t)$ respectively. On the event $F(\ux, t)$, we can replace $\hat f_{(x_i,t)}$, $1\leq i\leq n$, by independent Brownian motions $(x_i(s))_{s\leq t}$, $1\leq i\leq n$, starting respectively at $x_i$ at time $t$ and {\em running backward in time}. Then
\be{rf2}
r_t(x_i, x_j) = 2t + r_0(\xi(x_i(0)), \xi(x_j(0))).
\ee
Let $\wt F(\ux, t)$ denote the event that $(x_i(s))_{s\leq t}$, $1\leq i\leq n$, do not intersect during the time interval $[0,t]$.
Then using (\ref{rf2}), we can rewrite (\ref{rf1}) as
\begin{equation}\label{EPhiCXt}
\begin{aligned}
\E[\Phi(\CX_t)] & = \int_{\R^n} \!\!\! g(\ux) \E[\,\phi(\underline{\underline{r_t}}) 1_{F^c(\ux, t)}]\, \dd \ux + \int_{\R^n} \!\!\! g(\ux) \E[\,\phi(\{2t+r_0(\xi(x_i(0)), \xi(x_j(0)))\}_{1\leq i<j\leq n})] \,\dd \ux \\
& \qquad - \int_{\R^n} g(\ux) \E[\,\phi(\{2t+r_0(\xi(x_i(0)), \xi(x_j(0)))\}_{1\leq i<j\leq n}) 1_{\wt F^c(\ux, t)}] \,\dd \ux.
\end{aligned}
\end{equation}
Note that for Lebesgue a.e.\ $x_1, \ldots, x_n$, $\P(F^c(\ux, t))=\P(\wt F^c(\ux, t))\to 0$ as $t\downarrow 0$. Therefore by the bounded convergence theorem, the first and third {{}term on the right hand side of} (\ref{EPhiCXt}) converges to $0$ as $t\downarrow 0$, uniformly in $\CX_0$.
For the second {{}term on the right hand side of}  (\ref{EPhiCXt}), we can make the change of variable $y_i:= x_i(0)$, $y_i(s):=x_i(s)$, to rewrite it as
\be{rs1}
\int_{\R^n} \E[ g(y_1(t), \ldots, y_n(t)) \,\phi(\{2t+r_0(\xi(y_i), \xi(y_j))\}_{1\leq i<j\leq n})] \,\dd \uy.
\ee
For each $\uy\in \R^n$, clearly the quantity inside the expectation converges a.s.\ to the analogue in (\ref{EPhiCX0}) as $t\downarrow 0$, and the
speed of convergence does not depend on $\CX_0$. Therefore the expectation in (\ref{rs1}) also converges, uniformly in $\CX_0$. Using the fact that $g$ has bounded support, while $g$ and $\phi$ are both bounded, we can easily dominate the integrand in (\ref{rs1}) w.r.t.\ $\dd \uy$  by an integrable function as $t\downarrow 0$; (\ref{limnPhinUnif2}), and {{}hence (\ref{limnPhinUnif}) follows}.
\medskip

{\em Proof of $(\ref{limnPhint})$.}
For each $x_i\in \R$, $1\leq i\leq n$, let us denote $\hat f_{(x_i, t)}\in \wh\CW(x_i, t)$ by $\hat f_i$. Then
\be{rf4}
\E[\Phi(\CX^{(m)}_t)] = \int_{\R^n} g(\ux) \E[\phi( \underline{\underline{r_t}}^{\!(m)})] \dd \ux,
\ee
where $\underline{\underline{r_t}}^{\!(m)}=\{r^{(m)}_t(x_i, x_j)\}_{1\leq i<j\leq n}$ depends on the realization of $(\hat f_i)_{1\leq i\leq n}$. Let $\hat \tau$ be the smallest time in $[0,t]$ when a coalescence occurs among the paths $(\hat f_i)_{1\leq i\leq n}$. Let $I_1, \ldots, I_k$ denote the partition of $\{1, \ldots, n\}$, where all $\hat f_i$ with $i$ in the same partition element $I_j$ have coalesced into a single path at time $\hat \tau$. Conditioned on $(\hat f_i)_{1\leq i\leq n}$ on the time interval $[\hat \tau, t]$, the distribution of the remaining $k$ coalescing Brownian motions on the time interval $[0,\hat\tau]$ is then given by the distribution of $k$ Brownian motions conditioned not to intersect on the time interval $[0,\hat\tau]$, and their positions at time $0$ has a probability density in $\R^k$, which we denote by $g_{\hat f,\hat \tau}(v_1, \ldots, v_k)$. Note that conditioned on $(\hat f_i)_{1\leq i\leq n}$ on the time interval $[\hat \tau,t]$ and their positions $v_1, \ldots, v_k$ at time $0$, $\phi( \{r^{(m)}_t(x_i, x_j)\}_{1\leq i<j\leq n})$ only depends on $\{r_0^{(m)}(\xi(v_i), \xi(v_j))\}_{1\leq i<j\leq k}$, cf.~(\ref{poly2}). We can therefore write
\be{rf5}
\begin{aligned}
\Phi_{\hat f, \hat \tau}(\CX_0^{(m)})& :=\E[\phi( \underline{\underline{r_t}}^{\!(m)})\,|\,(\hat f_i(s))_{1\leq i\leq n, s\in [\hat\tau, t]}]\\
& = \int_{\R^k} \!\! g_{\hat f, \hat\tau}(v_1, \ldots, v_k) \E[\phi_{\hat f, \hat\tau}(\{r_0^{(m)}(\xi(v_i), \xi(v_j))\}_{1\leq i<j\leq k})]\, \dd\uv,
\end{aligned}
\ee
where given the realization of $(\hat f_i(s))_{1\leq i\leq n, s\in [\hat\tau, t]}$ and $\{\hat f_i(0)\}_{1\leq i\leq n}=\{v_1, \ldots, v_k\}$,
\be{rf6}
\phi_{\hat f, \hat\tau}(\{r_0^{(m)}(\xi(v_i), \xi(v_j))\}_{1\leq i<j\leq k})= \phi( \{r_t^{(m)}(x_i, x_j)\}_{1\leq i<j\leq n}).
\ee
Note that $\Phi_{\hat f, \hat \tau}$ is a polynomial of order $k$ on $\U^\R$, defined from the bounded continuous functions $g_{\hat f, \hat \tau}$ and $\phi_{\hat f, \hat\tau}$, except that $g_{\hat f, \hat \tau}$ does not have bounded support. Nevertheless, $g_{\hat f, \hat \tau}$ is integrable and can be approximated by continuous functions with bounded support. Therefore from the assumption $\CX^{(m)}_0\to \CX_0$ in $\U^\R_1$, we deduce that $\Phi_{\hat f, \hat \tau}(\CX^{(m)}_0)\to \Phi_{\hat f, \hat \tau}(\CX_0)$ for Lebesgue a.e.\ $x_1, \ldots, x_n$ and a.e.\ realization of $(\hat f_i(s))_{1\leq i\leq n, s\in [\hat \tau,t]}$. It then follows from the bounded convergence theorem that
\begin{equation}
\E[\Phi(\CX^{(m)}_t)] = \int_{\R^n} g(\ux) \E[\Phi_{\hat f, \hat \tau}(\CX^{(m)}_0)] \dd \ux \asto{m} \int_{\R^n} g(\ux) \E[\Phi_{\hat f, \hat \tau}(\CX_0)] \dd \ux = \E[\Phi(\CX_t)],
\end{equation}
which concludes the proof of the Feller property.
\medskip

{\bf (d):} This follows readily from the construction of the CSSM genealogy process $\CX^{\rm CS}_t$. For any $N>0$, the genealogical distances among individuals with spatial locations in $[-N,N]$ are determined by coalescing Brownian motions in the dual Brownian web $\widehat\CW$. The initial condition $\CX^{\rm CS}_0$ affects the genealogical distances only on the event that the backward coalescing Brownian motions starting from $[-N,N]$ at time $t$ do not coalesce into a single path by time $0$. As $t\to\infty$, the probability of this event tends to $0$, and therefore $(\CX^{\rm CS}_{t+s})_{s\geq 0}$ converges in distribution to the CSSM genealogy process constructed from $\widehat\CW$ as in Section \ref{S:IntroCSSM}, with initial condition being at time $-\infty$.
\qed

\section{Proof of convergence of Rescaled IFV Genealogies}\label{S:conv}
In this section, we prove Theorem~\ref{T:conv}, that under diffusive scaling of space and time as well as rescaling of measure,
the genealogies of the interacting Fleming-Viot process converges to those of a CSSM. In Section \ref{ss.convfin}
we prove f.d.d.-convergence, and in Section \ref{ss.tight} tightness in path space. In Section \ref{ss.conv2},
{we prove} Theorem~\ref{T:conv2} on the measure-valued process, which is needed to prove tightness in Section \ref{ss.tight}.

As in Theorem~\ref{T:conv}, let
$(\CX^{\rm FV, \eps}_t)_{t\geq 0}=\overline{(X_t^{\rm FV, \eps}, r_t^{\rm FV,\eps}, \mu_t^{\rm FV, \eps})}_{t\geq 0}$
be the family of IFV genealogy processes on $\Z$ indexed by $\eps>0$, such that $S_\eps \CX_0^{\rm FV, \eps}\to \CX_0^{\rm CS}\in \U^\R_1$
as $\eps\to 0$, and let $(\CX^{\rm CS}_t)_{t\geq 0}=\overline{(X_t^{\rm CS}, r_t^{\rm CS}, \mu_t^{\rm CS})}_{t\geq 0}$
be the CSSM genealogy process with initial condition $\CX^{\rm CS}_0$.

\subsection{Convergence of Finite-dimensional Distributions}
\label{ss.convfin}
In this subsection we prove the convergence $(S_\eps\CX^{\rm FV, \eps}_{\eps^{-2}t})_{t\geq 0}\Rightarrow (\CX^{\rm CS}_t)_{t\geq 0}$ in finite-dimensional distribution, i.e.,
\begin{equation}\label{fddconv}
(S_\eps\CX^{\rm FV, \eps}_{\eps^{-2}t_1}, \ldots, S_\eps\CX^{\rm FV, \eps}_{\eps^{-2}t_k}) \Astoo{\eps} (\CX^{\rm CS}_{t_1}, \ldots, \CX^{\rm CS}_{t_k}) \qquad \forall\ 0\leq t_1<t_2<\cdots<t_k,
\end{equation}
where $\Rightarrow$ denotes weak convergence of $(\U^\R)^k$-valued random variables. By \cite[Prop.~3.4.6]{EK86}
on convergence determining class for product spaces, it suffices to show that for any
$\Phi_i:=\Phi^{n_i, \phi_i, g_i} \in \Pi^{1,2}$, $1\leq i\leq k$, we have (recall (\ref{Phingphi}))
\begin{equation}\label{fddconvPhi}
\E\Big[\prod_{i=1}^k \Phi_i(S_\eps\CX^{\rm FV, \eps}_{\eps^{-2}t_i})\Big] {\astoo{\eps}} \E\Big[\prod_{i=1}^k \Phi_i(\CX^{\rm CS}_{t_i})\Big].
\end{equation}
For notational convenience we assume first that the initial tree is the trivial one
(all distances are zero) and
we shall see at the end of the argument that this easily generalizes.
We will first rewrite both sides of this convergence relation in terms of the dual coalescents and then apply the invariance principle
for coalescing random walks.
\medskip

\noindent
{
{\bf Step~1 (Claim rephrased in terms of coalescents).}
}
We can by the definition of the polynomial in (\ref{poly1}) rewrite the left hand side of (\ref{fddconvPhi}) as
\begin{equation}\label{EprodPhiFV}
\begin{aligned}
& \E\Big[\prod_{i=1}^k \Phi_i(S_\eps \CX^{\rm FV, \eps}_{\eps^{-2}t_i})\Big] \\
=\ &  \E\Big[\prod_{i=1}^k (\eps\sigma^{-1})^{n_i} \!\!\!\!\!\!\!\!\!
\sum_{x^i_1, \ldots, x^i_{n_i}\in \Z}\!\!\!\!\!\! g_i\big(\eps\sigma^{-1}(x^i_a)_{1\leq a\leq n_i}\big) \phi_i
\big(\eps^2 (r^{\rm FV, \eps}_{\eps^{-2}t_i}(\xi^\eps_{t_i}(x^i_a), \xi^\eps_{t_i}(x^i_b)))_{1\leq a<b\leq n_i}\big)\Big],
\end{aligned}
\end{equation}
where for each time $t_i$, we sample $n_i$ individuals in $X^{\rm FV, \eps}_{\eps^{-2}t_i}$ at respective spatial positions
 $(x^i_a)_{1\leq a\leq n_i}$, with
 \be{rs2}
 \xi^\eps_{t_i}(x_a) \mbox{ being sampled from } X^{\rm FV, \eps}_{\eps^{-2}t_i} \mbox{ according to }
 \mu^{\mathrm{FV}}_{\ve^{-2} t_i} (\cdot |x_a),
 \ee
the conditional distribution of $\mu^{\rm FV, \eps}_{\eps^{-2}t_i}$ on $X^{\rm FV, \eps}_{\eps^{-2}t_i}$ conditioned on
the spatial mark being equal to $x_a\in\Z$.

By the space-time duality relation (\ref{gs9}) for the IFV genealogy processes,
every summand of the R.H.S.\ of (\ref{EprodPhiFV}) can be calculated in terms of coalescing random walks. Namely, the joint law of the space-time genealogies of the sampled individuals $\xi^\eps_{t_i}(x^i_a)\in X^{\rm FV, \eps}_{\eps^{-2}t_i}$, $1\leq a\leq n_i$ and $1\leq i\leq k$, is equal to that of a collection of coalescing random walks $(x^i_a(s))_{s\leq \eps^{-2}t_i}$ (recall here time $s$ runs {\em backwards}), starting respectively at $x^i_a$ at time $\eps^{-2} t_i$, where each walk evolves {\em backward} in time as rate $1$
continuous time random walk on $\Z$ with transition probability kernel $\bar a$, and two walks at the same location coalesce at rate $\gamma$.
{{}From the duality relation we get the following the stochastic representation:}
\begin{equation}\label{rFVeps}
r^{\rm FV, \eps}_s(x,y) :=
\left\{
\begin{aligned}
2 (s-\hat \tau) \qquad \quad & \qquad \mbox{ if } \hat\tau\geq 0, \\
2 s + r^{\rm FV, \eps}_0(\xi^\eps_0(u), \xi^\eps_0(v))  & \qquad \mbox{ if } \hat\tau< 0,
\end{aligned}
\right.
\end{equation}
where $\hat \tau$ denotes the time of coalescence between the two coalescing walks starting at $x, y\in\Z$ at time $s$,
while $u,v$ are the positions of the two walks at time $0$.

We observe next that the continuum population is represented by
$A_t \cup E^+_t \cup E^-_t$ which is a version of $\R$ marked on $E_t$ by $+,-$,
the geographic marks are the reals and since the sampling measure is Lebesgue measure, we
can write polynomials based on integration over $\R$ instead of $X \times V$ as:
\be{ag26}
\intl_X \intl_V \mu^{\otimes n} (d(\ux, \uu)) g(\uu) \varphi (\uur(\ux))
= \intl_\R \lambda^{\otimes n} (d \uu) g(\uu) \varphi(\uur(\uu)).
\ee

We can rewrite using the duality of Corollary \ref{C.dual}, see also (\ref{rxtys}), the R.H.S.\ of (\ref{fddconvPhi}) in the same form as in (\ref{EprodPhiFV}):
\begin{equation}\label{EprodPhiCS}
\E\Big[\prod_{i=1}^k \Phi_i(\CX^{\rm CS}_{t_i})\Big] = \E\Big[\prod_{i=1}^k\int_{\R^{n_i}} g_i\big((y^i_a)_{1\leq a\leq n_i}\big) \phi_i\big((r^{\rm CS}_{t_i}(y^i_a, y^i_b))_{1\leq a<b\leq n_i} \big)\dd \uy^i\Big],
\end{equation}
where at each time $t_i$, we sample $n_i$ individuals from $X^{\rm CS}_{t_i}$ according to $\mu^{\mathrm{CS}}_{t_i}$ (which is Lebesgue measure on $\R$) at positions $(y^i_a)_{1\leq a\leq n_i}$,
 and their joint space-time genealogy lines are by construction distributed as a collection of coalescing Brownian motions
  $(y^i_a(s))_{s\leq t_i}$, evolving {\em backward} in time.

To link (\ref{EprodPhiFV}) with (\ref{EprodPhiCS}), we note that in (\ref{EprodPhiFV}), we can regard
\be{rf7}
\prod_{i=1}^k (\eps\sigma^{-1})^{n_i} \!\!\!\!\!\!\!\!\! \sum_{x^i_1, \ldots, x^i_{n_i}\in \Z}\!\!\!\!\!\! g_i\big(\eps\sigma^{-1}(x^i_a)_{1\leq a\leq n_i}\big) \delta_{x^i_1}\cdots \delta_{x^i_{n_i}}
\ee
as a finite signed sampling measure (recall that $g$ has bounded support), which is easily seen to converge weakly to the finite signed sampling measure
appearing in (\ref{EprodPhiCS}), namely
\be{rf8}
\prod_{i=1}^k g_i\big((y^i_a)_{1\leq a\leq n_i}\big) \dd y^i_1\cdots \dd y^i_{n_i}.
\ee

To prove (\ref{fddconvPhi}), it then suffices to show that (having (\ref{rf7})-(\ref{rf8}) in mind): If for each $1\leq a\leq n_i$ and $1\leq i\leq k$, $x^{i, \eps}_a\in\Z$ and $\eps\sigma^{-1} x^{i,\eps}_a \to y^i_a$ as $\eps\to 0$, then
\begin{equation}\label{EPhiconv}
\E\Big[\prod_{i=1}^k \!\phi_i\big( (\eps^2r^{\rm FV, \eps}_{\eps^{-2}t_i}(\xi^\eps_{t_i}(x^{i,\eps}_a), \xi^\eps_{t_i}(x^{i,\eps}_b)))_{1\leq a<b\leq n_i}\big)\Big]
\astoo{\eps} \E\Big[\prod_{i=1}^k \!\phi_i\big((r^{\rm CS}_{t_i}(y^i_a, y^i_b))_{1\leq a<b\leq n_i} \big)\Big]. \!\!\!\!\!\!
\end{equation}
\medskip

\noindent
{
{\bf Step~2 (Invariance principle for coalescents).}
}
We next prove (\ref{EPhiconv}) by means of an invariance principle for coalescing random walks.

This {\em invariance principle} reads as follows.
Given a collection of backward coalescing random walks starting at $\eps$-dependent positions $x^{i,\eps}_a\in\Z$ at time
$\eps^{-2}t_i$, $1\leq a\leq n_i$, $1\leq i\leq k$, such that $\eps\sigma^{-1} x^{i,\eps}_a \to y^i_a$ as $\eps\to 0$, the collection of coalescing random walks
$(x^{i, \eps}_a(s))_{s\leq \eps^{-2}t_i}$, rescaled diffusively as $(\eps\sigma^{-1}x^{i,\eps}_a(\eps^{-2}s))_{s\leq t_i}$,
converges in distribution to the collection of coalescing Brownian motions $(y^i_a(s))_{s\leq t_i}$ evolving backward in time. Furthermore,
the times of coalescence between the coalescing random walks, scaled by $\eps^2$, converge in distribution to the times
of coalescence between the corresponding Brownian motions. The proof of such an invariance principle can be easily
adapted from \cite[Section 5]{NRS05}, which considered discrete time random walks with instantaneous coalescence.
We will omit the details.

Let $\delta>0$ be small. Note that the collection of rescaled coalescing random walks $(\eps\sigma^{-1}x^{i, \eps}_a(\eps^{-2}s))$ restricted to the time interval $s\in [\delta, t_k]$, together with their times of coalescence, converge in joint distribution to the collection of coalescing Brownian motions $(y^i_a(s))$ restricted to the time interval $s\in [\delta, t_k]$, together with their times of coalescence. Using Skorohod's representation theorem (see e.g.~\cite{B89}), we can couple $(\eps\sigma^{-1}x^{i, \eps}_a(\eps^{-2}s))_{s\in [\delta, t_k]}$ and $(y^i_a(s))_{s\in [\delta, t_k]}$ such that the paths and their times of coalescence converge almost surely. {\em Let us assume such a coupling} from now on.

By the same argument as in the proof of Theorem~\ref{T:MarkovCont}~(c), we can rewrite the expectations in (\ref{EPhiconv}) in terms of the backward coalescing random walks $x^{i,\eps}_a\in\Z$ and coalescing Brownian motions $y^i_a$. Furthermore, we can condition on the coalescing random walks $x^{i,\eps}_a(s)$ on the time interval $[\delta\eps^{-2}, \eps^{-2}t_k]$ and condition on the coalescing Brownian motions $y^i_a(s)$ on the time interval $[\delta, t_k]$, coupled as above.

Given the locations $u^\eps_1, \ldots, u^\eps_l\in\Z$ of the remaining coalescing random walks at time $\delta\eps^{-2}$, we now make an approximation and replace them by {\em independent} random walks on the remaining time interval $[0, \delta \eps^{-2}]$, and make a similar replacement for the coalescing Brownian motions. Note that the error we introduce to the two sides of
(\ref{EPhiconv}) is bounded by a constant (determined only by $|\phi_i|_\infty$, $1\leq i\leq k$) times the probability that there is a coalescence among the random walks (resp.\ Brownian motions) in the time interval $[0, \delta \eps^{-2}]$ (resp.\ $[0, \delta]$), which tends to $0$ as $\delta\downarrow 0$ uniformly in $\eps$ by the properties of Brownian motion and the invariance principle. Therefore to prove (\ref{EPhiconv}), it suffices to prove its analogue where we make such an approximation for a fixed $\delta>0$, replacing coalescing random walks (resp.\ Brownian motions) on the time interval $[0, \delta\eps^{-2}]$ (resp.\ $[0, \delta]$) by independent ones. Let us fix such a $\delta>0$ from now on.

By conditioning on the coalescing random walks and the coalescing Brownian motions on the macroscopic time interval $[\delta, t_k]$ and using the a.s.\ coupling between them, we note that the analogue of (\ref{EPhiconv}) discussed above follows readily if we show:
\begin{lemma}\label{L.eps1}{}
If $u^\eps_1, \ldots, u^\eps_l\in \Z$ satisfy $\eps\sigma^{-1}u^\eps_i\to u_i$ as $\eps\to 0$, then for any bounded continuous function $\psi: \R^{l\choose 2} \to\R$, we have
\begin{equation}\label{gepsconv}
\begin{aligned}
&\sum_{x_1, \ldots, x_l \in \Z} g_\delta^\eps(x_1, \ldots, x_l)
\E\big[\psi\big((\eps^2r^{\rm FV, \eps}_0(\xi^\eps_0(x_i), \xi^\eps_0(x_j)))_{1\leq i<j\leq l}\big)\big] \\
\astoo{\eps} & \quad \int_{\R^l} g_\delta(y_1, \ldots, y_l) \E\big[\psi\big((r^{\rm CS}_0(\xi(y_i), \xi(y_j)))_{1\leq i<j\leq l}\big)\big] \dd \uy,
\end{aligned}
\end{equation}
where $g_\delta^\eps(\ux)$ is the probability mass function of $l$ independent random walks at time $\delta\eps^{-2}$, starting at $u^\eps_1, \ldots, u^\eps_l$; while $g_\delta(\uy)$ is the probability density function of $l$ independent Brownian motions at time $\delta$, starting at
$u_1, \ldots, u_l$.
\end{lemma}
\noindent
{\em Proof.} If we can replace $g_\delta^\eps(x_1, \ldots, x_l)$ in (\ref{gepsconv}) by $(\eps\sigma^{-1})^l g_\delta(\eps\sigma^{-1}x_1, \ldots, \eps\sigma^{-1}x_l)$, then (\ref{gepsconv}) follows immediately by applying the polynomial $\Phi^{l, \psi, g_\delta}$ to the states $S_\eps\CX^{\mathrm{FV},\ve}_0$ and $\CX^{\mathrm{CS}}_0$, using the assumption $S_\eps\CX^{\mathrm{FV},\ve}_0\to \CX^{\mathrm{CS}}_0$. The only problem is that $g_\delta$ does not have bounded support as we require for a polynomial. However, it is continuous and integrable, and hence can be approximated by continuous functions with bounded support. Therefore the above reasoning is still valid.

To see why we can replace $g_\delta^\eps(\ux)$ by $(\eps\sigma^{-1})^l g_\delta(\eps\sigma^{-1}\ux)$, note that by the {\em local central limit theorem} (see e.g.~\cite{S76}),
\be{rf11}
(\eps\sigma^{-1})^{-l} g^\eps_\delta(\lceil\eps^{-1}\sigma\uy\rceil) \astoo{\eps} g_\delta(\uy)
\ee
uniformly in $\uy \in [-L, L]^l$ for any $L>0$. Therefore when we restrict the summation in (\ref{gepsconv}) to $\ux\in [-\eps^{-1} L, \eps^{-1} L]^l$, the replacement induces an error that tends to $0$ as $\eps\to 0$. By the central limit theorem, the contribution to the sum in (\ref{gepsconv}) from $\ux \notin [-\eps^{-1} L, \eps^{-1} L]^l$ can be made arbitrarily small (uniformly in $\eps$) by choosing $L$ large, and hence can be safely neglected if we first let $\eps\to 0$ and then let $L\to\infty$.
\qed

\subsection{Tightness}
\label{ss.tight}

In this {subsection} we prove the tightness of the family of rescaled IFV genealogy processes, $(S_\eps \CX^{\rm FV, \eps})_{\eps>0}$, regarded as $C([0,\infty), \U^\R)$-valued random variables.

First we note that it is sufficient to prove the tightness of $(S_\eps \CX^{\rm FV, \eps})_{\eps>0}$ as random variables taking values in the Skorohod space
$D([0,\infty), \U^\R)$. Indeed, the tightness of $(S_\eps \CX^{\rm FV, \eps})_{\eps>0}$ in the Skorohod space, together with the convergence of $S_\eps \CX^{\rm FV, \eps}$ to $\CX^{\rm CS}$ in finite-dimensional distributions, imply that $S_\eps \CX^{\rm FV, \eps}\Rightarrow \CX^{\rm CS}$ as $D([0,\infty), \U^\R)$-valued random variables. In particular, $(\CX^{\rm CS}_t)_{t\geq 0}$ admits a version which is a.s.\ c\`adl\`ag. Together with the fact that $\CX^{\rm CS}_t$ is a.s.\ continuous in $t>0$, which was established in the proof of Theorem~\ref{T:MarkovCont}~(b), it follows that $\CX^{\rm CS}_t$ must be a.s.\ continuous in $t\geq 0$. {\em Note that this concludes the proof of Theorem~\ref{T:MarkovCont}~(b).}

Using Skorohod's representation theorem (see e.g.~\cite{B89}) to couple $(S_\eps \CX^{\rm FV, \eps})_{\eps>0}$ and $\CX^{\rm CS}$ such that the convergence in $D([0,\infty), \U^\R)$ is almost sure, and using the a.s.\ continuity of $(\CX^{\rm CS}_t)_{t\geq 0}$, we can then easily conclude that $S_\eps \CX^{\rm FV, \eps}\to \CX^{\rm CS}$ a.s.\ in $C([0,\infty), \U^\R)$, which implies the tightness of $(S_\eps \CX^{\rm FV, \eps})_{\eps>0}$ as $C([0,\infty), \U^\R)$-valued random variables.
\medskip

By Jakubowski's criterion (see e.g.\ \cite[Theorem 3.6.4]{D93}), to show that $(S_\eps \CX^{\rm FV, \eps})_{\eps>0}$ is a tight family of random variables in the Skorohod space $D([0,\infty), \U^\R)$, it suffices to show that the following two conditions are satisfied:

\begin{itemize}
 \item[{\bf (J1)}] (Compact Containment) For each $T >0$ and $\delta >0$, there exists a compact set $K_{T,\delta} \subset  \U^\R$ such that for all $\eps>0$,
\begin{equation} \label{compactcont}
\P\big(S_\eps \CX^{\rm FV, \eps}_{\eps^{-2}t} \in K_{T,\delta} \,\forall \, 0 \leq t \leq T \big)  \geq 1 - \delta;
\end{equation}

\item[{\bf (J2)}] (Tightness of Evaluations)
For each $f\in \wt\Pi^{1,2}$, $(f(S_\eps \CX^{\rm FV, \eps}_{\eps^{-2}t}))_{t\geq 0}$, indexed by $\eps>0$, is a tight family of $D([0,\infty), \R)$-valued random variables.
\end{itemize}
Note that $\wt\Pi^{1,2}$ (recall from (\ref{gs1a})) separates points in $\U^\R$ by Theorem~\ref{T:polydet}, and is closed under addition.

We first prove {\bf (J2)},  following an approach used in~\cite{AS11}, where a family of rescaled measure-valued processes induced by the voter model
 on $\Z$ is shown to be tight and converge weakly to the measure-valued CSSM as $\eps\to 0$.
We will verify {\bf (J2)} via Aldous' tightness criterion (see e.g.~\cite[Theorem 3.6.5]{D93}), reducing {\bf (J2)} to the following conditions:
\begin{itemize}
 \item[{\bf (A1)}] For each rational $t\geq 0$, $\{f(S_\eps \CX^{\rm FV, \eps}_{\eps^{-2}t})\}_{\eps>0}$ is a tight family of $\R$-valued random variables.
\item[{\bf (A2)}] Given $T>0$ and any stopping time $\tau_\eps\leq T$, if $\delta_\eps\downarrow 0$ as $\eps\downarrow 0$, then for each $\eta>0$,
\begin{equation}\label{aldous}
\lim_{\eps\downarrow 0} \P\big(\big| f(S_\eps \CX^{\rm FV, \eps}_{\eps^{-2}(\tau_\eps+\delta_\eps)}) - f(S_\eps \CX^{\rm FV, \eps}_{\eps^{-2}\tau_\eps}) \big| >\eta \big) =0.
\end{equation}
\end{itemize}
\bigskip

\noindent
{\bf Proof of (J2) via (A1)--(A2).} Note that each $f\in \wt\Pi^{1,2}$ can be written as $f=\sum_{i=1}^k c_i \Phi_i$
for finitely many $\Phi_i\in \Pi^{1,2}$ and $c_i\in\R$. It then follows by the triangle inequality that to verify
 {\bf (A1)--(A2)}, it suffices to consider $f=\Phi\in \Pi^{1,2}$. We abbreviate for  $\Phi\in\Pi^{1,2}$,
 $\Phi^\eps_t:= \Phi(S_\eps \CX^{\rm FV, \eps}_{\eps^{-2}t})$.

{\bf (A1).}
Note that for $\Phi=\Phi^{n,\phi, g} \in \Pi^{1,2}$, $(\Phi^\eps_t)_{\eps>0}$ is uniformly bounded, because $\phi$ is bounded,
$g$ is bounded with bounded support, and the projection of $S_\eps \mu^{\rm FV, \eps}_{\eps^{-2}t}$
on the mark space $\R$ converges to the Lebesgue measure as $\eps\downarrow 0$, giving {\bf (A1)}.

{\bf (A2).} We will use the duality between interacting Fleming-Viot processes and coalescing random walks (see Theorem~\ref{T.dual}
and Corollary~\ref{C.tsdual}). Recall that $\Phi^\eps_t=\Phi(S_\eps \CX^{\rm FV, \eps}_{\eps^{-2}t})$. First, we bound
\begin{eqnarray}
\P(|\Phi^\eps_{\tau_\eps+\delta_\eps} - \Phi^\eps_{\tau_\eps}|>\eta) &\leq& \frac{1}{\eta^2} \E[ (\Phi^\eps_{\tau_\eps+\delta_\eps} - \Phi^\eps_{\tau_\eps})^2]
= \frac{1}{\eta^2} \E\big[ \E[(\Phi^\eps_{\tau_\eps+\delta_\eps} - \Phi^\eps_{\tau_\eps})^2 \,|\, \CX^{\rm FV, \eps}_{\eps^{-2}\tau_\eps}]\, \big] \nonumber \\
&\leq& \frac{2}{\eta^2} \E\big[{\rm Var}\big( \Phi^\eps_{\tau_\eps+\delta_\eps}\,|\,  \CX^{\rm FV, \eps}_{\eps^{-2}\tau_\eps}\big) \big] + \frac{2}{\eta^2}\E\big[ \big(\E[\Phi^\eps_{\tau_\eps+\delta_\eps}\,|\,  \CX^{\rm FV, \eps}_{\eps^{-2}\tau_\eps}]- \Phi^\eps_{\tau_\eps}\big)^2\big], \qquad \label{tight1}
\end{eqnarray}
where in the second inequality, we added and subtracted $\E[\Phi^\eps_{\tau_\eps+\delta_\eps}\,|\,  \CX^{\rm FV, \eps}_{\eps^{-2}\tau_\eps}]$
from $\Phi^\eps_{\tau_\eps+\delta_\eps} - \Phi^\eps_{\tau_\eps}$ and used $(a+b)^2 \leq 2a^2+2b^2$.
We treat the two terms on the r.h.s.\ separately.

{\em First term in \eqref{tight1}.}
{{}We bound this term} by bounding ${\rm Var}\big( \Phi^\eps_{\tau_\eps+\delta_\eps}\,|\,  \CX^{\rm FV, \eps}_{\eps^{-2}\tau_\eps}\big)$ uniformly in $\CX^{\rm FV, \eps}_{\eps^{-2}\tau_\eps}$. First note that $\CX^{\rm FV, \eps}$ is a strong Markov process by Theorem~\ref{T.wpp}. Therefore $\CX^{\rm FV, \eps}_{\eps^{-2}(\tau_\eps+\delta_\eps)}$ can be seen as the IFV genealogy process $\CX^{\rm FV, \eps}$ at time $\eps^{-2}\delta_\eps$ with initial condition $\CX^{\rm FV, \eps}_{\eps^{-2}\tau_\eps}$. In particular, it suffices to bound ${\rm Var}( \Phi^\eps_{\delta_\eps})$ uniformly in the initial condition $\CX^{\rm FV, \eps}_0$, which we can assume to be deterministic.

Let $\Phi=\Phi^{n,\phi,g}$, and denote $\ux:=(x_1, \ldots, x_n)\in \Z^n$, $\uy:=(y_1, \ldots, y_n)\in \Z^n$. Then by the definition of the scaling map $S_\eps$ in (\ref{Seps}), we have
\begin{eqnarray}
{\rm Var}( \Phi^\eps_{\delta_\eps}) &=& {\rm Var}\big( \Phi^{n,\phi,g}(S_\eps \CX^{\rm FV, \eps}_{\eps^{-2}\delta_\eps}) \big) = \E[(\Phi^\eps_{\delta_\eps})^2] - \E[\Phi^\eps_{\delta_\eps}]^2  \nonumber \\
&=& (\eps\sigma^{-1})^{2n}\sum_{\ux, \uy\in\Z^n} g(\eps \sigma^{-1}\ux) g(\eps \sigma^{-1}\uy)\ {\rm Cov}\big(\phi(\eps^2r^{\rm FV,\eps}_{\eps^{-2}\delta_\eps}(\ux)),\ \phi(\eps^2r^{\rm FV,\eps}_{\eps^{-2}\delta_\eps}(\uy)) \big), \label{tight2}
\end{eqnarray}
where $r^{\rm FV,\eps}_{\eps^{-2}\delta_\eps}(\ux)$ denotes the distance matrix $r^{\rm FV,\eps}_{\eps^{-2}\delta_\eps}(\xi(x_i), \xi(x_j))_{1\leq i<j\leq n}$ of $n$ individuals $\xi(x_1), \ldots, \xi(x_n)$ sampled independently from $\CX^{\rm FV, \eps}_{\eps^{-2}\delta_\eps}$ at positions $x_1, \ldots, x_n$ respectively. In order to evaluate the r.h.s.\ {{}of (\ref{tight2})} we represent the quantity using the duality in terms of a collection of coalescing random walks as follows.

Let $(X^{x_i}_t)_{1\leq i\leq n}$ and $(X^{y_i}_t)_{1\leq i\leq n}$ denote a family of rate $1$ continuous time random walks on $\Z$ with transition kernel $\bar a$ as in (\ref{abar}), and every pair of walks at the same site coalesce at rate $\gamma$. The coalescence gives a partition of the set of coalescing random walks at time $\eps^{-2}\delta_\eps$, and independently for each partition element, say at position $z\in\Z$, we sample an individual from $\CX^{\rm FV, \eps}_0$ at position $z$. Let $r^{\rm FV,\eps}_{0}(X^{\ux}_{\eps^{-2}\delta_\eps})$ denote the distance matrix of the collection of sampled individuals associated with the walks $X^{x_1}_{\eps^{-2}\delta_\eps},\ldots, X^{x_n}_{\eps^{-2}\delta_\eps}$ at time $\eps^{-2}\delta_\eps$, and let $r^{\rm FV,\eps}_{0}(X^{\uy}_{\eps^{-2}\delta_\eps})$ {be} defined similarly. We can further construct $(\wt X^{y_i}_t)_{1\leq i\leq n}$, a copy of $(X^{y_i}_t)_{1\leq i\leq n}$, which coincides with $(X^{y_i}_t)_{1\leq i\leq n}$ up to time $\eps^{-2}\delta_\eps$ on the event
$$
G_{\eps^{-2}\delta_\eps}(\ux,\uy):=\{ \mbox{none of } (X^{x_i})_{1\leq i\leq n} \mbox{ coalesces with any } (X^{y_i})_{1\leq i\leq n} \mbox{ before time } \eps^{-2}\delta_\eps\},
$$
such that $(\wt X^{y_i}_t)_{1\leq i\leq n}$ is independent of $(X^{x_i}_t)_{1\leq i\leq n}$. Let $r^{\rm FV,\eps}_{0}(\wt X^{\uy}_{\eps^{-2}\delta_\eps})$ be the associated distance matrix, which is independent of $r^{\rm FV,\eps}_{0}(X^{\ux}_{\eps^{-2}\delta_\eps})$.
By the duality relation (see Theorem~\ref{T.dual}), we have
\begin{equation}\label{tight2.5}
\begin{aligned}
& {\rm Cov}\big(\phi(\eps^2r^{\rm FV,\eps}_{\eps^{-2}\delta_\eps}(\ux)),\ \phi(\eps^2r^{\rm FV,\eps}_{\eps^{-2}\delta_\eps}(\uy)) \big)
=\  {\rm Cov}\big(\phi(\eps^2r^{\rm FV,\eps}_0(X^{\ux}_{\eps^{-2}\delta_\eps})),\ \phi(\eps^2r^{\rm FV,\eps}_0(X^{\uy}_{\eps^{-2}\delta_\eps})) \big) \\
=\ & \E\big[\phi(\eps^2r^{\rm FV,\eps}_0(X^{\ux}_{\eps^{-2}\delta_\eps})) \phi(\eps^2r^{\rm FV,\eps}_0(X^{\uy}_{\eps^{-2}\delta_\eps})) \big]
- \E\big[\phi(\eps^2r^{\rm FV,\eps}_0(X^{\ux}_{\eps^{-2}\delta_\eps}))\big] \E\big[\phi(\eps^2r^{\rm FV,\eps}_0(X^{\uy}_{\eps^{-2}\delta_\eps})) \big] \\
=\ & \E\big[\phi(\eps^2r^{\rm FV,\eps}_0(X^{\ux}_{\eps^{-2}\delta_\eps})) \big\{\phi(\eps^2r^{\rm FV,\eps}_0(X^{\uy}_{\eps^{-2}\delta_\eps}))- \phi(\eps^2r^{\rm FV,\eps}_0(\wt X^{\uy}_{\eps^{-2}\delta_\eps}))\big\} \big] \\
\leq\ & 2|\phi|_\infty^2 \P(G_{\eps^{-2}\delta_\eps}(\ux, \uy)^c) \leq 2|\phi|_\infty^2 \sum_{1\leq i,j\leq n} \P(\tau_{x_i, y_j} \leq \eps^{-2}\delta_\eps),
\end{aligned}
\end{equation}
where $\tau_{x_i,y_i}$ denotes the time it takes for the two walks $X^{x_i}_\cdot$ and $X^{y_j}_\cdot$ to meet. Note that this bound is uniform w.r.t.\ $\CX^{\rm FV, \eps}_0$. Substituting it into (\ref{tight2}) then gives
\begin{eqnarray}
{\rm Var}( \Phi^\eps_{\delta_\eps}) &\leq& |\phi|_\infty^2 \sum_{1\leq i,j\leq n} (\eps\sigma^{-1})^{2n}\sum_{\ux, \uy\in\Z^n} g(\eps \sigma^{-1}\ux) g(\eps \sigma^{-1}\uy)  \P(\tau_{x_i, y_j} \leq \eps^{-2}\delta_\eps) \nonumber \\
&=& |\phi|_\infty^2 \sum_{1\leq i,j\leq n} \sum_{\tilde\ux, \tilde\uy\in\eps\sigma^{-1}\Z^n} (\eps\sigma^{-1})^{2n}g(\tilde\ux) g(\tilde\uy)  \P(\tau_{\eps^{-1}\sigma \tilde x_i, \eps^{-1}\sigma y_j} \leq \eps^{-2}\delta_\eps) \label{tight3},
\end{eqnarray}
where $\tilde \ux := \eps\sigma^{-1}\ux$ and $\tilde \uy := \eps\sigma^{-1}\uy$.

We claim that the r.h.s.\ of (\ref{tight3}) tends to $0$ as $\eps\downarrow 0$. Indeed, the measure
$$
(\eps\sigma^{-1})^{2n} \sum_{1\leq i,j\leq n}\sum_{\tilde \ux, \tilde \uy\in\eps\sigma^{-1}\Z^n} \delta_{\tilde x_i} \delta_{\tilde y_j} g(\tilde \ux) g(\tilde\uy)
$$
converges weakly to the finite measure $g(\tilde\ux)g(\tilde\uy) \dd \tilde\ux \dd \tilde\uy$ on $\R^{2n}$ as $\eps\downarrow 0$. By Donsker's invariance principle and the fact that $\delta_\eps\to 0$ as $\eps\downarrow 0$, we note that for any $\lambda>0$,
\[
\P(\tau_{\eps^{-1}\sigma \tilde x_i, \eps^{-1}\sigma y_j} \leq \eps^{-2}\delta_\eps) \astoo{\eps} 0
\]
uniformly in $\tilde x_i$ and $\tilde y_j$ with $|\tilde x_i-\tilde y_j|> \lambda$. It follows that when restricted to $\tilde x_i$ and $\tilde y_j$ with $|\tilde x_i-\tilde y_j|> \lambda$, the inner sum in (\ref{tight3}) tends to $0$ as $\eps\downarrow 0$. On the other hand, when restricted to $\tilde x_i$ and $\tilde y_j$ with $|\tilde x_i-\tilde y_j|\leq \lambda$, the inner sum in (\ref{tight3}) can be bounded from above by replacing $\P(\cdot)$ with $1$, which then converges to the integral of the finite measure $g(\tilde\ux)g(\tilde\uy) \dd \tilde\ux \dd \tilde\uy$ over the subset of $\R^{2n}$ with $|\tilde x_i-\tilde y_j|\leq \lambda$, and can be made arbitrarily small by choosing $\lambda>0$ small.

This proves that ${\rm Var}( \Phi^\eps_{\delta_\eps})$ tends to $0$ uniformly in $\CX^{\rm FV, \eps}_0$ as $\eps\downarrow 0$, and hence the first term in (\ref{tight1}) tends to $0$ as $\eps\downarrow 0$.
\medskip

{\em Second term in \eqref{tight1}}. By the strong Markov property of $\CX^{\rm FV, \eps}$, it suffices to bound
$\big|\E[\Phi^\eps_{\delta_\eps}]- \Phi^\eps_0\big|$ uniformly in the (deterministic) initial condition $\CX^{\rm FV, \eps}_0$.

Let $(X^{x_i}_t)_{1\leq i\leq n}$,
$r^{\rm FV, \eps}_0(X^{\ux}_{\eps^{-2}\delta_\eps})$, and $r^{\rm FV, \eps}_t(\ux)$ be defined as before (\ref{tight2.5}). Let $(\wt X^{x_i}_t)_{1\leq i\leq n}$
be a collection of independent random walks, such that $(\wt X^{x_i}_t)_{1\leq i\leq n}$ coincides with $(X^{x_i}_t)_{1\leq i\leq n}$ up to time $\eps^{-2}\delta_\eps$ on the event
$$
G_{\eps^{-2}\delta_\eps}(\ux):=\{ \mbox{no coalescence has taken place among } (X^{x_i})_{1\leq i\leq n} \mbox{ before time } \eps^{-2}\delta_\eps\}.
$$
Then we have
\begin{eqnarray}\label{rf20}
&& \big|\E[\Phi^\eps_{\delta_\eps}]- \Phi^\eps_0\big| \nonumber \\
&\!\!\!\!\!\!=& (\eps\sigma^{-1})^{n} \Big|\sum_{\ux \in \Z^n} g(\eps\sigma^{-1}\ux) \E\big[\phi(\eps^2r^{\rm FV,\eps}_{\eps^{-2}\delta_\eps}(\ux))\big] - \sum_{\uy \in \Z^n} g(\eps\sigma^{-1}\uy) \E\big[\phi(\eps^2r^{\rm FV,\eps}_0(\uy))\big] \Big| \nonumber\\
&\!\!\!\!\!\!=&  (\eps\sigma^{-1})^{n} \Big|\sum_{\ux \in \Z^n} g(\eps\sigma^{-1}\ux) \E\big[\phi(\eps^2r^{\rm FV,\eps}_0(X^{\ux}_{\eps^{-2}\delta_\eps}))\big] - \sum_{\uy \in \Z^n} g(\eps\sigma^{-1}\uy) \E\big[\phi(\eps^2r^{\rm FV,\eps}_0(\uy))\big] \Big| \nonumber\\
&\!\!\!\!\!\!\leq& (\eps\sigma^{-1})^{n} \Big|\sum_{\ux \in \Z^n} g(\eps\sigma^{-1}\ux) \E\big[\phi(\eps^2r^{\rm FV,\eps}_0(\wt X^{\ux}_{\eps^{-2}\delta_\eps}))\big] - \sum_{\uy \in \Z^n} g(\eps\sigma^{-1}\uy) \E\big[\phi(\eps^2r^{\rm FV,\eps}_0(\uy))\big] \Big| \nonumber \\
&&\!\!\!\!\!\! + \ (\eps\sigma^{-1})^{n} \Big|\sum_{\ux \in \Z^n} g(\eps\sigma^{-1}\ux) \E\Big[\Big\{\phi(\eps^2r^{\rm FV,\eps}_0(X^{\ux}_{\eps^{-2}\delta_\eps}))- \phi(\eps^2r^{\rm FV,\eps}_0(\wt X^{\ux}_{\eps^{-2}\delta_\eps}))\Big\}1_{\{ G_{\eps^{-2}\delta_\eps}^c(\ux)  \}}\Big] \Big|, \nonumber
\end{eqnarray}
Note that the second term in the bound above tends to $0$ as $\eps\downarrow 0$ by the same argument as the one showing that the bound for
${\rm Var}( \Phi^\eps_{\delta_\eps})$ in (\ref{tight3}) tends to $0$ as $\eps\downarrow 0$. To bound the first term in the bound above,
we decompose according to the positions of the random walks and rewrite it as follows, where $p_t(x)$ denotes the transition probability kernel of $\wt X^0_t$:
\begin{equation}
\begin{aligned}
\label{e:rewrite}
& (\eps\sigma^{-1})^{n} \Big|\sum_{\ux \in \Z^n} g(\eps\sigma^{-1}\ux) \E\big[\phi(\eps^2r^{\rm FV,\eps}_0(\wt X^{\ux}_{\eps^{-2}\delta_\eps}))\big] - \sum_{\uy \in \Z^n} g(\eps\sigma^{-1}\uy) \E\big[\phi(\eps^2r^{\rm FV,\eps}_0(\uy))\big] \Big|
\\
  &= (\eps\sigma^{-1})^{n} \Big|\sum_{\uy \in \Z^n} \sum_{\ux\in\Z^n} \prod_{i=1}^n p_{\eps^{-2}\delta_\eps}(y_i-x_i) g(\eps\sigma^{-1}\ux) \E\big[\phi(\eps^2r^{\rm FV,\eps}_0(\uy)\big] - \sum_{\uy \in \Z^n} g(\eps\sigma^{-1}\uy) \E\big[\phi(\eps^2r^{\rm FV,\eps}_0(\uy))\big] \Big|
  \\
&\leq |\phi|_\infty \sum_{\uy\in\Z^n} (\eps\sigma^{-1})^{n} \Big|\sum_{\ux\in\Z^n} \prod_{i=1}^n p_{\eps^{-2}\delta_\eps}(y_i-x_i) \big(g(\eps\sigma^{-1}\ux)-g(\eps\sigma^{-1}\uy)\big) \Big|
  \\
&= |\phi|_\infty \sum_{\uy\in\Z^n} (\eps\sigma^{-1})^{n} \Big|\sum_{\ux\in\Z^n} \prod_{i=1}^n p_{\eps^{-2}\delta_\eps}(y_i-x_i) \Big\{\eps\sigma^{-1} \langle \ux-\uy, \nabla g (\eps\sigma^{-1}\uy)\rangle \\
& \hspace{6cm} +\ \frac{(\eps\sigma^{-1})^2}{2}\langle \ux-\uy, \nabla^2g(\eps\sigma^{-1}\uu(\ux,\uy)) (\ux-\uy)\rangle\Big\} \Big| \\
&\leq C |\phi|_\infty |\nabla^2 g|_\infty (\eps\sigma^{-1})^{n+2} \!\!\! \sum_{\ux, \uy\in\Z^n} \prod_{i=1}^n p_{\eps^{-2}\delta_\eps}(y_i-x_i)\big(1_{\{\eps\sigma^{-1}\ux \in {\rm supp}(g)\}}+ 1_{\{\eps\sigma^{-1}\uy \in {\rm supp}(g)\}}\big) \sum_{i=1}^n (y_i-x_i)^2 \\
&\leq 2nC |\phi|_\infty |\nabla^2 g|_\infty (\eps\sigma^{-1})^n \sum_{\uy\in\Z^n} 1_{\{\eps\sigma^{-1}\uy \in {\rm supp}(g)\}}  \ \delta_\eps,
\end{aligned}
\end{equation}
where $C$ is a constant depending only on $n$. In the derivation above, we Taylor expanded $g(\eps\sigma^{-1}\ux)$ around $\eps\sigma^{-1} \uy$ when either $\eps\sigma^{-1}\ux$ or $\eps\sigma^{-1}\uy$ is not in the support of $g$, $\nabla g$ and $\nabla^2 g$ denote the first and second derivatives of $g$, and $\uu(\ux, \uy)$ is some point on the line segment connecting $\ux$ and $\uy$. Lastly, we used the fact that $\sum_{z\in\Z} z p_t(z)=0$ and $\sum_{z\in\Z} z^2p_t(z)= t\sigma^2$. Since $g$ has bounded support, the bound we obtained above is bounded by $C'\delta_\eps$ for some $C'$ depending only on $n$, $\phi$ and $g$, and hence tends to $0$ as $\eps\downarrow 0$.

This verifies that the second term in (\ref{tight1}) also tends to $0$ as $\eps\downarrow 0$, which concludes the proof of {\bf (A2)}.
\qed
\medskip

We have verified {\bf (J2)} above and hence to
conclude the proof of tightness of $(S_\eps \CX^{\rm FV, \eps})_{\eps>0}$ as a family of random variables in the Skorohod space $D([0,\infty), \U^\R)$, it only remains to verify the compact containment condition {\bf (J1)}.
Some technical difficulties arise. Because the geographical space is unbounded, truncation in space is needed. We also need to control how the sizes of different families fluctuate over time, as well as how the population flux across the truncation boundaries affect the family sizes. Our strategy is to enlarge the mark space by assigning types to different families. Using a weaker convergence result, Theorem~\ref{T:conv2}, for measure-valued IFV processes with types (but no genealogies), we can control the evolution of family sizes as well as their dispersion in space, which can then be strengthened to control the genealogical structure of the population. As we will point out later, Theorem~\ref{T:conv2} can be proved by adapting what we have done so far, because condition {\bf (J1)} is trivial in that context. Therefore invoking Theorem~\ref{T:conv2} to prove {\bf (J1)} for the genealogy processes is justified.

\medskip

\noindent
{\bf Proof of (J1).} As noted in Remark~\ref{R:GromovHash}, we can regard $\U^\R$ as a subset of $(\U^\R_f)^{\N}$, endowed with the product $\R$-marked Gromov-weak topology. {{}Since the product of compact sets gives a compact set in the product space $(\U^\R_f)^{\N}$, }
to prove {\bf (J1)}, it suffices to show that for each $k\in\N$, the restriction of $(S_\eps \CX^{\rm FV, \eps}_{\eps^{-2}t})_{t\geq 0, \eps>0}$ to the subset of marks $(-k,k) \subset \R$, i.e.,
$$
S^{(k)}_\eps \CX^{\rm FV, \eps}_{\eps^{-2}t} := \overline{(X^{\rm FV,\eps}_{\eps^{-2}t}, S_\eps r^{\rm FV, \eps}_{\eps^{-2}t}, 1_{\{|v|< k\}}S_\eps \mu^{\rm FV,\eps}_{\eps^{-2}t}(\dd x \dd v))}, \qquad t\geq 0, \eps>0,
$$
satisfies the compact containment condition (\ref{compactcont}). More precisely, it suffices to show that for each $T>0$ (which we will assume to be $1$ for simplicity), and for each $\delta>0$, there exists a compact $K_\delta \subset \U^\R_f$, such that for all $\eps>0$ sufficiently small,
\begin{equation} \label{compactcont2}
\P\big(S^{(k)}_\eps \CX^{\rm FV, \eps}_{\eps^{-2}t} \in K_{\delta} \ \forall \, 0 \leq t \leq 1 \big)  \geq 1 - \delta,
\end{equation}
where $k\in\N$ will be fixed for the rest of the proof.

We will construct $K^1_\delta, K^2_\delta, K^3_\delta \subset \U^\R_f$, which satisfy respectively conditions (i)--(iii) in Theorem~\ref{T:mmmrc} for the relative compactness of subsets of $\U^\R_f$. We can then take $K_\delta:= \overline{ K^1_\delta \cap K^2_\delta \cap K^3_\delta}$, which is a compact subset of $\U^\R_f$. To prove (\ref{compactcont2}), it then suffices to prove the same inequality
but with $K_\delta$ replaced by $K^i_\delta$ for each $1\leq i\leq 3$, which we do in {\bf (1)--(3)} below.

Later when we construct $K_\delta^2$ and $K_\delta^3$, we will keep track of the mass of different individuals having some specified
properties. The way to do this is to introduce additional marks. We will enlarge the mark space from $\R$ to $\R\times [0,1]$,
where $[0,1]$ is the space of the additional types, and $\CX^{\rm FV, \eps}_t$ and $\CX^{\rm CS}_t$ then become random variables taking values in the space
$\U^{\R\times [0,1]}$. For each $\overline{(X,r,\mu)}\in \U^{\R\times [0,1]}$, $\mu(\dd x\dd v\dd \tau)$ is then a measure on
$X\times \R\times [0,1]$. The types of individuals in $\CX^{\rm FV, \eps}_0$ and $\CX^{\rm CS}_0$ will be assigned later as we see fit.
\sm

{\bf (1)}
First let $K_\delta^{1}$ be the subset of $\U^\R_f$, such that for each $\overline{(X, r, \mu)} \in K^1_\delta$, $\mu(X\times \cdot)$ is supported on $[-k,k]$, with total mass bounded by $4k$. Since the family of measures on $[-k, k]$, with total mass bounded by $4k$, is relatively compact w.r.t.\ the weak topology, $K_\delta^1$ satisfies condition (i) in Theorem~\ref{T:mmmrc}. We further note that a.s., $S^{(k)}_\eps \CX^{\rm FV, \eps}_{\eps^{-2}t} \in K^1_\delta$ for all $t\geq 0$, and hence (\ref{compactcont2}) holds with $K_\delta$ replaced by $K^1_\delta$.

{\bf (2)}
For each $n\in\N$, we will find below $L(n)$ such that if $K_\delta^{2,n}$ denotes the subset of $\U^\R_f$ with
\begin{equation}\label{K2deltan}
\iint_{(X\times \R)^2} 1_{\{ r(x,y)> L(n)  \}} \mu(\dd x\dd u) \mu(\dd y \dd v) <  \frac{1}{n} \qquad \mbox{for each } \overline{(X,r,\mu)}\in K_\delta^{2,n},
\end{equation}
then uniformly in $\eps>0$ sufficiently small, we {{}will} have
\begin{equation}\label{comcont00}
\P\big( S^{(k)}_\eps \CX^{\rm FV, \eps}_{\eps^{-2}t} \notin K^{2,n}_{\delta} \ \mbox{for some } 0 \leq t \leq 1 \big) \leq \frac{\delta}{2^n}=:\delta_n.
\end{equation}
We can then take $K^2_\delta := \cap_{n\in\N} K^{2,n}_\delta$, which clearly satisfies condition (ii) in Theorem~\ref{T:mmmrc}, while (\ref{comcont00}) implies that
(\ref{compactcont2}) holds with $K_\delta$ replaced by $K_\delta^2$.

In order to find $L(n)$, we proceed in two steps. First we find an analogue of $L(n)$ for the limiting CSSM genealogies, and then in a second step, we use the convergence of the measure-valued IFV to obtain $L(n)$.

Fix $n\in\N$. To find $L(n)$ such that (\ref{comcont00}) holds, we first prove an analogue of (\ref{comcont00}) for the continuum limit
 $\CX^{\rm CS}$ by utilizing the types of the individuals. Given $\eta \geq 0$ and $\gamma \in [0,1]$, let
\begin{equation}\label{Geta}
G_{\eta, \gamma}:= \big\{\overline{(X, r, \mu)} \in \U^{\R\times [0,1]}_f: \mu(X\times [-k, k] \times [0,\gamma]) \leq \eta  \big\}.
\end{equation}
We claim that we can find $A$ sufficiently large, such that if all individuals in $\CX^{\rm CS}_0$ with spatial mark outside $[-A, A]$ are assigned type $0$, and all other individuals are assigned type $1$, then
\begin{equation}\label{prehatCXt}
\P\big( \CX^{\rm CS}_t \in G_{0, 0} \ \forall\, 0\leq t\leq 1 \big) \geq 1-\frac{\delta_n}{4}.
\end{equation}
In other words, with probability at least $1-\delta_n/4$, the following event occurs: for all $0\leq t\leq 1$, no individual in $\CX^{\rm CS}_t$ with spatial mark in $[-k,k]$ can trace its genealogy back to some individual at time $0$ with spatial mark outside $[-A, A]$.  This will allow us to restrict our attention to descendants of the population in $[-A, A]$ at time $0$.

Indeed, by the construction of the CSSM (Section~\ref{S:IntroCSSM}) using the Brownian web, the measure-valued process
$\widehat \CX^{\rm CS}_t$, which is the measure $\mu^{\rm CS}_t$ projected on the geographic and type space, is given by
\begin{equation}\label{hatCXt}
\widehat \CX^{\rm CS}_t (\dd v\dd \tau) = \delta_0(\dd \tau) \big(1_{\{v>f_{(A,0)}(t)\}} + 1_{\{v<f_{(-A,0)}(t)\}}\big) \dd v + \delta_1(\dd \tau) 1_{\{f_{(-A,0)}(t)\leq v\leq f_{(A,0)}(t)\}} \dd v,
\end{equation}
where $f_{(\pm A,0)}(\cdot)$ are the two coalescing Brownian motions in the Brownian web $\CW$, starting respectively at $(\pm A, 0)$. The event in (\ref{prehatCXt}) occurs if $f_{(\pm A, 0)}$ do {not} enter $[-k, k]$ before time $1$, the probability of which can be made arbitrarily close to $1$ by choosing $A$ large.

Let $A>0$ be chosen such that (\ref{prehatCXt}) holds. Since $\CX^{\rm CS}_0=\overline{(X^{\rm CS}_0, r^{\rm CS}_0,
\mu^{\rm CS}_0)} \in \U^{\R\times [0,1]}$, the population in $X^{\rm CS}_0$ with spatial marks in $[-A, A]$ can be partitioned
into a countable collection of disjoint balls of radius $\frac{1}{n}$. Let us label these balls by $B_{n,1}, B_{n,2}, \ldots \subset X^{\rm CS}_0$
in decreasing order of their total measure, and all individuals in $B_{i}$ is given type $\frac{1}{i} \in [0,1]$,
while all individuals with spatial mark outside $[-A, A]$ are given type $0$. Note that by choosing $M$ large,
	the total measure of individuals in $\CX^{\rm CS}_0$ with types in $(0, \frac{1}{M}]$ can be made arbitrarily small.
	It then follows by the Feller continuity of the measure-valued process $\widehat \CX^{\rm CS}$ stated in Remark~\ref{R:reg},
	that as $M\to\infty$, the total measure of individuals in $\CX^{\rm CS}_t$, $t\in [0, 1]$, with spatial marks in $[-k,k]$
	and types in $(0, \frac{1}{M}]$, converge in probability to the stochastic process which is identically $0$ on the time interval
	$[0,1]$. Combined with (\ref{prehatCXt}), this implies that
we can choose $M$ large, such that

\begin{equation}\label{comcont1}
\P\big(\CX^{\rm CS}_t \notin G_{\frac{1}{4kn}, \frac{1}{M}} \ \mbox{for some } 0\leq t\leq 1 \big) \leq \frac{\delta_n}{2}.
\end{equation}
Let $\widetilde L(n)$ denote the maximal distance between the balls $B_{n,1}, \ldots, B_{n, M}$. Note that on the event that {$\{\CX^{\rm CS}_t \in G_{\frac{1}{4kn}, \frac{1}{M}}$ for all $0\leq t\leq 1\}$}, for each $t\in [0,1]$, the set of individuals in $\CX^{\rm CS}_t$ with spatial marks in $[-k, k]$ and types in $[0, \frac{1}{M}]$ have total measure at most $\frac{1}{4kn}$; while the rest of the population with spatial marks in $[-k,k]$ trace their genealogies back at time $0$ to $B_{n,1},\ldots, B_{n, M}$, and their mutual distance is bounded by $\widetilde L(n)+2$. Letting $L(n) = \widetilde L(n)+3$ in (\ref{K2deltan}) then readily implies that $\overline{(X^{\rm CS}_t, r^{\rm CS}_t, 1_{\{ |v|\leq k\}} \mu^{\rm CS}_t (\dd x\dd v))}\in K_\delta^{2,n}$ on the event $\CX^{\rm CS}_t \in G_{\frac{1}{4kn}, \frac{1}{M}}$, and (\ref{comcont1}) implies that
\begin{equation}\label{comcont01}
\P\big( \overline{(X^{\rm CS}_t, r^{\rm CS}_t, 1_{\{ |v|\leq k\}} \mu^{\rm CS}_t (\dd x\dd v))} \notin K^{2,n}_{\delta} \ \mbox{for some } 0 \leq t \leq 1 \big) \leq \frac{\delta_n}{2},
\end{equation}
which is the analogue of (\ref{comcont00}) for $\CX^{\rm CS}$.

Next we turn to $\CX^{\mathrm{FV},\ve}$ and we will deduce (\ref{comcont00}) from (\ref{comcont01}) by exploiting Theorem~\ref{T:conv2}
on the convergence of the measure-valued processes $S_\eps \widehat \CX^{\rm FV, \eps}$ to $\widehat \CX^{\rm CS}$. Let $A$ be the same as above.
Let individuals in $\CX^{\rm CS}_0$ be assigned types in $\{0\}\cup\{\frac{1}{i}: i \in\N\}$ as before (\ref{comcont1}),
where $\widetilde L(n)$ denotes the maximal distance between any pair of individuals in $\CX^{\rm CS}_0$ with types in $[\frac{1}{M},1]$.
The assumption that $S^\eps \CX^{\rm FV, \eps}_0$ converges to $\CX^{\rm CS}_0$ in $\U^\R$ as $\eps\downarrow 0$ allows a coupling
between $S_\eps X^{\rm FV,\eps}_0$ and $X^{\rm CS}_0$, such that for most of the individuals in $S_\eps X^{\rm FV,\eps}_0$ with
spatial marks in $[-A,A]$, their genealogical distances and spatial marks are close to their counterparts in $X^{\rm CS}_0$.
We can then assign types to individuals in $S_\eps X^{\rm FV,\eps}_0$ in such a way that: individuals with spatial marks outside $[-A,A]$
are assigned type $0$ while those with spatial marks in $[-A, A]$ are assigned types in $\{\frac{1}{i}: i\in\N\}$.
The distance between any pair of individuals with types in $[\frac{1}{M}, 1]$ is bounded by $\widetilde L(n)+1$,
and the measure on geographic and type space, $S_\eps \widehat\CX^{\rm FV, \eps}_0$, converges vaguely to $\widehat \CX^{\rm CS}_0$
as $\eps\downarrow 0$. By Theorem~\ref{T:conv2}, $(S_\eps\widehat\CX^{\rm FV, \eps}_{\eps^{-2}t})_{0\leq t\leq 1}$ converges weakly to
$(\widehat\CX^{\rm CS}_t)_{0\leq t\leq 1}$ as random variables in $C([0,1], \CM(\R\times [0,1]))$. Applying this weak convergence result to
(\ref{comcont1}) then implies that for $\eps>0$ sufficiently small, we have the following analogue of (\ref{comcont1}):
\begin{equation}\label{comcont0}
\P\big( S_\eps \CX^{\rm FV, \eps}_{\eps^{-2}t} \notin G_{\frac{1}{2kn}, \frac{1}{M}} \ \mbox{for some } 0\leq t\leq 1 \big) \leq \delta_n.
\end{equation}

Note that on the event $S_\eps \CX^{\rm FV, \eps}_{\eps^{-2}t} \in G_{\frac{1}{2kn}, \frac{1}{M}}$ for all $0\leq t\leq 1$, for each $0\leq t\leq 1$, the set of individuals in $S_\eps \CX^{\rm FV, \eps}_{\eps^{-2}t}$ with spatial marks in $[-k,k]$ and types in $[0, \frac{1}{M}]$ have total measure at most $\frac{1}{2kn}$; while the rest of the population with spatial marks in $[-k, k]$ trace their genealogies back to an individual at time $0$ with type in $[\frac{1}{M}, 1]$, and hence their pairwise distance is bounded by $L(n)=\widetilde L(n)+3$.
It follows that $S_\eps^{(k)}\CX^{\rm FV, \eps}_{\eps^{-2}t} \in K^{2,n}_\delta$ on the event $S_\eps \CX^{\rm FV, \eps}_{\eps^{-2}t} \in G_{\frac{1}{2kn}, \frac{1}{M}}$. Together with (\ref{comcont0}), this implies
(\ref{comcont00}).
\medskip

{\bf (3)}
Our procedure for constructing $K^3_\delta$ is similar to that of $K^2_\delta$. For each $n\in\N$, we will find $M=M(n)$
such that if $K_\delta^{3,n}$ denotes the subset of $\U^\R_f$ with the property that for each
$\overline{(X,r,\mu)}\in K_\delta^{3,n}$, we can find $M$ balls of radius $\frac{1}{n}$ in $X$, say $B_1, \ldots, B_M$ with
$B:=\cup_{i=1}^M B_i$, such that $\mu(X\backslash{B} \times \R) < \frac{1}{n}$, then uniformly in $\eps>0$ sufficiently small, we have
\begin{equation}\label{comcont02}
\P\big( S^{(k)}_\eps \CX^{\rm FV, \eps}_{\eps^{-2}t} \notin K^{3,n}_{\delta} \ \mbox{for some } 0 \leq t \leq 1 \big) \leq \frac{\delta}{2^n}=:\delta_n.
\end{equation}
We can then take $K^3_\delta := \cap_{n\in\N} K^{2,n}_\delta$, which clearly satisfies condition (iii) in Theorem~\ref{T:mmmrc}, while (\ref{comcont02}) implies that (\ref{compactcont2}) holds with $K_\delta$ replaced by $K_\delta^3$.

To find $M(n)$ such that (\ref{comcont02}) holds, we partition the time interval $[0,1]$ into $[\frac{i-1}{2n}, \frac{i}{2n}]$ for $1\leq i\leq 2n$. It then suffices to show that for each $1\leq i\leq 2n$, we can find $M_i(n)$ such that if $K^{3,n}_\delta$ is defined using $M_i$, then uniformly in $\eps>0$ small,
\begin{equation}\label{comcont03}
\P\big( S^{(k)}_\eps \CX^{\rm FV, \eps}_{\eps^{-2}t} \notin K^{3,n}_{\delta} \ \mbox{for some } \frac{i-1}{2n} \leq t \leq \frac{i}{2n} \big) \leq \frac{\delta_n}{2n}
= \frac{\delta}{2n 2^n}.
\end{equation}

Again we first determine $M(n)$ for CSSM and then use the convergence of measure-valued IFV to measure-valued CSSM.
We now prove an analogue of (\ref{comcont03}) for $\CX^{\rm CS}$. Since $\CX^{\rm CS}_{\frac{i-1}{2n}}\in \U^{\R\times [0,1]}$ almost surely, we can condition on its realization and partition the population in $\CX^{\rm CS}_{\frac{i-1}{2n}}$ with spatial marks in $[-A, A]$ into disjoint balls of radius $\frac{1}{3n}$, $B_1, B_2, \ldots$, as we did in the argument leading to (\ref{comcont1}). Repeating the same argument there and assigning type $\frac{1}{j}$ to individuals in $B_j$, we readily obtain the following analogue of (\ref{comcont1}): we can choose $M_i$ large enough such that
\begin{equation}\label{comcont11}
\P\big(\CX^{\rm CS}_t \notin G_{\frac{1}{2n}, \frac{1}{M_i}} \ \mbox{for some } \frac{i-1}{2n}\leq t\leq \frac{i}{2n} \big) \leq \frac{\delta_n}{4n}.
\end{equation}
Note that on the event $\CX^{\rm CS}_t \in G_{\frac{1}{2n}, \frac{1}{M_i}}$ for all $\frac{i-1}{2n}\leq t\leq \frac{i}{2n}$, for each $t\in [\frac{i-1}{2n}, \frac{i}{2n}]$, the set of individuals in $\CX^{\rm CS}_t$ with spatial marks in $[-k, k]$ and types in $[0, \frac{1}{M_i}]$ have total measure at most $\frac{1}{2n}$; while the rest of the individuals in $\CX^{\rm CS}_t$ with spatial marks in $[-k, k]$ trace their genealogies back to an individual at time $\frac{i-1}{2n}$ in $B_1\cup B_2\cdots \cup B_{M_i}$, and hence they are contained in $M_i$ balls of radius $\frac{1}{3n}+t\leq \frac{5}{6n}$. Therefore $\overline{(X^{\rm CS}_t, r^{\rm CS}_t, 1_{\{ |v|\leq k\}} \mu^{\rm CS}_t (\dd x\dd v))} \in K^{3, n}_\delta$, and the analogue of (\ref{comcont03}) holds for $\CX^{\rm CS}$.

To establish (\ref{comcont03}) uniformly in small $\eps>0$, we again apply the convergence result in Theorem~\ref{T:conv2}. Note that by the f.d.d.\ convergence established in Section~\ref{ss.convfin}, $S_\eps \CX^{\rm FV, \eps}_{\eps^{-2}\frac{i-1}{2n}}$ converges in distribution to $\CX^{\rm CS}_{\frac{i-1}{2n}}$ as $\eps\downarrow 0$. Following the same argument as those leading to (\ref{comcont0}), we can then assign types to individuals in $S_\eps \CX^{\rm FV, \eps}_{\eps^{-2}\frac{i-1}{2n}}$ such that the associated measure-valued process, $(S_\eps\widehat \CX^{\rm FV, \eps}_{\eps^{-2}t})_{\frac{i-1}{2n}\leq t\leq \frac{i}{2n}}$ converges weakly to $(\widehat \CX^{\rm CS}_t)_{\frac{i-1}{2n}\leq t\leq \frac{i}{2n}}$, and individuals in $S_\eps \CX^{\rm FV, \eps}_{\eps^{-2}\frac{i-1}{2n}}$ with spatial marks in $[-A, A]$ and types in $[\frac{1}{M_i}, 1]$ are contained in $M_i$ balls with radius at most $\frac{1}{2n}$. Applying this convergence result to (\ref{comcont11}) then implies that for $\eps>0$ sufficiently small, we have the following analogue of (\ref{comcont11}):
\begin{equation}\label{comcont12}
\P\big(S_\eps\CX^{\rm FV, \eps}_t \notin G_{\frac{1}{n}, \frac{1}{M_i}} \ \mbox{for some } \frac{i-1}{2n}\leq t\leq \frac{i}{2n} \big) \leq \frac{\delta_n}{2n}.
\end{equation}
This is then easily seen to imply (\ref{comcont03}).

Combining parts (1)-(3) concludes the proof of (\ref{compactcont2}) and hence establishes the
compact containment condition {\bf (J1)}.
\qed

\subsection{Proof of convergence of rescaled IFV processes (Theorem~\ref{T:conv2})}
\label{ss.conv2}

In~\cite[Theorem~1.1]{AS11}, a convergence result similar to Theorem~\ref{T:conv2} was proved for the voter model on $\Z$, where the type space consists of only $\{0,1\}$, and a special initial condition was considered, where the population to the left of the origin all have type $1$, and the rest of the population have type $0$. The proof consists of two parts: proof of tightness, and convergence of finite-dimensional distribution.

In~\cite{AS11}, the proof of tightness for the voter model does not depend on the initial condition, and is based on the verification of Jakubowski's criterion and Aldous' criterion as we have done in Sections~\ref{ss.convfin} and \ref{ss.tight}
for the genealogical process.
Because the IFV process ignores the genealogical distances, the verification of the compact containment condition {\bf (J1)}
in Jakubowski's criterion is trivial, as in the case for the voter model. Using the duality between the IFV process and coalescing random walks
 with {\em delayed} coalescence,
 Aldous' criterion on the tightness of evaluations can be verified by exactly the same calculations as that for the voter model in~\cite{AS11},
  which uses the duality between the voter model and coalescing random walks with {\em instantaneous} coalescence.
  Recall here that in the rescaling the difference between instantaneous and delayed coalescence disappears because
  of recurrence of the difference walk.
   Lastly, the convergence of the finite-dimensional distribution for rescaled IFV process follows the same calculations as in
   Section~\ref{ss.convfin}, where we can simply enlarge the mark space to $\R\times [0,1]$ and suppress the genealogical
   distances by choosing $\phi_i\equiv 1$.
\qed

\section{Martingale Problem for CSSM Genealogy Processes}\label{S:CSSMmart}

In this section, we show that the CSSM genealogy process constructed in Section \ref{S:IntroCSSM} solves the martingale problem formulated in Theorem~\ref{T.mp}. We will first identify the generator action on regular test functions evaluated at regular states, and then extend it to more general test functions and verify the martingale property. Complications arise mainly from the singular nature of the resampling component of the generator, which are only well-defined a priori on regular test functions evaluated at regular states. Fortunately by Proposition~\ref{P:reg}, the CSSM genealogy process enters these regular states as soon as $t>0$, even though the initial state may not be regular.

\subsection{Generator Action on Regular Test Functions}\label{ss:generator1}
In this section, we identify the generator of the CSSM genealogy process $\CX^{\mathrm{CS}}$, acting on $\Phi \in \Pi^{1,2}_{\rm r}$ and evaluated at $\CX^{\rm CS}_t\in \U^\R_{\rm r}$, where $\Pi^{1,2}_{\rm r}$ and $\U^\R_{\rm r}$ are introduced in Section \ref{ss.mpgen}. The advantage of working with such regular $\Phi$ and $\CX^{\rm CS}_t$ is that, the relevant resamplings only occur at the boundary points of balls of radius at least $\delta$, which is a locally finite set. In Section~\ref{ss:generator2}, we will extend it to the case $\Phi\in \Pi^{1,2}$, where we need to consider the boundaries of all balls in $\CX^{\rm CS}_t$, which is a locally infinite set.

\bp\label{P:genaction}{\bf [Generator action on regular test functions]}\\
Let $\CX^{\rm CS}:=(\CX^{\rm CS}_t)_{t\geq 0}$ be the CSSM genealogy process with $\CX^{\rm CS}_0\in \U^\R_{\rm r}$. Let $\Phi=\Phi^{n, \phi, g}\in \Pi^{1,2}_{\rm r}$, defined as in Definition~\ref{D:PI12r}. Let $L^{\rm CS}=L^{\rm CS}_{\rm d} + L^{\rm CS}_{\rm a} + L^{\rm CS}_{\rm r}$ be defined as in \eqref{LCS}--\eqref{LCSr}.
Then we have
\be{generator1}
\lim_{t\downarrow 0} \frac{\E[\Phi(\CX^{\rm CS}_t)] - \Phi(\CX^{\rm CS}_0)}{t} = L^{\rm CS}\Phi({\CX^{\rm CS}_0}).
\ee
Furthermore,
\be{inteqn1}
\E[\Phi(\CX^{\rm CS}_t)] = \Phi(\CX^{\rm CS}_{0}) +\int_0^t \E[L^{\rm CS} \Phi(\CX^{\rm CS}_s)]\dd s,
\ee
where $\E[L^{\rm CS}\Phi(\CX^{\rm CS}_t)]$ is continuous in $t\geq 0$.
\ep
The proof is fairly long and technical and will be broken into parts, with (\ref{generator1}) and (\ref{inteqn1}) proved respectively in Sections \ref{sss.prooft9} and \ref{sss.proofPgena}.

\subsubsection{Proof of \eqref{generator1} in Proposition \ref{P:genaction}}
\label{sss.prooft9}

For each $t\geq 0$, denote $\CX^{\rm CS}_t =\overline{(X^{\rm CS}_t, r^{\rm CS}_t, \mu^{\rm CS}_t)}$. Let $L>0$ be chosen such that the support of $g(\ux)$ is contained in $(-L,L)^n$. Let $\delta>0$ be determined by $\phi$ as in (\ref{gs13}), so that $\phi((r_{i,j})_{1\leq i<j\leq n})$ is constant when any coordinate $r_{i,j}$ varies on the interval $r_{i,j}\in [0,\delta]$.

We proceed in five steps, first giving a suitable representation {{}of} $\Phi (X^{\mathrm{CS}}_0)$ and $\Phi (X^{\mathrm{CS}}_t)$,
and then calculating actions that lead to the different parts of the generator.
\medskip

\noindent
{{\bf Step~1 (Representation of $\Phi (X^{\mathrm{CS}}_0)$). }} We derive here a representation  {{}of} $\Phi (X^{\mathrm{CS}}_0)$ by partitioning $X^{\rm CS}_0$ into disjoint balls of radius {{}at least} $\delta/4$ and utilizing the fact that $\Phi\in \Pi^{1,2}_{\rm r}$.

Since $\CX^{\rm CS}_0\in \U^\R_{\rm r}$, by Remark \ref{R:Ur} and (\ref{XCSsplit}), we can identify $X^{\rm CS}_t$ for each $t\geq 0$ with
\be{u1}
\R^*_t:= \bigcup_{x\in E_t}\{x^+, x^-\} \bigcup \, (\R\backslash E_t),
\ee
where $x^+$ and $x^-$ are duplicates of the point $x$ in
\be{Elt22}
E_t=\cup_{l>0} E^{l}_t,
\ee
where $E^l_t$ is the set of points in $\R$ that lie at the boundary of two disjoint balls {{}of radius  $l$} in $X^{\rm CS}_t$, which is consistent with the definition in (\ref{Elt}).

Denote
\be{EdL1}
\{y_1< y_2< \cdots <y_m\} :=E^{\frac{\delta}{4}}_0 \cap [-2L, 2L],
\ee
and let $y_0$, resp.\ $y_{m+1}$, be the point in $E^{\frac{\delta}{4}}_0$ adjacent to $y_1$, resp.\ $y_m$. Note that the intervals $[y_i^+, y_{i+1}^-]:=(y_i, y_{i+1})\cup\{y_i^+, y_{i+1}^-\}$ form disjoint open balls of radius $\frac{\delta}{4}$ in $X^{\rm CS}_0$. Therefore for all $x_1\in [y_{i}^+, y_{i+1}^-]$ and
$x_2 \in [y_{j}^+, y_{j+1}^-]$ with $i\neq j$, $d^{ij}_0:=r^{\rm CS}_0(x_1, x_2)$ is constant, and
$d_0 := (d_0^{ij})_{0\leq i\leq j\leq m}$ forms a distance matrix, where we set $d_0^{ii}:=0$ for $0\leq i\leq m$. We can then write
\be{CU0}
\Phi(\CX^{\rm CS}_0) = \Phi^{n, \phi, g}(\CX^{\rm CS}_0)= \sum_{\uk \in \{0,\cdots, m\}^n} \phi(d_0^{\uk\times\uk}) \ G(\uy, \uk),
\ee
where $\uy=(y_0,\cdots, y_{m+1})$, $d_0^{\uk\times \uk}$ is the distance matrix with
$d_{0,ij}^{\uk\times \uk}=d_0^{k_ik_j}$, {{}$1\leq i\le j\leq n$}, and
\be{intG}
G(\uy, \uk) =\idotsint\limits_{y_{k_i}< x_i <y_{k_i+1} \atop 1\leq i\leq n} g(\ux) \, d\ux.
\ee
\medskip

\noindent
{{\bf Step~2 (Representation of $\Phi (X^{\mathrm{CS}}_t)$). }}
We next write $\Phi(\CX^{\rm CS}_t)$, for $t>0$, in terms of coalescing Brownian motions
running forward in time (in contrast to the spatial genealogies which run backward in time), which determine
the evolution of boundaries between disjoint balls in $\CX^{\rm CS}_t$.

Let $\{y_i+B_i\}_{0\leq i\leq m+1}$ be independent  Brownian motions starting from
$\{y_i\}_{0\leq i\leq m+1}$, from which we construct a family of coalescing Brownian motions
$\{y_i+\wt B_i\}_{0\leq i\leq m+1}$. Namely, let $y_0+\wt B_0 := y_0+B_0$ for all time, and let
$y_1+\wt B_1(s) := y_1+B_1(s)$ until the first time $y_1+B_1(s)$ hits $y_0+\wt B_0$. From
that time onward, define $y_1+\wt B_1$ to coincide with $y_0+\wt B_0$. In the same way, we successively
define $\{y_i+\wt B_i\}_{2\leq i\leq m+1}$ from $\{y_i+B_i\}_{2\leq i\leq m+1}$ by adding one path at a time.
Without loss of generality, we may assume that $y_i+\wt B_i$ is the a.s.\ unique path in the Brownian web $\CW$
starting from $(y_i,0)$.

To write $\Phi(\CX^{\rm CS}_t)$ in terms of the forward coalescing Brownian motions $\wt B$, we observe that
\be{Uspartition}
E^\delta_s \cap [y_0+ \wt B_0(s), y_{m+1}+ \wt B_{m+1}(s)] \subset \bigcup_{i=0}^{m+1}\{ y_i+ \wt B_i(s)\}
\quad \forall\ 0\leq s< \min\Big\{\frac{\delta}{4}, t\Big\},
\ee
since by our construction of $\CX^{\rm CS}_s$ in Section~\ref{S:IntroCSSM}, $(y_i+\wt B_i(s))_{0\leq i\leq m+1}$ are boundaries
of disjoint balls in $\CX^{\rm CS}_s$, and  $[(y_i+ \wt B_i(s))^+, (y_{i+1}+ \wt B_{i+1}(s))^-]$, $0\leq i\leq m$,
consists of either empty sets if $y_i+\wt B_i(s)=y_{i+1}+\wt B_{i+1}(s)$, or disjoint open balls of radius
$\frac{\delta}{4}+s$ in $X^{\rm CS}_s$.

Let $\{{\rm meet}\}_t$ denote the event that either $y_0+ B_0(s)$ reaches level $-L$
before time $t$, or $y_{m+1}+ B_{m+1}(s)$ reaches level $L$ before time $t$, or one of the pair
$(y_i+ B_i(s), y_{i+1}+ B_{i+1}(s))$ meet before time $t$. On the complementary event
$\{{\rm meet}\}_t^c$, $\wt B_i=B_i$ for all $i$. Therefore, for all $0\leq t < \frac{\delta}{4}$, we can write
\be{rs10}
\Phi(\CX^{\rm CS}_t) = 1_{\{{\rm meet}\}_t} \Phi(\CX^{\rm CS}_t) + 1_{\{{\rm meet}\}^c_t} \wt \Phi(\CX^{\rm CS}_t)
\ee
with
\be{phitilde}
\wt \Phi(\CX^{\rm CS}_t) = \sum_{\uk \in \{0,\cdots, m\}^n} \phi(d_t^{\uk\times\uk})
\ G(\uy+\uB(t), \uk),
\ee
where $d_t^{\uk\times \uk}$ is the distance matrix with $d_{t,ij}^{\uk\times \uk}= d^{k_ik_j}_0 + 2t(1-\delta_{k_i,k_j})$ for $1\leq i \leq j\leq n$, and $\uB(t)=(B_0(t),\cdots,B_{m+1}(t))$.

Since
\be{meeterror}
|\Phi(\CX^{\rm CS}_t) - \wt\Phi(\CX^{\rm CS}_t)| = 1_{\{{\rm meet}\}_t} |\Phi(\CX^{\rm CS}_t)-\wt\Phi(\CX^{\rm CS}_t)|
\leq C_{\Phi} 1_{\{{\rm meet}\}_t}
\ee
for some $C_{\Phi}$ depending on $\Phi$, and the probability of the event $\{{\rm meet}\}_t$ decays exponentially in $t^{-1}$ by elementary estimates for Brownian motions, we have
$\lim_{t\to 0} \frac{\E[\Phi(\CX^{\rm CS}_t)] - \E[\wt\Phi(\CX^{\rm CS}_t)]}{t}=0$.
Thus, we may replace $\Phi(\CX^{\rm CS}_t)$ by $\wt\Phi(\CX^{\rm CS}_t)$ up to error $o(t)$ as $t \to 0$.

By (\ref{CU0}) and (\ref{phitilde}), we can write
\begin{eqnarray}
&& \wt\Phi(\CX^{\rm CS}_t) - \Phi(\CX^{\rm CS}_0) \label{Diff} \\
\!\!\!\!\!\!
&=& \!\!\!\!\!\!\!\!\!\!\!\!\!    \sum_{\uk \in \{0,\cdots, m\}^n} \!\!\!\!\!\!\!\!\!\!
 \big(\phi(d_t^{\uk\times\uk}) -\phi(d_0^{\uk\times\uk})\big)\, G(\uy+\uB(t), \uk)
+ \!\!\!\!\!\!\!\!\!   \sum_{\uk \in \{0,\cdots, m\}^n} \!\!\!\!\!\!\!\!\!\!
\phi(d_0^{\uk\times\uk}) \big(G(\uy+\uB(t), \uk) - G(\uy, \uk)\big). \nonumber
\end{eqnarray}
\medskip

\noindent
{{\bf Step~3 (Aging). }} We first identify the aging term $L^{\rm CS}_{\rm a}$ in the generator, defined in (\ref{LCSa}). For each term in the first sum in (\ref{Diff}), by Taylor expanding $\phi(d_t^{\uk\times\uk})$ in $t$, it is easy to see that

\be{aginglimit}
\lim_{t\to 0} t^{-1}
\E[\big(\phi(d_t^{\uk\times\uk}) -\phi(d_0^{\uk\times\uk})\big)\, G(\uy+\uB(t),\uk)]
= 2 \!\!\!\!  \sum_{1\leq i<j\leq n}
\frac{\partial \phi}{\partial r_{ij}}(d_0^{\uk\times\uk})\, G(\uy, \uk),
\ee
where we note that $\frac{\partial \phi}{\partial r_{ij}}(d_0^{\uk\times\uk})=0$ if $k_i=k_j$, since in this case the $ij$-th argument of $\phi$, $d_0^{k_i k_j}$, is less than $\delta$.
Summing the above limit over $\uk \in \{0,\cdots, m\}^n$ and
using the definition of $G$, we find
\be{aginglimit2}
\begin{aligned}
&\lim_{t\to 0}\ \frac{1}{t}   \sum_{\uk \in \{0,\cdots, m\}^n}
\E\left[\big(\phi(d_t^{\uk\times\uk}) -\phi(d_0^{\uk\times\uk})\big)\, G(\uy+\uB(t),\uk)\right]\\[2ex]
& \hspace{6cm}
= 2 \int_{\R^n} g(\ux) \;\; \sum_{1\leq i<j\leq n} \frac{\partial \phi}{\partial r_{ij}}
(\uur)\, \dd\ux,
\end{aligned}
\ee
where $\uur := (r^{\rm CS}_0(x_i, x_j))_{1\leq i,j\leq n}$.  This gives the aging term $L^{\rm CS}_{\rm a}$.
\medskip

\noindent
{{\bf Step~4 (Resampling). }} We next identify the resampling term $L^{\rm CS}_{\rm r}$, defined in (\ref{LCSr}). For the second sum in (\ref{Diff}), we need to compute
\be{difGlim}
\begin{aligned}
   &\quad  \lim_{t\to 0} \ \frac{1}{t}   \sum_{\uk \in \{0,\cdots, m\}^n} \!\!\!
\phi(d_0^{\uk\times\uk}) \E[G(\uy+\uB(t), \uk) - G(\uy, \uk)] \\
= & \sum_{\uk \in \{0,\cdots, m\}^n} \!\!\!\!
\phi(d_0^{\uk\times\uk})\ \lim_{t\to 0} \frac{1}{t}\ \E[G(\uy+\uB(t), \uk) - G(\uy, \uk)].
\end{aligned}
\ee
For each $\uk \in \{0,\cdots, m\}^n$,
\be{rs11}
G(\uy+\uB(t), \uk) - G(\uy, \uk) = \!\!\!  \int_{\R^n} \!\!\!\!
g(\ux) \Big(\prod_{i=1}^n 1_{[y_{k_i}+B_{k_i}(t), y_{k_i+1}+B_{k_i+1}(t)]}(x_i)
- \prod_{i=1}^n 1_{[y_{k_i}, y_{k_i+1}]}(x_i)  \Big) \dd\ux.
\ee
We can rewrite the difference of the product of indicators as
\be{rs12}
\prod_{i=1}^n \big( 1_{[y_{k_i}, y_{k_i+1}]}(x_i) - 1_{[y_{k_i}, y_{k_i}+B_{k_i}(t)]}(x_i)
+ 1_{[y_{k_i+1},y_{k_i+1}+B_{k_i+1}(t)]}(x_i) \big)
- \prod_{i=1}^n 1_{[y_{k_i}, y_{k_i+1}]}(x_i) .
\ee
We can expand the first product above and sort the result into three groups of terms,
(G$_1$), (G$_2$) and (G$_3$), depending on whether each term contains one, two, or more factors of the form
$1_{[y_i, y_i+B_i(t)]}$ for some $0\leq i\leq m+1$. If $h(\ux)$ denotes a term in (G$_3$), then necessarily
$\int_{\R^n} g(\ux) h(\ux) \leq C_g |B_{i_1}(t)B_{i_2}(t) B_{i_3}(t)|$ for some $C_g$ depending only on $g$
and some $i_1, i_2, i_3 \in \{0,\cdots, m+1\}$. Since $\E|B_{i_1}(t)B_{i_2}(t)B_{i_3}(t)| \leq ct^{\frac{3}{2}}$, terms in (G$_3$) do not contribute to the limit in (\ref{difGlim}).

Each term in (G$_2$) is of the following form, where $1\leq i\neq j\leq n$,
\be{G2a}
\begin{aligned}
1_{[y_{k_i}, y_{k_i}+B_{k_i}(t)]}(x_i) 1_{[y_{k_j}, y_{k_j}+B_{k_j}(t)]}(x_j)
& \prod_{1\leq \tau\leq n \atop \tau \neq i,j} 1_{[y_{k_\tau}, y_{k_\tau+1}]}(x_\tau),  \\
1_{[y_{k_i+1}, y_{k_i+1}+B_{k_i+1}(t)]}(x_i) 1_{[y_{k_j+1}, y_{k_j+1}+B_{k_j+1}(t)]}(x_j)
& \prod_{1\leq \tau\leq n \atop \tau \neq i,j} 1_{[y_{k_\tau}, y_{k_\tau+1}]}(x_\tau),  \\
- 1_{[y_{k_i+1}, y_{k_i+1}+B_{k_i+1}(t)]}(x_i) 1_{[y_{k_j}, y_{k_j}+B_{k_j}(t)]}(x_j)
& \prod_{1\leq \tau\leq n \atop \tau \neq i,j} 1_{[y_{k_\tau}, y_{k_\tau+1}]}(x_\tau),  \\
- 1_{[y_{k_i}, y_{k_i}+B_{k_i}(t)]}(x_i) 1_{[y_{k_j+1}, y_{k_j+1}+B_{k_j+1}(t)]}(x_j)
& \prod_{1\leq \tau\leq n \atop \tau \neq i,j} 1_{[y_{k_\tau}, y_{k_\tau+1}]}(x_\tau).
\end{aligned}
\ee
Denote the four terms in (\ref{G2a}) respectively by $h^{(1)}_{ij}(\ux)$, $h^{(2)}_{ij}(\ux)$,
$h^{(3)}_{ij}(\ux)$ and $h^{(4)}_{ij}(\ux)$. Then
\be{G2error}
\int_{\R^n} h^{(1)}_{ij}(\ux) g(\ux) \dd\ux = \idotsint\limits_{y_{k_\tau}<x_\tau<y_{k_\tau+1} \atop \tau\neq i,j}
\!\!\!
\int\limits_{x_i=y_{k_i}}^{y_{k_i}+B_{k_i}(t)}\int\limits_{x_j=y_{k_j}}^{y_{k_j}+B_{k_j}(t)}
\!\!\!\!\!\!\!\!
\big(g(x_1,\cdots, y_{k_i},\cdots, y_{k_j},\cdots, x_n) + o(1)) \dd\ux,
\ee
where we replaced $g(\ux)$ by $g(x_1,\cdots, y_{k_i},\cdots, y_{k_j},\cdots, x_n)$, with an error of $o(1)$ as $t\downarrow 0$. Therefore
\be{rs14}
\lim_{t\to 0} \frac{1}{t} \E\left[\int_{\R^n} h^{(1)}_{ij}(\ux)g(\ux)\dd\ux\right]
= \delta_{k_i, k_j} \idotsint\limits_{y_{k_\tau}<x_\tau<y_{k_\tau+1} \atop \tau\neq i,j} \!\!\!
g(x_1,\cdots, y_{k_i},\cdots, y_{k_j},\cdots x_n) \prod_{1\leq \tau\leq n\atop \tau\neq i,j} \dd x_\tau.
\ee
We obtain similar results for $h^{(2)}_{ij}$, $h^{(3)}_{ij}$ and $h^{(4)}_{ij}$.

For a fixed pair $i\neq j$,
when we sum over $\uk$ and all contributions from $h^{(\cdot)}_{ij}$ in (\ref{difGlim}), we obtain an
integral for $\phi(\uur)g(\ux)$, where $x_\tau$, $\tau\neq i,j$, are still integrated over
$\R^n$, however the integration for $x_i$ and $x_j$ are replaced by summation over
$\{y_{\sigma}\}_{1\leq \sigma\leq m}$. Contributions only come from $x_i=x_j$, and is positive when
$k_i=k_j$, and negative when $k_i=k_j+1$ or $k_j=k_i+1$. Writing everything in terms of $\CX^{\rm CS}_0$,
we easily verify that the contribution of terms in (G$_2$) to the limit in (\ref{difGlim}) is exactly
\be{rs15}
\sum_{1\leq i\neq j\leq n} \int_{\R^{n-2}} \sum_{x_i, x_j \in \{y_\sigma^+, y_\sigma^- : 1\leq\sigma\leq m\}}
\!\!\!\!\!\!
1_{\{x_i=x_j\}}g(\ux) (\theta_{ij}\phi-\phi)(\uur)
\prod_{1\leq \tau\leq n\atop \tau\neq i,j} \dd x_\tau,
\ee
where $\theta_{ij}\phi$ is defined as in (\ref{gs3}), and $\uur := (r^{\rm CS}_0(x_i, x_j))_{1\leq i,j\leq n}$. This gives the resampling term $L^{\rm CS}_{\rm r}$ defined in (\ref{LCSr}).
\medskip

\noindent
{{\bf Step~5 (Migration).  }} Lastly we identify the diffusion (migration) term $L^{\rm CS}_{\rm d}$, defined in (\ref{LCSd}).
We note that each term in group (G$_1$) is of the following form, where $1\leq j \leq n$,
\be{rs16}
\begin{aligned}
-  1_{[y_{k_j}, y_{k_j}+B_{k_j}(t)]}(x_j) & \prod_{1\leq i\leq n \atop i\neq j} 1_{[y_{k_i}, y_{k_i+1}]}(x_i),  \\
1_{[y_{k_j+1}, y_{k_j+1}+B_{k_j+1}(t)]}(x_j) & \prod_{1\leq i\leq n \atop i\neq j} 1_{[y_{k_i}, y_{k_i+1}]}(x_i),
\end{aligned}
\ee
Denote the two terms respectively by $h^{(1)}_j(\ux)$ and $h^{(2)}_j(\ux)$. Then
\begin{eqnarray}
&\quad & \int_{\R^n} \!\!\! h^{(1)}_j(\ux) g(\ux) \dd\ux  \nonumber\\
&\quad =& - \idotsint\limits_{y_{k_i}<x_i<y_{k_i+1} \atop i\neq j} \!\!\!
\int\limits_{x_j=y_{k_j}}^{y_{k_j}+B_{k_j}(t)}  \!\!\!\!\!\!
\left(g(\ux\big|_{x_j=y_{k_j}}) + (x_j-y_{k_j})\frac{\partial g}{\partial x_j}(\ux\big|_{x_j=y_{k_j}}) + o(|x_j-y_{k_j}|)\right)
\dd\ux. \label{G1error}
\end{eqnarray}
Therefore,
\be{rs17}
\lim_{t\to 0} \frac{1}{t} \E\left[\int_{\R^n} \!\!\! h^{(1)}_j(\ux) g(\ux) \dd\ux \right] =
- \frac{1}{2} \idotsint\limits_{y_{k_i}<x_i<y_{k_i+1} \atop i\neq j}
\!\!\! \frac{\partial g}{\partial x_j}(\ux\big|_{x_j=y_{k_j}}) \prod_{1\leq i\leq n \atop i\neq j} \dd x_i.
\ee
A similar result holds for $h^{(2)}_j$. Combining the two, we see that
\be{rs18}
\lim_{t\to 0} \frac{1}{t} \E\left[\int_{\R^n} \!\!\! \big(h^{(1)}_j(\ux)+h^{(2)}_j(\ux)\big)
g(\ux) \dd\ux\right]
= \frac{1}{2} \idotsint\limits_{y_{k_i}<x_i<y_{k_i+1}}
\!\!\! \frac{\partial^2 g}{\partial^2 x_j}(\ux) \prod_{1\leq i\leq n} \dd x_i.
\ee
Therefore, the contribution from terms in (G$_1$) to the limit in (\ref{difGlim}) gives precisely
the migration term $L^{\rm CS}_{\rm d}$ defined in (\ref{LCSd}). This establishes the generator formula (\ref{generator1}).
\qed

\subsubsection{Proof of \eqref{inteqn1} in Proposition \ref{P:genaction}}
\label{sss.proofPgena}

The complications in proving \eqref{inteqn1} arise from trying to prove uniform integrability for
various quantities. We proceed in three steps. First we show that, for each $t>0$,
\be{rderiv}
\lim_{h\downarrow 0} \frac{\E[\Phi(\CX^{\rm CS}_{t+h})] - \E[\Phi(\CX^{\rm CS}_t)]}{h} = \E[{{}L^{\mathrm{CS}}\Phi}(\CX^{\rm CS}_t)].
\ee
By the Markov property of $(\CX^{\rm CS}_t)_{t\geq 0}$,
\be{rderiv2}
\lim_{h\downarrow 0} \frac{\E[\Phi(\CX^{\rm CS}_{t+h})] - \E[\Phi(\CX^{\rm CS}_t)]}{h}
= \lim_{h\downarrow 0} \E\Bigg[\frac{\E[\Phi(\CX^{\rm CS}_{t+h}) | \CX^{\rm CS}_t]-\Phi(\CX^{\rm CS}_t)}{h}\Bigg].
\ee
Since $\CX^{\rm CS}_t \in \U^\R_{\rm r}$ a.s.\ by Proposition \ref{P:reg},
$\lim_{h\downarrow 0} \frac{\E[\Phi(\CX^{\rm CS}_{t+h}) | \CX^{\rm CS}_t]-\Phi(\CX^{\rm CS}_t)}{h} = {{}L^{\mathrm{CS}}\Phi}(\CX^{\rm CS}_t)$ almost surely.
Therefore, the first step is to interchange limit and expectation in (\ref{rderiv2}) and to deduce (\ref{rderiv}). We need to show that
\be{sun25}
\Big(\frac{\E[\Phi(\CX^{\rm CS}_{t+h}) | \CX^{\rm CS}_t]-\Phi(\CX^{\rm CS}_t)}{h}\Big)_{h>0}
\mbox{ is uniformly integrable for } h>0 \mbox{ small, say } 0<h < \frac{\delta\wedge t}{2}.
\ee

Once the uniform integrability has been verified, in Step 2 {{}we prove} that $\E[{{}L^{\mathrm{CS}}\Phi}(\CX^{\rm CS}_t)]$ is continuous in $t$,
and then in Step 3 {we} put things together and  prove (\ref{inteqn1}).
\medskip

\noindent
{\bf Step 1 (Uniform integrability).} This step constitutes the bulk of the proof of (\ref{inteqn1}).
Let us fix $\CX^{\rm CS}_t$ and examine the error terms in our earlier calculation of the generator formula (\ref{generator1}) with $\CX^{\rm CS}_0$ replaced by $\CX^{\rm CS}_t$. Instead of partitioning $X^{\rm CS}_t$ into disjoint open balls of radius $\delta/4$ as done in (\ref{EdL1}), we set
\be{EdL2}
\{y_1< \cdots <y_{m_t}\} := E^{\frac{\delta\wedge t}{4}}_t \cap [-2L,2L],
\ee
which is determined by the dual Brownian web $\widehat \CW$ as seen from the definition in (\ref{Elt}). This choice of partitioning of $X^{\rm CS}_t$ removes the dependence of $\uy=(y_1,\ldots, y_{m_t})$ on the initial condition $\CX^{\rm CS}_0$. Note that $m_t$ depends on $\CX^{\rm CS}_t$.

Let $\E_t[\cdot]$ denote the conditional expectation $\E[\cdot |\CX^{\rm CS}_t]$. Following the arguments leading to (\ref{meeterror}),
\begin{eqnarray}
&& h^{-1}\big|\E_t[\Phi(\CX^{\rm CS}_{t+h})] - \Phi(\CX^{\rm CS}_t)\big| \nonumber \\
&=& h^{-1}\Big|\E_t[\Phi(\CX^{\rm CS}_{t+h}) - \wt \Phi(\CX^{\rm CS}_{t+h})]+
h^{-1}\E_t[\wt\Phi(\CX^{\rm CS}_{t+h}) -\Phi(\CX^{\rm CS}_t)]\Big| \nonumber \\
&\leq& h^{-1}\E_t\big[1_{\{{\rm meet}\}_h} |\Phi(\CX^{\rm CS}_{t+h}) - \wt \Phi(\CX^{\rm CS}_{t+h})|\big]
+ h^{-1}\big|\E_t[\wt\Phi(\CX^{\rm CS}_{t+h}) -\Phi(\CX^{\rm CS}_t)]\big|, \label{errordec}
\end{eqnarray}
where $\{{\rm meet}\}_h$ is the event that either $y_0+B_0(s)$ hits level $-L$ before time $h$, or $y_{m_t+1}+B_{m_t+1}(s)$ hits level $L$ before time $h$, or one of the pair $(y_i+B_i(s), y_{i+1}+B_{i+1}(s))$ coalesces before time $h$. We now estimate the two terms in (\ref{errordec}).
\medskip

{\bf (i)} We start with the second term in (\ref{errordec}).  Based on the decomposition (\ref{Diff}), $\wt\Phi(\CX^{\rm CS}_{t+h}) -\Phi(\CX^{\rm CS}_t)=H(h, \uB(h))$ for some function
$H(h, z_0, \cdots, z_{m_t+1})$ which is continuously differentiable in $h$ and three times continuously differentiable in $z_0,\ldots, z_{m_t+1}$ with uniformly bounded derivatives. Since the generator formula (\ref{generator1}) is derived by Taylor expanding $H(h,\uB(h))$, it is not hard to see that uniformly in $h\in (0, \frac{\delta\wedge t}{4})$,
\be{rs19}
\Big|\E_t[\wt\Phi(\CX^{\rm CS}_{t+h}) -\Phi(\CX^{\rm CS}_t)] - hL^{\rm CS}\Phi(\CX^{\rm CS}_t)\Big|
= \big|\E_t[H(h,\uB(h))]- hL^{\rm CS}\Phi(\CX^{\rm CS}_t)\big| \leq C_{g,\phi} (m_t+3) h
\ee
for some $C_{g,\phi}$ depending on $g$ and $\phi$, and $m_t+3$ is the number of variables in $H(h, z_0,\ldots, z_{m_t+1})$. Therefore
\be{rs20}
h^{-1}|\E_t[\wt\Phi(\CX^{\rm CS}_{t+h}) -\Phi(\CX^{\rm CS}_t)]| \leq C_{g,\phi}(m_t+3) + |L^{\rm CS}\Phi(\CX^{\rm CS}_t)|
\ee
uniformly for $h \in (0, \frac{\delta\wedge t}{4})$. From the definition of $L^{\rm CS}\Phi$, we note that
\be{rs21}
|L^{\rm CS}\Phi(\CX^{\rm CS}_t)| \leq C_{g,\phi}(1 + m_t) <\infty.
\ee
By the definition of $m_t$ in (\ref{EdL2}) and by Lemma~\ref{L:density1}, we have $\E[m_t]\leq \frac{4L}{\sqrt{\pi \frac{\delta \wedge t}{4}}} <\infty$. This implies the uniform integrability of the second term in (\ref{errordec}) for $h\in (0, \frac{\delta\wedge t}{4})$.
\medskip

{\bf (ii)}
We now consider the first term in (\ref{errordec}). For $0\leq i \leq m_t$, let $\tau_{i,i+1}$ be the first hitting time
between $y_i+B_i(s)$ and $y_{i+1}+B_{i+1}(s)$. Let $\tau_0$ be the first hitting time of level $-L$ by $y_0+B_0(s)$, and $\tau_{m_t+1}$ the first hitting time of level $L$ by $y_{m_t+1}+B_{m_t+1}(s)$. Then
\be{hitsplit}
\begin{aligned}
& h^{-1}\E_t\big[1_{\{{\rm meet}\}_h} |\Phi(\CX^{\rm CS}_{t+h})\! -\! \wt \Phi(\CX^{\rm CS}_{t+h})|\big]
\\
\leq \ \ & 2|\Phi|_\infty h^{-1} \E_t\big[1_{\{\tau_0\leq h\}}\! +\! {{}1_{\{\tau_{m_t+1}\leq h\}}}\big] +
h^{-1}\E_t\Big[\!\sum_{i=0}^{m_t} 1_{\{\tau_{i,i+1}\leq h\}} \big|\Phi(\CX^{\rm CS}_{t+h})\! -\! \wt \Phi(\CX^{\rm CS}_{t+h})\big|\Big].
\end{aligned}
\ee
Since $y_0<-2L$ and $y_{m_t+1}>2L$, the probability of the events $\{\tau_0\leq h\}$ and $\{\tau_{m_t+1}\leq h\}$ decay exponentially fast
in $h^{-1}$ and the events are independent of $\CX^{\rm CS}_t$. Therefore the first term in
(\ref{hitsplit}) is uniformly bounded in $h>0$.

Bounding the second term in (\ref{hitsplit}) is more delicate, because it remains of order $1$ as $h\downarrow 0$. We will need to use negative correlation inequalities for coalescing Brownian motions established in Appendix~\ref{A:Corr}.

Let us recall the definition of $\wt \Phi(\CX^{\rm CS}_{t+h})$ from (\ref{phitilde}), where we replace the pair $\{0,t\}$ by $\{t, t+h\}$. By integrating over the population at time $t+h$, we note that both $\Phi(\CX^{\rm CS}_{t+h})$ and $\wt \Phi(\CX^{\rm CS}_{t+h})$ can be
written as integrals of $g(\ux)\phi(\uur(\ux))$ integrated over $\ux=(x_1,\cdots, x_n)\in (-L,L)^n$, except that: for a given $\ux$, the distance matrix $\uur(\ux)$ may be different for $\Phi$ and $\wt\Phi$, and for $\wt\Phi$, the same point $\ux$ may be integrated over several times with different distance matrix $\uur(\ux)$ due to the fact that $\{y_i+B_i(h)\}_{0\leq i\leq m_t+1}$ may not have the same order as $\{y_i\}_{0\leq i\leq m+1}$. However, for $\ux$ not in
\be{rs23}
D:=\Big\{\ux\in\R^n : x_i  \in \bigcup_{s\in [0,h]}\{y_j+B_j(s)\} \mbox{ for some } 1\leq i\leq n\mbox{ and } 0\leq j\leq m_t+1\Big\},
\ee
$\ux$ is integrated over exactly once in $\wt\Phi$, and the associated distance matrix $\uur (\ux)$ is the same for both $\Phi$ and $\wt\Phi$. Therefore, contributions from $\ux\notin D$ cancel out in
$|\Phi(\CX^{\rm CS}_{t+h})-\wt\Phi(\CX^{\rm CS}_{t+h})|$. Since $g$ has compact support, the contribution
from $\ux\in D$ to $|\Phi(\CX^{\rm CS}_{t+h})-\wt\Phi(\CX^{\rm CS}_{t+h})|$, including multiple integrations
over the same $\ux$ by $\wt\Phi$, is at most
$C_{\phi, g}n |g|_\infty|\phi|_\infty \sum_{i=0}^{m_t+1}R_i(h)$, where
$R_i(h):=(\sup_{0\leq s\leq h}B_i(s)-\inf_{0\leq s\leq h}B_i(s))$. Therefore,
\begin{eqnarray}
h^{-1} \E_t\Big[\!\sum_{i=0}^{m_t} 1_{\{\tau_{i,i+1}\leq h\}} \big|\Phi(\CX^{\rm CS}_{t+h})\! -\! \wt
\Phi(\CX^{\rm CS}_{t+h})\big|\Big]
&\leq&
C_{n,\phi, g} h^{-1}\E_{\uB}\Big[\sum_{i=0}^{m_t} 1_{\{\tau_{i,i+1}\leq h\}} \sum_{j=0}^{m_t+1} R_j(h)\Big] \nonumber \\
&\leq& C_{n,\phi, g} h^{-1} \sum_{i=0}^{m_t}\sum_{j=0}^{m_t+1} \P_{\uB}(\tau_{i,i+1}\leq h)^{\frac{1}{2}}\E_{\uB}[R_j(h)^2]^{\frac{1}{2}} \nonumber\\
&\leq& C'_{n,\phi, g} h^{-\frac{1}{2}} (m_t+2) \sum_{i=0}^{m_t} \P_{\uB}(\tau_{i,i+1}\leq h)^{\frac{1}{2}}, \label{yptbd}
\end{eqnarray}
where $\E_{\uB}$ denotes expectation w.r.t.\ the Brownian motions $\uB=(B_0,\ldots, B_{m_t+1})$, and we applied H\"older inequality
and the fact that $\E[R_j(h)^2]^{\frac{1}{2}} = c\sqrt{h}$ for some $c>0$. It only remains to prove the uniform integrability of the r.h.s.\ of (\ref{yptbd}) w.r.t.\ the law of $\CX^{\rm CS}_t$ for $0<h<\frac{\delta\wedge t}{4}$.

Note that the r.h.s.\ of (\ref{yptbd}) depends on $y_0$ and $y_{m_t+1}$, which lie outside $[-2L, 2L]$. To control the dependence on $y_0$ and $y_{m_t+1}$, we enlarge the interval and let $\{z_1,\cdots,z_{M+1}\}:=E_t^{\frac{\delta \wedge t}{4}} \cap (-2L-1, 2L+1)$ as in (\ref{EdL2}), which contains $\{y_1, \cdots, y_{m_t}\}$ as a subset. Denote
$$
\psi\big(\frac{y_{i+1}-y_i}{\sqrt h}\big) :=\P(\tau_{i,i+1}\leq h)^{\frac{1}{2}} = \Big(2\int^\infty_{\frac{y_{i+1}-y_i}{\sqrt{2h}}} e^{-\frac{x^2}{2}}dx\Big)^{\frac{1}{2}},
$$
we can then replace the r.h.s.\ of (\ref{yptbd}) by
\be{rs24}
F_h(\CX^{\rm CS}_t) := \frac{M}{\sqrt h}  \sum_{i=1}^{M}  \psi\big(\frac{z_{i+1}-z_i}{\sqrt h}\big),
\ee
because $C'_{n,\phi, g} F_h(\CX^{\rm CS}_t)$ dominates the r.h.s.\ of (\ref{yptbd}) except for possible missing terms $\frac{1}{\sqrt h} \psi\big(\frac{y_{1}-y_0}{\sqrt h}\big)(m_t+2)$, resp.\ $\frac{1}{\sqrt h} \psi\big(\frac{y_{m_t+1}-y_{m_t}}{\sqrt h}\big)(m_t+2)$, on the event $y_0\leq -2L-1$, resp.\ $y_{m_t+1}\geq 2L+1$. Since $y_1\geq -2L$ and $y_{m_t}\leq 2L$ by definition, both missing terms are bounded by $C_1m_t+C_2$ uniformly for $h>0$, and is thus uniformly integrable.

It only remains to prove the uniform integrability of $F_{h}(\CX^{\rm CS}_t)$ uniformly in $0<h<\frac{\delta\wedge t}{4}$. We achieve this by bounding its second moment. Note that by Cauchy-Schwarz,
\begin{eqnarray}
\E[F_{h}(\CX^{\rm CS}_t)^2] &=& \frac{1}{h}\, \E\left[M^2 \Big(\sum_{i=1}^{M}
\psi\big(\frac{{{}z_{i+1}-z_{i}}}{\sqrt h}\big)\Big)^2 \right] \nonumber \\
&\leq& \frac{1}{h}\, \E\left[M^3 \sum_{i=1}^{M}
\psi\big(\frac{{{}z_{i+1}-z_{i}}}{\sqrt h}\big)^2 \right]
\leq \frac{1}{h}\, \E\left[\sum_{i_1, i_2, i_3=1}^{M+1} \sum_{i_4=1}^{M}
\psi\big(\frac{z_{i_4+1}-z_{i_4}}{\sqrt h}\big)^2 \right]. \label{rs25}
\end{eqnarray}
{{}Note that the above summation can be seen as a summation of $x_1:=z_{i_1}, x_2:=z_{i_2}, x_3:=z_{i_3}$ in $\{z_1, \ldots, z_{M+1}\}$ and $x_4:=z_{i_4}$ in $\{z_1, \ldots, z_M\}$. By integrating over $x_1, x_2, x_3, x_4$ over $[-2L-1, 2L+1]$ and the position $x_5:=z_{i_4+1} \in (x_4, 2L+1)$}, we can rewrite the r.h.s.\ of \eqref{rs25}
 in terms of the correlation functions of the translation invariant simple point process $\xi:=E^{\frac{\delta \wedge t}{4}}_t\subset \R$. More precisely, it can be written as
\begin{equation}\label{rs26}
\idotsint\limits_{x_4<x_5\atop |x_1|, |x_2|, |x_3|, |x_4|, |x_5|< 2L+1} \frac{1}{h}\psi\big(\frac{x_5-x_4}{\sqrt h}\big)^2 K^{\rm c}_{4,5}(x_1, x_2, x_3, x_4, x_5) \dd x_1 \dd x_2 \dd x_3 \dd x_4 \dd x_5,
\end{equation}
where
\be{sun25-2}
K^{\rm c}_{4,5} (x_1, \cdots, x_5) := \lim_{\epsilon\downarrow 0} \epsilon^{-5}\P\big([x_j, x_j+\epsilon]\cap \xi \neq\emptyset
\ \forall\ 1\leq j\leq 5, (x_4, x_5)\cap \xi=\emptyset \big)
\ee
is the density of finding a point at $x_i$ for each $1\leq i\leq 5$, with no point in $(x_4, x_5)$. By the definition of $E^l_t$ in (\ref{Elt}) and the duality between the forward and dual Brownian web $(\CW, \widehat \CW)$, we see that $\xi$ is the point process generated on $\R$ at time $t$ by coalescing Brownian motions in the Brownian web $\CW$ starting from every point in $\R$ at time $t-\frac{\delta \wedge t}{4}$. By the negative correlation inequality in Lemma~\ref{L:negcor2}, we can bound
$$
\begin{aligned}
& K^{\rm c}_{4,5} (x_1, \cdots, x_5) \\
\leq\ & \lim_{\epsilon\downarrow 0} \epsilon^{-5}\P\big([x_j, x_j+\epsilon]\cap \xi \neq\emptyset,
\ j=1,2,3\big) \P\big([x_j, x_j+\epsilon]\cap \xi \neq\emptyset, j=4,5; \ (x_4, x_5)\cap \xi=\emptyset \big) \\
=\ & K(0)^3 K^{\rm c}(x_4, x_5),
\end{aligned}
$$
where by Lemma~\ref{L:density1},
$$
K(0) = \lim_{\epsilon\downarrow 0} \epsilon^{-1} \P\big([x,x+\epsilon] \cap \xi\neq\emptyset) = \frac{2}{\sqrt{\pi \delta\wedge t}} \qquad \mbox{ for all } x\in\R,
$$
and
$$
K^{\rm c}(x, y) = \lim_{\epsilon\downarrow 0} \epsilon^{-2} \P\big([x,x+\epsilon] \cap \xi\neq\emptyset, (x+\epsilon,y)\cap \xi=\emptyset, [y, y+\epsilon]\cap \xi\neq\emptyset\big) \quad \mbox{ for } x<y.
$$
By Lemma \ref{L:c2point}, for $x<y$,
$$
K^{\rm c}(x, y) \leq C_{\delta,t}(y-x)\wedge 1
$$
for some $C_{\delta,t}>0$ depending only on $\delta$ and $t$. Substituting the above bounds into (\ref{rs26}), using the definition of $\psi$, and separating the integration into two regions depending on whether $0<x_5-x_4<\sqrt{h}$ or $x_5-x_4\geq \sqrt{h}$, it is easily seen that the integral in (\ref{rs26}) is uniformly bounded for $0<h<\frac{\delta\wedge t}{4}$. This implies the uniform integrability of $F_h(\CX^{\rm CS}_t)$ for $0<h<\frac{\delta\wedge t}{4}$, and hence that of $\frac{\E_t[\Phi(\CX^{\rm CS}_{t+h})]-\Phi(\CX^{\rm CS}_t)}{h}$.

\bigskip

\noindent
{\bf Step 2 (Continuity of $\E[L^{\rm CS}\Phi(\CX^{\rm CS}_t)]$).}
Recall from (\ref{LCS})--(\ref{LCSr}) that $L^{\rm CS}=L^{\rm CS}_{\rm d}+L^{\rm CS}_{\rm a}+L^{\rm CS}_{\rm r}$.
We will prove the continuity for each component of the generator.

By our assumptions on $g$ and $\phi$, $|L^{\rm CS}_{\rm d}\Phi|_\infty, |L^{\rm CS}_{\rm a}\Phi|_\infty \leq C_{\Phi}<\infty$. It was shown in (\ref{rcont}) that for any $t> 0$, for Lebesgue a.e.\ $\ux=(x_1,\cdots, x_n)\in\R^n$, the distance matrix $r^{\rm CS}_s(x_i, x_j)$ converges a.s.\ to $r^{\rm CS}_t(x_i, x_j)$ as $s\to t$, and this conclusion is also easily seen to hold for $t=0$ by the assumption $\CX^{\rm CS}_0\in \U^\R_{\rm r}$. Therefore, almost surely w.r.t.\ $(\CX^{\rm CS}_s)_{s\geq 0}$,
\be{sun26}
\lim_{s\to t}L^{\rm CS}_{\rm d}\Phi(\CX^{\rm CS}_{s})=L^{\rm CS}_{\rm d}\Phi(\CX^{\rm CS}_t)\quad  \mbox{ and } \quad
\lim_{s\to t}L^{\rm CS}_{\rm a}\Phi(\CX^{\rm CS}_s) =L^{\rm CS}_{\rm a}\Phi(\CX^{\rm CS}_t).
\ee
Therefore $\E[L^{\rm CS}_{\rm d}\Phi(\CX^{\rm CS}_t)]$ and $\E[L^{\rm CS}_{\rm a}\Phi(\CX^{\rm CS}_t)]$ are continuous in $t\geq 0$ by the bounded convergence theorem.

We now turn to the continuity of $\E[L^{\rm CS}_{\rm r}\Phi(\CX^{\rm CS}_t)]$. We first prove the continuity  at $t=0$ and later point out how it extends to $t >0$. Let $E_0^{\frac{\delta}{4}}\subset \R$ be defined from $\CX^{\rm CS}_0$ as in (\ref{Elt}), which are the boundaries of balls of radius $\frac{\delta}{4}$ in $\CX^{\rm CS}_0$. Let us follow paths in the Brownian web $\CW$ starting from points in $E_0^{\frac{\delta}{4}}\subset \R$ at time $0$, which determine the evolution of these boundaries, and denote the point set on $\R$ generated at time $s>0$ by $E_0^{\frac{\delta}{4}}(s)$. Then $E^{\frac{\delta}{4}+s}_s \subset E_0^{\frac{\delta}{4}}(s)$. By our regularity assumption on $\phi$, only resampling at the boundaries of balls of radius $\delta$ or more has an effect on $\Phi(\CX^{\rm CS}_s)$. Therefore for $s\in (0, \delta/4)$, we have
\be{L3Phil0}
L^{\rm CS}_{\rm r}\Phi(\CX^{\rm CS}_s) = \sum_{1\leq i\neq j\leq n} \int_{\R^{n-2}} \sum_{x_i, x_j \in \{y^+, y^- : y \in E_0^{\delta/4}(s)\}}
\!\!\!\!\!\!
1_{\{x_i=x_j\}}g(\ux) (\theta_{ij}\phi-\phi)(\uur)
\prod_{1\leq \tau\leq n\atop \tau\neq i,j} \dd x_\tau,
\ee
where $\theta_{ij}\phi$ is defined as in (\ref{gs3}), and $\uur := (r^{\rm CS}_s(x_i, x_j))_{1\leq i,j\leq n}$. By our assumptions on $g$ and $\phi$, the fact $\CX^{\rm CS}_0\in \U^\R_{\rm r}$, and our construction of $\CX^{\rm CS}_s$ in terms of the (dual) Brownian web, it is then easily seen by dominated convergence that
\be{sun27}
\lim_{s\downarrow 0}L_{\rm r}\Phi(\CX^{\rm CS}_s) = L_{\rm r}\Phi(\CX^{\rm CS}_0)
\quad \mbox{ almost surely}.
\ee

To prove the continuity of $\E[L^{\rm CS}\Phi(\CX^{\rm CS}_t)]$ at $t=0$, it only remains to verify the uniform integrability of $L^{\rm CS}_{\rm r}\Phi(\CX^{\rm CS}_s)$ for $s$ close to $0$, say $s\in [0,\delta/4]$. We will achieve this by showing that $L^{\rm CS}_{\rm r}\Phi(\CX^{\rm CS}_s)$ has uniformly bounded second moments.

Note that because $g$ is assumed to be supported on $[-L, L]^n$, we have
\be{rs27}
|L^{\rm CS}_{\rm r}\Phi(\CX^{\rm CS}_s)| \leq C_{n, \phi, g} \big|E^{\delta/4}_0(s) \cap (-L,L)\big|.
\ee
By Lemma \ref{L:negcor}, $E^{\delta/4}_0(s)$ is negatively correlated, and hence by Lemma \ref{L:negmom}, we have
\be{rs28}
\E\Big[\big|E^{\delta/4}_0(s) \cap (-L,L)\big|^2\Big] \leq
2\E\Big[\big|E^{\delta/4}_0(s) \cap (-L,L)\big|\Big]+
\Big[\E\big|E^{\delta/4}_0(s) \cap (-L,L)\big|\Big]^2.
\ee
Thus it suffices to bound $\E\big[\big|E^{\delta/4}_0(s)\cap(-L,L)\big|\big]$ uniformly in $s\in [0, \delta/4]$.

Since $E^{\delta/4}_0(s)$ is obtained by evolving coalescing Brownian motions starting from $E^{\delta/4}_0$, we can bound
\be{Zsplit}
\E\big[\big|E^{\delta/4}_0(s) \cap (-L,L)\big|\big]
\leq \big|E^{\delta/4}_0\cap[-2L,2L]\big| + 2
\sum_{i\geq 2L} \E\big[\big|\xi_s^{[i,i+1]\times\{0\}} \cap (-L,L)\big|\big],
\ee
where $\xi^{[i,i+1]\times\{0\}}_s\subset \R$ is the point set generated at time $s$ by coalescing Brownian motions in the Brownian web $\CW$ starting from everywhere in $[i, i+1]$ at time $0$.

Note that the first term in (\ref{Zsplit}) is finite and independent of $s\geq 0$. We now treat the second term. For each $i\geq 2L$, let us denote $\xi_s^{[i,i+1]\times\{0\}} \cap (-L,L)$ by $\xi^i_{s,L}$, which is also a point process on $(-L,L)$ satisfying negative correlation. In particular,
\be{Zsplit2}
\E[|\xi^i_{s,L}|] = \sum_{n=1}^\infty \P(|\xi^i_{s,L}|\geq n) \leq \sum_{n=1}^\infty \P(|\xi^i_{s,L}\geq 1)^n
\leq \sum_{n=1}^\infty \Big(\frac{1}{\sqrt{2\pi s}}\int_{i-L}^\infty e^{-\frac{x^2}{2s}}dx\Big)^n,
\ee
where we used Lemma \ref{L:ocucor} and the observation that,
$\xi^i_{s,L}\neq \emptyset$ implies that the Brownian motion
starting at $(i,0)$ in the Brownian web must be to the left of $L$ at time $s$. Since $L$ is fixed and large,
$\frac{1}{\sqrt{2\pi s}}\int_{i-L}^\infty e^{-\frac{x^2}{2s}}dx \leq \alpha$ for some $\alpha\in (0,1)$ uniformly
in $i\geq 2L$ and $s\in [0, \delta/4]$, and hence
\be{rs29}
\E[|\xi^i_{s,L}|]\leq \frac{\int_{i-L}^\infty  e^{-\frac{x^2}{2s}}dx}{(1-\alpha)\sqrt{2\pi s}}
= \frac{\int_{\frac{i-L}{\sqrt s}}^\infty  e^{-\frac{x^2}{2}}dx}{(1-\alpha)\sqrt{2\pi }}.
\ee
It is then clear that $\sum_{i\geq 2L} \E[|\xi^i_{s,L}|]$ tends to $0$ as $s\downarrow 0$ and is uniformly
bounded for $s\in [0,\delta/4]$, which concludes the proof of the uniform integrability of $L_{\rm r}\Phi(\CX^{\rm CS}_s)$ for $s\in [0,\delta/4]$. Therefore $\E[L^{\rm CS}_{\rm r}\Phi(\CX^{\rm CS}_t)]$ is continuous at $t=0$.

The continuity of $\E[L^{\rm CS}_{\rm r}\Phi(\CX^{\rm CS}_t)]$ at $t>0$ can be proved similarly. For
$s\in [t-\frac{\delta\wedge t}{4}, t+\frac{\delta\wedge t}{4}]$, we can use the same representation as in
(\ref{L3Phil0}) by treating $t-\frac{\delta\wedge t}{4}$ as the starting time $0$ in (\ref{L3Phil0}). The a.s.\ convergence $\lim_{s\to t}L_{\rm r}\Phi(\CX^{\rm CS}_s)=L_{\rm r}\Phi(\CX^{\rm CS}_t)$ then follows as before. The uniform integrability of $L_{\rm r}\Phi(\CX^{\rm CS}_s)$ follows easily from the density estimate, Lemma \ref{L:density1}, if we restrict to $s\in [t-\frac{\delta\wedge t}{8}, t+\frac{\delta\wedge t}{8}]$, since
only resampling at locations occupied at time $s$ by Brownian web paths starting from $\R$ at time $t-\frac{\delta\wedge t}{4}$ has an effect on $\Phi(\CX^{\rm CS}_s)$.
\bigskip

\noindent
{\bf Step 3 (Proof of (\ref{inteqn1})).} We have thus far shown that $\E[\Phi(\CX^{\rm CS}_t)]$ has a continuous right derivative
$\E[L^{\rm CS}\Phi(\CX^{\rm CS}_t)]$ on $[0,\infty)$. Note that $\E[\Phi(\CX^{\rm CS}_t)]$ is also continuous on $[0,\infty)$ since $(\CX^{\rm CS}_s)_{s\geq 0}$ has a.s.\ continuous sample paths by Theorem~\ref{T:MarkovCont}. Equation (\ref{inteqn1}) then follows from a variant of the fundamental theorem of calculus, formulated as Lemma~\ref{L:FTC} below. This finally completes the proof of (\ref{inteqn1}) and that of Prop.~\ref{P:genaction}.
\qed

\begin{lemma}\label{L:FTC}{\bf [Fundamental Theorem of Calculus]} \\
Let $f:\R\to\R$ be continuous, with a continuous right derivative $D^+f(x):=\lim_{h\downarrow 0} \frac{1}{h}(f(x+h)-f(x))$. Then
$$
f(x) = f(0)+\int_0^x D^+f(y) {\rm d}y \qquad \mbox{for all } x\in \R.
$$
\end{lemma}
\noindent
{\bf Proof.} Let $g(x):= f(0)+\int_0^x D^+f(y) {\rm d}y$, which is clearly differentiable with continuous {{}right} derivative $D^+f$. Then $h:=f-g$ is continuous with $h(0)=0$ and right derivative $D^+h\equiv 0$. It suffices to show that $h\equiv 0$.

For $\epsilon>0$, let $h_\eps(x)=h(x)+\eps x$. Then $D^+h_\eps\equiv \eps>0$. This implies that $h_\eps$ is non-decreasing on $\R$, since otherwise if $h_\eps(x) > h_\eps(y)$ for some $x<y$, then at the point $x_0\in [x,y]$ where $h_\eps$ achieves its maximum on $[x,y]$, we have $D^+h_\eps(x_0)\leq 0$, contradicting $D^+h_\eps\equiv\eps$. Letting $\eps\downarrow 0$, the monotonicity of $h_\eps$ then implies that $h$ is also non-decreasing. Applying the same argument to $-h$ shows that $h$ is also non-increasing. Thus $h\equiv 0$.
\qed

\subsection{Generator Action on General Test Functions}\label{ss:generator2}

In this section, we extend the validity of the integral equation (\ref{inteqn1}) in Prop.~\ref{P:genaction} to general $\Phi \in \Pi^{1,2}$, assuming $\CX^{\rm CS}_0\in\U^\R_{\rm r}$. The complication is that the generator identified in (\ref{generator1}), which acts on regular test functions $\Phi \in \Pi^{1,2}_{\rm r}$ and evaluated at regular states $\CX^{\rm CS}_t\in \U^\R_{\rm r}$, can only be extended to general test functions $\Phi\in \Pi^{1,2}$ provided that we restrict to the regular subclass of states $\CX^{\rm CS}_t\in \U^\R_{\rm rr}$ introduced in Definition~\ref{D:Urr}. Fortunately for each $t>0$, almost surely $\CX^{\rm CS}_t\in \U^\R_{\rm rr}$, which makes it possible to extend from $\Phi \in \Pi^{1, 2}_{\rm r}$ to $\Phi \in \Pi^{1,2}$.

\bp\label{P:genaction2}{\bf [Generator action on general test functions]}\\
Let $\CX^{\rm CS}$ be the CSSM genealogy process with $\CX^{\rm CS}_0\in \U^\R_{\rm r}$, and let $\Phi=\Phi^{n,\phi, g}\in \Pi^{1,2}$. Then
\begin{itemize}
\item[\rm (a)] The integral equation \eqref{inteqn1} still holds, and $\E[L^{\rm CS}\Phi(\CX^{\rm CS}_t)]$ is continuous in $t>0$.

\item[\rm (b)] If $\CX^{\rm CS}_0\in \U^\R_{\rm rr}$, defined as in Definition~\ref{D:Urr}, then the generator equation \eqref{generator1} still holds, and $\E[L^{\rm CS}\Phi(\CX^{\rm CS}_t)]$ is continuous in $t\geq 0$.
\end{itemize}
\ep
\noindent
{\bf Proof.} We will approximate $\Phi=\Phi^{n,\phi, g}\in \Pi^{1,2}$ by $\Phi_l\in \Pi^{1,2}_{\rm r}$ as follows. Let $\rho: [0,\infty)\to\R$ be continuously differentiable with $\rho(x)=0$ for $x\in [0,2]$, $\rho'(x)>0$ for $x\in (2,3]$, and $\rho(x)=x$ for $x\in [3,\infty)$. For each $l>0$, denote $\rho_l(x) = l\rho(x/l)$. Then $\rho_l$ is a smooth truncation with $\sup_{x\geq 0}|\rho_l(x)-x| \to 0$
and $\rho_l'(x)\to 1$ for each $x> 0$ as $l\downarrow 0$. Given $\phi=\phi((r_{ij})_{1\leq i<j\leq n})$, let $\phi_l:=\phi\circ \rho_l$. We then define the truncated version of $\Phi$ by $\Phi_l:=\Phi^{n, \phi_l, g}$.

It is easily seen that $\Phi_l\in \Pi^{1,2}_{\rm r}$, and hence Proposition~\ref{P:genaction} can be applied. We will then deduce Proposition~\ref{P:genaction2} by taking the limit $l\downarrow 0$.
\medskip

{\bf (a)} We will take the limit $l\downarrow 0$ in the integral equation (\ref{inteqn1}). In {\bf Step 1}, we will show that (\ref{inteqn1}) also holds for $\Phi$. In {\bf Step 2}, we will prove the continuity of $\E[L^{\rm CS}\Phi(\CX^{\rm CS}_t)]$ in $t>0$. Note that without assuming
$\Phi\in \Pi^{1,2}_{\rm r}$ or $\CX^{\rm CS}_0 \in \U^\R_{\rm rr}$, $L^{\rm CS}\Phi(\CX^{\rm CS}_0)$ may not be well-defined.
\medskip

\noindent
{\bf Step 1}  First note that $|\Phi_l-\Phi|_\infty \to 0$ as $l\downarrow 0$ by our assumption on $g$ and $\phi$ and the fact that $\sup_{x\geq 0}|\rho_l(x)-x|\to 0$ as $l\downarrow 0$. Therefore,
\be{LPhi1}
\lim_{l\downarrow 0}\E[\Phi_l(\CX^{\rm CS}_t)]  = \E[\Phi(\CX^{\rm CS}_t)]\qquad \mbox{and}\qquad
\lim_{l\downarrow 0} \Phi_l(\CX^{\rm CS}_0)=\Phi(\CX^{\rm CS}_0).
\ee

Recall the decomposition of $L^{\rm CS}=L^{\rm CS}_{\rm d}+L^{\rm CS}_{\rm a} + L^{\rm CS}_{\rm r}$ in Prop.~\ref{P:genaction}.
Note that $|L^{\rm CS}_{\rm d}\Phi|_\infty$, $|L^{\rm CS}_{\rm d}\Phi_l|_\infty$, $|L^{\rm CS}_{\rm a}\Phi|_\infty$, $|L^{\rm CS}_{\rm a}\Phi_l|_\infty$ are all uniformly bounded by some $C_{\Phi}<\infty$ independent of $l$. Also, by our assumptions on $g$, $\phi$ and $\rho_l$,
for each $s>0$ and a.s.\ every realization of $\CX^{\rm CS}_s$, we have
\be{sun29}
L^{\rm CS}_{\rm d}\Phi_l(\CX^{\rm CS}_s)\to L^{\rm CS}_{\rm d}\Phi(\CX^{\rm CS}_s)
\mbox{ and } L^{\rm CS}_{\rm a}\Phi_l(\CX^{\rm CS}_s) \to L^{\rm CS}_{\rm a}\Phi(\CX^{\rm CS}_s) \mbox{ as } l\downarrow 0.
\ee
Therefore, by the bounded convergence theorem,
\be{LPhi2}
\lim_{l\downarrow 0} \int_0^t \E[(L^{\rm CS}_{\rm d}\Phi_l+L^{\rm CS}_{\rm a}\Phi_l)(\CX^{\rm CS}_s)]{\rm d}s
= \int_0^t \E[(L^{\rm CS}_{\rm d}\Phi+L^{\rm CS}_{\rm a}\Phi)(\CX^{\rm CS}_s)] {\rm d}s.
\ee

For the resampling generator $L^{\rm CS}_{\rm r}$, note that for any $s>0$ and a.s.\ every realization of $\CX^{\rm CS}_s$,
\begin{equation}\label{LRXsbd}
\begin{aligned}
|L^{\rm CS}_{\rm r}\Phi(\CX^{\rm CS}_s)| &\leq
\sum_{1\leq i\neq j\leq n} \int_{\R^{n-2}} \sum_{x_i, x_j \in \{x^+, x^-: x\in E_s\}}
1_{\{x_i=x_j\}} |g(\ux)|\,
\big|(\theta_{ij}\phi - \phi)(\uur^{\rm CS}_s)\big| \prod_{1\leq \tau\leq n\atop \tau\neq i,j} {\rm d}x_\tau \\
&\leq \ C_{n,g,\phi} \sum_{x\in (-M,M)\cap E_s} d_x \wedge 1,
\end{aligned}
\end{equation}
where $E_s$ is defined as in (\ref{Elt22}), which is the subset of $\R$ that are points of multiplicity in the dual Brwonian web $\widehat \CW$, $d_x=r^{\rm CS}_s(x^+, x^-)$ is the distance between the two points in $X^{\rm CS}_s$ with the same spatial location as $x\in E_s$, $M>0$ is chosen such that ${\rm supp}(g) \subset (-M, M)^n$, and we used the assumption that $\phi$ has bounded derivative. Such a bound holds for $\phi_l$ uniformly in $l>0$. If $\sum_{x\in (-M,M)\cap E_s} d_x\wedge 1 <\infty$, then by dominated convergence, we have
\be{sun30}
\lim_{l\downarrow 0} L^{\rm CS}_{\rm r}\Phi_l(\CX^{\rm CS}_s) = L^{\rm CS}_{\rm r}\Phi(\CX^{\rm CS}_s).
\ee
To prove the analogue of (\ref{LPhi2}) for $L^{\rm CS}_{\rm r}\Phi$, by dominated convergence, it only remains to show
\be{dxfinite}
\int_0^t \E\Big[\sum_{x\in (-M,M)\cap E_s} d_x\wedge 1\Big] {\rm d}s  <\infty.
\ee
Note that for each $s\in (0,t)$,
{
\begin{equation}
\begin{aligned}
& \E\Big[\sum_{x\in (-M,M)\cap E_s} d_x\wedge 1\Big]  \\
&= \E\Big[\sum_{x\in (-M,M)\cap E_s} \int_0^1 1_{\{d_x\wedge 1>u\}}{\rm d}u\Big]
= \int_0^1 \E\Big[\Big|\{x \in (-M,M)\cap E_s : d_x >u\} \Big|\Big] {\rm d}u  \\
& \leq      \int\limits_0^{2s\wedge 1}
\E\Big[\Big|\{x \in (-M,M)\cap E_s : d_x >u\} \Big|\Big] {\rm d}u
+   \int\limits_{2s\wedge 1}^1
\E\Big[\Big|\{x \in (-M,M)\cap E_s : d_x >2s\} \Big|\Big] {\rm d}u  \\
&= \int_0^{2s\wedge 1} \E\Big[\,\Big|\xi^{\R\times\{s-u/2\}}_s \cap (-M,M)\Big|\,\Big] {\rm d}u
+ \int_{2s \wedge 1}^1 \E\Big[\,\Big|\xi^{\R\times\{0\}}_s \cap (-M,M)\Big|\,\Big] {\rm d}u  \\
& = CM\sqrt{2s\wedge 1} + CM\frac{1-2s\wedge 1}{\sqrt s} \leq CM\big(1+\tfrac{1}{\sqrt s}\big), \label{dxint}
\end{aligned}
\end{equation}
}
where $\xi^A_s$ denotes the point set on $\R$ generated at time $s$ by the collection of paths in the
Brownian web $\CW$ starting from the space time set $A\subset \R^2$, and we used $\E[\,|\xi^{\R\times\{s-u\}}_s \cap [a,b]|\,] = \frac{b-a}{\sqrt{\pi u}}$ by Lemma \ref{L:density1}. Inequality (\ref{dxfinite}) then follows. To summarize, we have thus shown that (\ref{inteqn1}) is also valid for
a general polynomial $\Phi\in \Pi^{1,2}$.
\medskip

\noindent
{\bf Step 2}  To prove the continuity of $\E[L^{\rm CS}\Phi(\CX^{\rm CS}_t)]$ in $t> 0$, we again decompose $L^{\rm CS}$ into its
three summands. First note that the continuity of $\E[L^{\rm CS}_{\rm d}\Phi(\CX^{\rm CS}_t)]$ and $\E[L^{\rm CS}_{\rm a}\Phi(\CX^{\rm CS}_t)]$ on $[0,\infty)$ follow by the same arguments as that for $\E[L^{\rm CS}_{\rm d}\Phi(\CX^{\rm CS}_t)]$ and
$\E[L^{\rm CS}_{\rm a}\Phi(\CX^{\rm CS}_t)]$ in the proof of Prop.~\ref{P:genaction}.

To prove the continuity of $\E[L^{\rm CS}_{\rm r}\Phi(\CX^{\rm CS}_t)]$ in $t>0$, we fix $t>0$ and a truncation parameter $\eps\in (0, t)$. For each $s\in (t-\eps, t+\eps)$, we decompose $L^{\rm CS}_{\rm r}\Phi(\CX^{\rm CS}_s)$ into $L_{\rm r}^{<\eps+s-t}\Phi(\CX^{\rm CS}_s)$ and $L_{\rm r}^{\geq\eps+s-t}\Phi(\CX^{\rm CS}_s)$, where both $L_{\rm r}^{<\eps+s-t}\Phi(\CX^{\rm CS}_s)$ and $L_{\rm r}^{\geq\eps+s-t}\Phi(\CX^{\rm CS}_s)$ are defined as in (\ref{L3Phil0}), except that resampling therein is carried out by summing over
\be{rs31}
\begin{aligned}
& \sum_{x_i, x_j \in \{y^+, y^-: y \in E_s^{\eps+s-t}\}} & 1_{\{x_i=x_j\}}
&\  \mbox{ for} \quad L_{\rm r}^{\geq\eps+s-t}\Phi(\CX^{\rm CS}_s), \\
& \sum_{x_i, x_j \in \{y^+, y^-: y \in E_s\backslash E_s^{\eps+s-t}\}} & 1_{\{x_i=x_j\}}
&\ \mbox{ for} \quad L_{\rm r}^{<\eps+s-t}\Phi(\CX^{\rm CS}_s),
\end{aligned}
\ee
where $E_s$ and $E_s^\delta$ are defined as in (\ref{Elt22}).

The same argument as in the proof of the continuity of $\E[L^{\rm CS}_{\rm r}\Phi(\CX^{\rm CS}_t)]$ in
Prop.~\ref{P:genaction} shows that
\be{LPhicont1}
\lim_{s\to t}\E[|L_{\rm r}^{\geq \eps+s-t}\Phi(\CX^{\rm CS}_s)-L_{\rm r}^{\geq \eps+s-t}\Phi(\CX^{\rm CS}_t)|]=0.
\ee
On the other hand, for $s\in [t-\frac{\eps}{2}, t+\frac{\eps}{2}]$, by the same calculations as those leading to (\ref{dxint}), we have
\begin{eqnarray}
&& \E[|L_{\rm r}^{<\eps+s-t}\Phi(\CX^{\rm CS}_s)|] \ \leq\ C_{n,g,\varphi}
\E\Big[\sum_{x\in (-M,M)\cap  E_s\backslash E_s^{\eps+s-t}}d_x\wedge 1\Big] \nonumber \\
&=& C_{n,g,\varphi} \E\Big[\sum_{x\in (-M,M)\cap E_s \atop d_x \leq 2(\eps+s-t)}d_x\wedge 1\Big]
\ \leq\ CM \sqrt{\eps+s-t} \ \leq\ C'M\sqrt{\eps}. \label{LPhicont2}
\end{eqnarray}
Since $\eps>0$ can be chosen arbitrarily small, (\ref{LPhicont1}) and (\ref{LPhicont2}) together imply that
$$
\lim_{s\to t} \E[L^{\rm CS}_{\rm r}\Phi(\CX^{\rm CS}_s)]=\E[L^{\rm CS}_{\rm r}\Phi(\CX^{\rm CS}_t)] \qquad \mbox{when } t>0.
$$

{\bf (b)} We now verify the continuity of $\E[L^{\rm CS}_{\rm r}\Phi(\CX^{\rm CS}_t)]$ (and hence of $\E[L^{\rm CS}\Phi(\CX^{\rm CS}_t)]$) at $t=0$ under the additional assumption that $\CX^{\rm CS}_0 \in \U^\R_{\rm rr}$. Together with (\ref{inteqn1}), this also implies that the generator
equation (\ref{generator1}) holds for a general polynomial $\Phi\in \Pi^{1,2}$.

As before, we separate $L^{\rm CS}_{\rm r}\Phi(\CX^{\rm CS}_s)$ into $L_{\rm r}^{\geq \eps+s}\Phi(\CX^{\rm CS}_s)$ and $L_{\rm r}^{<\eps-s}\Phi(\CX^{\rm CS}_s)$, where the truncation parameter $\eps>0$ is fixed and small.
Equation (\ref{LPhicont1}) also holds with $t=0$ by the same argument as that for the continuity of $\E[L^{\rm CS}_{\rm r}\Phi(\CX^{\rm CS}_t)]$ at $t=0$, when $\Phi \in \Pi^{1,2}_{\rm r}$. Indeed, in both cases, only resampling between individuals with sufficiently large genealogical distance contribute to the generator action on $\Phi$.

It remains to show that
\be{LPhicont3}
\sup_{0\leq s\leq \eps/2}\E[|L_{\rm r}^{<\eps+s}\Phi(\CX^{\rm CS}_s)|]\underset{\eps\downarrow 0}{\longrightarrow} 0.
\ee

For $0<s<\eps/2$, we can separate the resampling terms according to whether the genealogies of the two resampled individuals merge above or below time $0$, and write
\begin{eqnarray}
\E[|L_{\rm r}^{<\eps+s}\Phi(\CX^{\rm CS}_s)|] &\leq& C_{n,g,\varphi}
\E\Big[\sum_{x\in (-M,M)\cap  E_s \atop d_x \leq 2(\eps+s)}d_x\Big]
=
C_{n,g,\varphi}\E\Big[\sum_{x\in (-M,M)\cap E_s \atop d_x \in (0, 2s]}d_x +
\sum_{x\in (-M,M)\cap  E_s \atop d_x \in (2s, 2(\eps+s))} d_x\Big] \nonumber \\
&\leq& C_{n,g,\varphi} M \sqrt{s}
+ C_{n,g,\varphi}\,\E\Big[\sum_{x\in (-M,M)\cap E_s \atop d_x \in (2s, 2(\eps+s))} d_x\Big]. \label{LPhicont4}
\end{eqnarray}
It remains to bound the expectation on the r.h.s. above, which originate from resampling between individuals whose genealogical distance depend on
the distance of their ancestors in $\CX^{\rm CS}_0$ at time $0$. Let $\hat \xi^{[-M,M]\times\{s\}}_0$ be the point set on $\R$ generated by the collection of paths in the dual Brownian web
$\widehat\CW$ starting from the space-time set $[-M, M]\times\{s\}$. Let us order the points in $\hat \xi^{[-M,M]\times\{s\}}_0$
and denote them by $z_1<z_2\cdots <z_{M_s+1}$, with a total of $M_s+1$ points. Then by our construction of $\CX^{\rm CS}_s$, we have
\begin{eqnarray}
\E\Big[\sum_{x\in (-M, M)\cap E_s \atop d_x \in (2s, 2(\eps+s))} d_x\Big]
&=& \E\Big[\sum_{i=1}^{M_s} 1_{\{r^{\rm CS}_0(z_i,z_{i+1})<2\eps\}} (2s+r^{\rm CS}_0(z_i,z_{i+1}))\Big]
\nonumber \\
&\leq& CM \sqrt{s} + \E\Big[\sum_{i=1}^{M_s} 1_{\{r^{\rm CS}_0(z_i,z_{i+1})<2\eps\}} r^{\rm CS}_0(z_i,z_{i+1})\Big],
\label{LPh5}
\end{eqnarray}
where we used $\E[M_s] = \frac{2M}{\sqrt{\pi s}}$ because by duality between the Brownian web and its dual, $M_s$ is exactly the number of points in the interval $[-M, M]$ occupied at time $s$ by Brownian web paths starting from $\R$ at time $0$, and hence Lemma~\ref{L:density1} can be applied. For each pair $(z_i,z_{i+1})$ with
$r^{\rm CS}_0(z_i,z_{i+1})<2\eps$, by the properties of $\CX^{\rm CS}_0$, there must exist a $y_i\in E_0$
with $z_i<y_i<z_{i+1}$ such that $r^{\rm CS}_0(z_i,z_{i+1}) = d_{y_i}$. Note that the $y_i$, $1\leq i\leq M_s$,
are all distinct with $d_{y_i}<2\eps$. We can now separate the contribution to the second term in (\ref{LPh5}) into
two groups.

The first group consists of contributions from pairs $(z_i,z_{i+1})$ with $-2M\leq z_i<z_{i+1}\leq 2M$.
The total contribution from these terms is uniformly dominated by
$\sum_{x\in [-2M,2M]\cap E_0\atop d_x\leq 2\eps}d_x$, which tends to $0$ as $\eps\downarrow 0$ by the assumption $\CX^{\rm CS}_0\in \U^\R_{\rm rr}$.

The second group consists of contributions from pairs $(z_i,z_{i+1})$ with either $z_i\leq -2M$ or
$z_{i+1}\geq 2M$. We bound these terms by $2\eps$ times the expected cardinality of such pairs. By the
duality between $\CW$ and $\widehat\CW$, the cardinality of such pairs of $(z_i,z_{i+1})$ is bounded by the cardinality
of $\xi_s^{\{x\in\R :|x|\geq 2M\}\times\{0\}} \cap (-M,M)$, the point set on $(-M,M)$ generated at time $s$
by paths in $\CW$ starting from the space-time set$\{x\in\R: |x|\geq 2M\}\times\{0\}$. As shown in (\ref{Zsplit2}),
the expected cardinality of $\xi_s^{\{x\in\R :|x|\geq 2M\}\times\{0\}} \cap (-M,M)$ is uniformly bounded in $s\in [0,l]$
for any $l>0$. Therefore, uniformly for $s\in [0,\eps/2]$, the second term in (\ref{LPh5}) is also bounded by a
function of $\eps$, which tends to $0$ as $\eps\downarrow 0$. This concludes the proof of (\ref{LPhicont3}) and
Proposition \ref{P:genaction2}.
\qed

\subsection{Proof of Theorem~\ref{T.mp}}

\label{ss:martingales}

We now prove Theorem~\ref{T.mp}. Note that the various path properties listed in Theorem~\ref{T.mp} has been verified for the CSSM genealogy process $\CX^{\rm CS}$ in Prop.~\ref{P:reg} and Theorem~\ref{T:MarkovCont}, except for the claim that for each $t>0$, $\CX^{\rm CS}_t\in \U^\R_{\rm rr}$ almost surely. To verify this last claim, we need to show that a.s.,
\begin{equation}\label{Urrbound1}
\sum_{x\in E_t\cap (-M, M)} r^{\rm CS}_t(x^+, x^-) <\infty \qquad \mbox{for each } M>0.
\end{equation}
where $E_t$ is defined as in (\ref{Elt22}). Since for each $\epsilon>0$, the set $\{x\in E_t \cap (-M, M) : r^{\rm CS}_t(x^+, x^-)\geq \eps\}$ is a finite set because $\CX^{\rm CS}_t \in \U^\R_{\rm r}$, (\ref{Urrbound1}) is reduced to showing that a.s.,
$$
\sum_{x\in E_t\cap (-M, M)} 1_{\{r^{\rm CS}_t(x^+, x^-) <\eps\}} r^{\rm CS}_t(x^+, x^-) <\infty \qquad \mbox{for each } M>0.
$$
The l.h.s.\ above has been shown in (\ref{dxint}) to have finite expectation, and hence is finite a.s. Therefore $\CX^{\rm CS}_t\in \U^\R_{\rm rr}$ a.s.

To conclude the proof of Theorem~\ref{T.mp}, it only remains to verify the martingale property, namely that given $\Phi \in \Pi^{1,2}$ and $\CX^{\rm CS}_0\in \U^\R_1$, for each $0\leq s<t$, we have a.s.\
\begin{equation}\label{marting}
\mathbb{E}\Big[\Phi(\CX^{\rm CS}_t)  - \Phi(\CX^{\rm CS}_s) - \int_s^t (L^{\rm CS}\Phi)(\CX^{\rm CS}_u)\dd u\,\Big|\, (\CX^{\rm CS}_u)_{0\leq u\leq s}\Big] =0.
\end{equation}

Consider first the case $s>0$. Then a.s.\ $\CX^{\rm CS}_s\in \U^\R_{\rm r}$. We can use the Markov property of $\CX^{\rm CS}$ and apply the integral equation (\ref{inteqn1}) for $\Phi\in \Pi^{1,2}$ and initial state $\CX^{\rm CS}_s$, as established in Prop.~\ref{P:genaction2}. Equation (\ref{marting}) then follows provided that we can apply Fubini's Theorem to interchange the integral with expectation in
\begin{equation}\label{fubini}
 \mathbb{E}\Big[\int_s^t (L^{\rm CS}\Phi)(\CX^{\rm CS}_u)\dd u\Big] = \int_s^t \E[(L^{\rm CS}\Phi)(\CX^{\rm CS}_u)]\dd u.
\end{equation}
As in the proof of Prop.~\ref{P:genaction2}, we can decompose $L^{\rm CS}$ as $L^{\rm CS}_{\rm a} + L^{\rm CS}_{\rm d} + L^{\rm CS}_{\rm r}$. The aging and diffusion part of the generator action on $\Phi$ can be uniformly bounded in time and in $\CX^{\rm CS}_u$, while the resampling part can be bounded as in (\ref{LRXsbd}), which is shown to be integrable w.r.t.\ $E[\int_0^t\cdot ]$ in (\ref{dxfinite}). Therefore Fubini can be applied, and (\ref{marting}) holds for $s>0$ a.s.

To treat the case $s=0$, we take expectation in (\ref{marting}) for $s>0$, i.e.,
$$
\mathbb{E}\Big[\Phi(\CX^{\rm CS}_t)  - \Phi(\CX^{\rm CS}_s) - \int_s^t (L^{\rm CS}\Phi)(\CX^{\rm CS}_u)\dd u\Big]=0,
$$
and then let $s\downarrow 0$. By the a.s.\ continuity of $\CX^{\rm CS}_s$ in $s\geq 0$, we have $\E[\Phi(\CX^{\rm CS}_s)] \to \Phi(\CX^{\rm CS}_0)$. The convergence of $\\mathbb{E}\Big[\int_s^t (L^{\rm CS}\Phi)(\CX^{\rm CS}_u)\dd u\big]$ as $s\downarrow 0$ follows by dominated convergence, using the same argument as that for (\ref{fubini}). Therefore (\ref{marting}) also holds for $s=0$, which proves the desired martingale property.
\qed

\appendix

\section{Proofs on marked metric measure spaces}\label{S:mmm}
In this section, we prove Theorems \ref{T:mmmPolish} and \ref{T:polydet}, derive a relative compactness condition for subsets of $\M^V$, and a tightness criterion for
laws on $\M^V$. Furthermore we formally treat pasting of trees. Our starting point for the first of these points is Remark~\ref{R:GromovHash}, that $\M^V$ can be identified
as a subspace of $(\M^V_f)^{\N}$, endowed with the product $V$-marked Gromov-weak topology.
\medskip

\noindent
{\bf Proof of Theorem~\ref{T:mmmPolish}.} As noted in Remark~\ref{R:GromovHash}, we can identify $\M^V$ as a subspace of $(\M^V_f)^{\N}$, endowed with the product $V$-marked Gromov-weak topology. Furthermore, under this identification, $\M^V$ is a closed subspace of $(\M^V_f)^{\N}$. It was shown in \cite[Theorem~2]{DGP11} that $\M^V_1$, the space of $V$-mmm spaces with probability measures, equipped with the $V$-marked Gromov-weak topology, is a Polish space. The same conclusion is easily seen to hold for $\M^V_f$. Therefore $(\M^V_f)^{\N}$ is also Polish, which implies that any closed subspace, including $\M^V$, is also Polish.
\qed
\bigskip


\noindent

{\bf Proof of Theorem~\ref{T:polydet}.} For each $k\in \N\cup\{0,\infty\}$, let
$$
\widetilde \Pi^k := \cup_{n\in\N_0} \widetilde \Pi^k_n := \cup_{n\in\N_0} \big\{ \Phi^{n,\phi} : \phi \in C^k( \R_+^{\binom{n}{2}} \times V^n, \R)\big\},
$$
where $C^k( \R_+^{\binom{n}{2}}\times V^n, \R)$ is the space of bounded continuous real-valued functions on $\R_+^{\binom{n}{2}}\times V^n$ that are $k$ times continuously differentiable in the first $\binom{n}{2}$ coordinates, and
$$
\Phi^{n,\phi} (\overline{(X, r, \mu)}) := \idotsint  \phi(\uur, \uv) \mu^{\otimes n}(d(\ux, \uv)) \qquad \mbox{for each } \overline{(X, r, \mu)}\in \M^V_1,
$$
with $\uur := (r(x_i, x_j))_{1\leq i<j\leq n}$ and $\uv:=(v_1, \ldots, v_n)$.

{We use once more Remark~\ref{R:GromovHash} to now identify {$\M^{(V,\le b)}$} as a subspace of
{$\prod_{k\in\mathbb{N}}\M^{V}_{\le c_k}$}, endowed with the product $V$-marked Gromov-weak topology.
Here {$\M^{V}_{\le c}$} denotes the  space of $V$-mmm spaces with a measure of total mass {at most} $c$, and e.g.\ $c_k:=\|\psi_k\|_\infty b(r)$ for any $r>0$ with $\mathrm{supp}(\psi_k)\subseteq B_r(o)$.}

By \cite[Theorem~5]{DGP11}, $\widetilde \Pi^k := \cup_n \widetilde \Pi^k_n$ is convergence determining in $\M^V_1$ and $\CM_1(\M_1^V)$,
and hence also in {$\M^V_{\le c}$ and $\CM_1(\M^V_{\le c})$}. We now want to argue that this holds as well for
$\bigcup_n \Pi^k \subseteq \bigcup_n \wt \Pi^k$. This follows immediately from \cite[Prop.~3.4.6]{EK86} for measures on product spaces.
\qed

\medskip

The following relative compactness criterion for subsets of $\M^V$ follow easily from the identification of $\M^V$ as a subspace of $(\M^V_f)^{\N}$, and the relative compactness criterion {formulated for subsets of $\mathbb{M}^{\{o\}}_1$ in~\cite[Prop.~6.1]{GPWmp13}, for subsets of
$\M^V_1$ in~\cite[Thm.~3]{DGP11} and for subsets in $\mathbb{M}^{\{o\}}$~\cite[Corollary~4.3]{AthreyaLohrWinter}.}
\bt[Relative compactness of subsets of $\M^V$]\label{T:mmmrc} Let $\Gamma \subset \M^V$, and let $o$ be any point in $V$. For each $k\in\N$, let $B_k(o)$ denote the open ball of radius $k$ centered at $o$. Then $\Gamma$ is relatively compact w.r.t.\ the $V$-marked Gromov-weak$^\#$ topology if for each $k\in\N$, $\{\overline{(X,r, 1_{\{v\in B_k(o)\}}\mu(\dd x\dd v) ) }: \overline{(X,r,\mu)} \in \Gamma\}$ is a relatively compact subset of $\M^V_f$, i.e.,
\begin{itemize}
\item[\rm (i)]  The family of finite measures on $V$,
$$
\Lambda_k := \{\mu (X\times \dd v)1_{\{v\in B_k(o)\}} : \overline{(X,r, \mu)}\in \Gamma\},
$$
is relatively compact w.r.t.\ the weak topology;

\item[\rm (ii)] For each $\eps>0$, there exists $L>0$ such that uniformly in $\overline{(X,r, \mu)}\in \Gamma$,
$$
\iint_{(X\times V)^2} 1_{\{ r(x,y)>L  \}} 1_{\{ u,v\in B_k(o) \}} \mu(\dd x\dd u) \mu(\dd y \dd v) \leq \eps;
$$

\item[\rm (iii)] For each $\eps>0$, there exists $M\in\N$ such that uniformly in $\overline{(X,r, \mu)}\in \Gamma$,
 we can find $M$ balls of radius $\eps$ in $X$, say $B_{\ve,1}, \ldots, B_{\ve,M} \subset X$ with $B=\cup_{i=1}^M B_{\ve,i}$,
 such that $\mu( (X\backslash B)\times B_k(o)) \leq \eps$.
\end{itemize}
\et

By the tightness criteria formulated {for random variables with values in $\mathbb{M}^{\{o\}}_1$ in~\cite[Theorem~3]{GPWmp13}, with values in
$\M^V_1$ in~\cite[Theorem~4]{DGP11} and with values  in $\mathbb{M}^{\{o\}}$ in~\cite[Corollary~4.6]{AthreyaLohrWinter},}
we obtain the following tightness criteria for $\M^V$-valued random variables.
\bt[Tightness of $\M^V$-valued random variables]\label{T:mmmProperties} Let $o$ be any point in $V$, and let $B_k(o)$ denote the open ball of radius $k$ centered at $o$.
A family of $\M^V$-valued random variables $\{\CX_i=\overline{(X_i, r_i, \mu_i)}\}_{i\in I}$ is tight if for each $k\in\N$, $\{\overline{X_i, r_i, 1_{\{v\in B_k(o)\}}\mu_i(\dd x\dd v)}\}_{i\in I}$ is a tight family of $\M^V_f$-valued random variables, i.e.,
\begin{itemize}
\item[\rm (i)] $\{1_{\{v\in B_k(o)\}}\mu_i(X\times \dd v) \}_{i\in I}$ is a tight family of random variables taking values in the space of finite measures on $V$ (equipped with the weak topology);

\item[\rm (ii)] $\{\overline{X_i, r_i, \mu_i(\dd x\times B_k(o))}\}_{i\in I}$ is a tight family of random variables taking values in the space of metric measure spaces (equipped with the Gromov-weak topology).
\end{itemize}
\et
Using the characterization of relatively compact sets in $\M^V$ given in Theorem~\ref{T:mmmrc}, one can also formulate more concrete conditions for the tightness of a family of $\M^V$-valued random variables, using concrete conditions for the tightness of a family of random metric measure spaces formulated in~\cite[Thm.~3]{GPW09}.
\bi

\section{The Brownian web}\label{S:web}
In this section, we recall the construction and basic properties of the Brownian web. For a recent comprehensive survey, see~\cite{SchertzerSunSwart}. Recall from Subsection \ref{S:IntroCSSM} the random variable $(\CW, \wh \CW)$ as constructed in \cite{FINR04, FINR06}, and in particular from (\ref{gs9b})-(\ref{gs9d}), the state spaces of $\Pi$ of $\CW$ and $\wh \Pi$ of $\wh \CW$. It has been shown~\cite[Theorem~2.1]{FINR04} that the Brownian web $\CW$ can be characterized as follows:

\bt[Characterization of the Brownian web]\label{T:webchar}

The Brownian web $\CW$ is a random closed subset of $\Pi$, whose law is uniquely determined by the following properties:
\begin{itemize}
\item[\rm (i)] For every $z \in \R^2$, almost surely $\CW(z)$ contains a unique path.

\item[\rm (ii)] For every finite $n$ and deterministic points $z_1, \cdots, z_n\in\R^2$, $ \{\CW(z_i) : i=1,\cdots,n\}$ are distributed as $n$ coalescing Brownian motions starting from $z_1, \cdots, z_n$.

\item[\rm (iii)] For every deterministic countable dense subset $\CD\subset \R^2$, $\CW$ is almost surely the closure of $\{\CW(z); z \in \CD\}$ in $\Pi$.
\end{itemize}
\et
The following result shows that every path in $\CW$ can be approximated by the countable set of paths $\{\CW(z); z \in \CD\}$ in a very strong sense (see e.g.~\cite[Lemma 3.4]{SS08}).

\bl[Convergence of Paths in $\CW$]\label{L:pathconv}
Almost surely, if $(f_n)_{n\in\N}$ and $f$ are paths in $\CW$, starting respectively at times $(s_n)_{n\in\N}$ and $s$, and $f_n\to f$ in $\Pi$, then $\sup\{ t : f_n(t)\neq f(t)\} \to s$ as $n\to\infty$.
\el

The dual Brownian web $\wh\CW$ can be characterized as follows~\cite[Theorem 3.7]{FINR06}:

\bt[Characterization of the dual Brownian web]\label{T:dwebchar}
Let $\CW$ be a Brownian web. Then there exists an almost surely uniquely determined $\wh\Pi$-valued random variable $\wh\CW$ defined on the same probability space as $\CW$, called the dual Brownian web, such that:
\begin{itemize}
\item[\rm (i)] Almost surely, paths in $\CW$ and $\wh\CW$ do not cross, i.e., there exist no $f\in \CW$ and $\hat f\in \wh\CW$ and $s\neq t$ such that $(f(s)-\hat f(s))(f(t)-\hat f(t))<0$;

\item[\rm (ii)] ${\CR}\wh\CW$ has the same law as $\CW$, where ${\CR}$ denotes the reflection map that maps each $\hat f\in \widehat\CW$ to an $f\in\Pi$ such that the graph of $f$ in $\R^2$ is the reflection of the graph of $\hat f$ with respect to the origin.
\end{itemize}
\et

Below we collect some basic properties of the Brownain web which we will use. For further details, see~\cite{FINR04, FINR06, SS08}. The first property concerns the configuration of paths in $\CW$ and $\wh\CW$ entering and leaving a point $z\in \R^2$. For each $z=(x,t)\in \R^2$, let $m_{\rm out}(z)$ denote the cardinality of $\CW(z)$. We will let $m_{\rm in}(z)$ denote the number of equivalence classes of paths in $\CW$ entering $z$, where a path $f\in \CW$ is said to enter $z$ if it starts before time $t$ and $f(t)=x$, while two paths $f,g\in\CW$ entering $z$ are called equivalent if they coalesce before time $t$. Note that $m_{\rm in}(z)=2$ if and only if $z$ is a point of coalescence between two paths in $\CW$.
Similarly, we can define $\hat m_{\rm in}(z)$ and $\hat m_{\rm out}(z)$, based on the configuration of paths in $\wh\CW$. The pair $(m_{\rm in}, m_{\rm out})$
is called the type of $z$ in $\CW$.

We cite the following result from \cite[Prop.~2.4]{TW98} and \cite[Thm.~3.11-3.14]{FINR06}.

\bl[Special points for the Brownian web]\label{L:spts}
Let $\CW$ and $\wh\CW$ be a Brownian web and its dual. Almost surely:
\begin{itemize}
\item[\rm (1)] The set of $z\in\R^2$ with $(m_{\rm in}(z), m_{\rm out}(z))= (\hat m_{\rm in}(z), \hat m_{\rm out}(z)) = (0, 1)$ has full
Lebesgue measure on $\R^2$.

\item[\rm (2)] For each $t\in\R$, the set of  {$z=(x',t')\in \R\times \{t\}$} with $\hat m_{\rm out}(z)\geq 2$ is a countable set, with $m_{\rm in}(z)\geq 1$ for each such $z$, i.e., $z$ lies on the graph of some path in $\CW$ starting before time $t$.
\end{itemize}
\el

Next we cite a result on the decay of the density of coalescing paths started at time $0$.

\bl[Density for the Brownian web]\label{L:density1}
For $t>0$, let $\xi^{\R\times \{0\}}_t:=\{f(t): f\in \cup_{x\in\R}\CW(x,0) \}$ denote the point set on $\R$ generated at time $t$ by the collection of coalescing paths in the Brownian web $\CW$ started at time $0$. Then for any $a<b$,
\be{webden}
\E[|\xi^{\R\times \{0\}}_t \cap [a,b]|] = \frac{b-a}{\sqrt{\pi t}}.
\ee
\el
This result can be easily derived by using the duality between $\CW$ and $\wh\CW$, namely that {$\xi^{\R\times\{0\}}_t\cap (x,x+\eps)\neq \emptyset$} if and only if the two paths in $\wh\CW$ starting from $(x,t)$ and $(x+\eps,t)$ do not collide on the time interval $[0,t]$. See e.g.~\cite[Prop.~1.12]{SS08}, where such a density calculation is carried out for a generalization of the Brownian web known as the {\em Brownian net}, which in addition allows branching of paths.

\section{Correlation Inequalities for Coalescing Brownian motions}\label{A:Corr}
In this section, we prove some negative correlation inequalities for a collection of coalescing Brownian motions, which are used in Section~\ref{S:CSSMmart}. {}Similar inequalities have previously been established in \cite{MRTZ06}, \color{black} see also \cite[Remark 7.5]{NRS05}. Here we deduce more general negative correlation inequalities from {\em Reimer's inequality} applied to coalescing random walks.

In van den Berg and Kesten \cite{BK02}, Reimer's inequality was applied to continuous time coalescing random walks with a generalized coalescing rule. Since we are interested in coalescing Brownian motions, discrete space-time coalescing random walks with instantaneous coalescing already provide an adequate approximation, and Reimer's inequality can be applied without any complication
to the latter.

First we recall Reimer's inequality \cite{R00}. For each $i\in I:=\{1,\cdots,n\}$, let $S_i$ be a finite set
with a probability measure $\mu_i$ on $S_i$. Let $\Omega = S_1\times S_2\cdots\times S_n$ and
$\mu=\mu_1\times\cdots\times\mu_n$. For $K\subset I$ and $\omega=(\omega_i)_{i\in I}$, define the cylinder
set $C(K, \omega) :=\{\omega'\in\Omega : \omega'_i = \omega_i \ \forall\ i\in K\}$. Given two events
$A, B\subset\Omega$, we say {\em $A$ and $B$ occur disjointly} for a configuration $\omega\in\Omega$ if there
exists $K\subset I$ such that $C(K,\omega) \subset A$ and $C(I\backslash K, \omega) \subset B$. The set
of $\omega\in\Omega$ for which $A$ and $B$ occur disjointly, which we call the disjoint intersection of $A$ and $B$,
is denoted by
$$
A \square B := \{\omega \in \Omega : \exists\, K\subset I\ \mbox{s.t.}\ C(K, \omega)\subset A \ \mbox{and}\
C(I\backslash K, \omega)\subset B\}.
$$
Then Reimer's inequality asserts that, for any two events $A, B\subset \Omega$,
\be{Reimer}
\mu(A\square B) \leq \mu(A) \mu(B).
\ee

Now we apply this inequality to coalescing random walks. We
recall first the construction of discrete space-time coalescing random walks. Let
$\Z^2_{\rm even}=\{(x,t)\in\Z^2: x+t\ \mbox{is even}\}$. Let $\{\omega_z\}_{z\in\Z^2_{\rm even}}$ be i.i.d.\
random variables taking values in $\{\pm1\}$. A directed edge is drawn from each $z=(x,t)\in\Z^2_{\rm even}$,
which ends at $(x+1,t+1)$ if $\omega_z=1$, and ends at $(x-1, t+1)$ if $\omega_z=-1$. This provides a graphical
construction of a collection of coalescing random walks, where the random walk path starting from each
$z\in\Z^2_{\rm even}$ is constructed by following the directed edges in $\Z^2_{\rm even}$ drawn according to
$\omega$.

To {}illustrate \color{black} how Reimer's inequality is applied, let $(X^{(x_i,t_i)}_j)_{j\geq t_i}$, $1\leq i\leq n$, be a collection
of coalescing random walks constructed as above with starting points $(x_i,t_i)\in\Z^2_{\rm even}$, and assume
for simplicity $t_i\leq 0$ for all $1\leq i\leq n$. For $t\in\N$, let
$\xi_t = \{x\in\Z : x=X^{(x_i,t_i)}_t\ \mbox{for some}\ 1\leq i\leq n\}$. Let $O_1,\cdots, O_k$ be disjoint
subsets of $\Z$, and let $A_i =\{\omega : \xi_t \cap O_i\neq \emptyset\}$. The crucial observation is that,
if the events $(A_i)_{1\leq i\leq n}$ occur simultaneously, then they must occur disjointly w.r.t.\
$(\omega_z)_{z\in\Z^2_{\rm even}}$ because of the coalescence. Namely,

\be{Reimer2}
\bigcap_{i=1}^k A_i = A_1 \square A_2\cdots \square A_k.
\ee
Reimer's inequality (\ref{Reimer}) then gives the negative correlation inequality

\be{Reimer3}
\P\big(\bigcap_{i=1}^k A_i\big) \leq \prod_{i=1}^k \P(A_i).
\ee
The same reasoning allows us to choose each $A_i$ to be an increasing event of the occupation configuration
$\xi_t\cap O_i$, i.e., given $\omega$ and $\omega'$ with respective occupation configurations
$\xi_t\cap O_i\subset \xi'_t\cap O_i$, if $\omega \in A_i$, then also $\omega'\in A_i$.
{}Note that Reimer's inequality may even be applied to the disjoint occurence of the event $ A_1 $ with itself, which we use later in the proof of Lemma \ref{L:c2point}. \color{black}

Using Reimer's inequality as illustrated above, together with the invariance principle for coalescing random walks,
we will deduce a host of negative correlation inequalities for coalescing Brownian motions, which we formulate next.

\bd[Negatively correlated point processes]\label{D:negcor}
We say a simple point process $\xi$ on $\R$ is negatively correlated, if for any $n\in\N$ and any disjoint open
intervals $O_1,\cdots, O_n$, we have
$\P(\cap_{i=1}^n \{\xi\cap O_i\neq\emptyset\}) \leq \prod_{i=1}^n \P(\xi\cap O_i\neq\emptyset)$.
\ed

\bl[Negative correlation for colaescing Brownian motions]\label{L:negcor}
Let $A\subset \R\times (-\infty,0]$, and let $\xi^A_t$ denote the point set on $\R$ generated at time $t>0$ by
the collection of coalescing Brownian motions in the Brownian web $\CW$ starting from $A$. Then $\xi^A_t$ is
negatively correlated.
\el
\noindent
{\bf Proof.}
By monotone convergence, it suffices to consider the case when $A$ consists of a finite
number of points $\{z_1,\cdots, z_k\}$. The fact that $\xi^A_t$ is negatively correlated then follows directly
from  the negative correlation inequality (\ref{Reimer3}) for coalescing random walks, and the distributional
convergence of coalescing random walks to coalescing Brownian motions in the local uniform topology (see e.g.~\cite[Section 5]{NRS05}).
\qed

\bl[Decoupling of correlation functions]\label{L:negcor2}
Let $\xi^A_t$ be as in Lemma \ref{L:negcor}. Let $a_1<b_1\leq a_2<b_2\cdots\leq a_n<b_n$. Then for
any $1\leq j\leq n-1$,
\be{corsplit}
\begin{aligned}
& \P\Big(\bigcap_{i=1}^n \{\xi^A_t \cap (a_i,b_i)\neq\emptyset\} \cap\{\xi^A_t \cap (b_{j},a_{j+1})=\emptyset\} \Big) \\
\leq\ &
\P\Big(\xi^A_t \cap (a_j, b_j)\neq \emptyset, \xi^A_t \cap (b_j, a_{j+1})= \emptyset, \xi^A_t \cap (a_{j+1}, b_{j+1})\neq \emptyset\Big)
\!\!\!\prod_{1\leq i\leq n \atop i\neq j, j+1}\!\!\!\! \P(\xi^A_t \cap (a_i, b_i)\neq\emptyset).
\end{aligned}
\ee
\el
\noindent
{\bf Proof.} As before, this follows from approximation by discrete space-time coalescing random walks and Reimer's
inequality. Note that for the discrete analogue of the events in the second line of (\ref{corsplit}), if they all
occur, then they must occur disjointly.
\qed

\bl[Negative correlation for occupation number]\label{L:ocucor}
Let $\xi^A_t$ be as in Lemma \ref{L:negcor}, and let $B\subset \R$ have finite Lebesgue measure. Then for any $k\in\N$,
\be{ocuineq}
\P(|\xi^A_t\cap B|\geq k) \leq \P(|\xi^A_t\cap B|\geq 1)^k.
\ee
\el
\noindent
{\bf Proof.}
By monotone convergence and approximation by open sets, it suffices to consider the case
when $A\subset \R\times(-\infty, 0]$ consists of a finite number of points, and $B$ is the finite union of disjoint
open intervals. In fact, the argument is the same if $B$ is a bounded open interval, say $(0,1)$.

We proceed by discrete approximation. Let $(Z_t)_{t\in\N}$ be the subset of $\Z$ occupied at time $t\in\N$ by a collection of
coalescing random walks on $\Z^2_{\rm even}$ starting from $z_1=(x_1,t_1),\cdots, z_n=(x_n,t_n)\in\Z^2_{\rm even}$ with
$t_i\leq 0$ for all $1\leq i\leq n$. Given $O\subset \Z$ and for $k\in\N$, let $A_k=\{\omega\in \Omega : |Z_t \cap O| \geq k\}$,
where $\omega=(\omega_z)_{z\in\Z^2_{\rm even}}$ are the i.i.d.\ $\{\pm1\}$-valued random variables underlying the graphical
construction of the coalescing random walks. We note that
\be{kfold}
A_k \subset \overbrace{A_1\square \cdots \square A_1}^{k}.
\ee
Indeed, if $A_k$ occurs, then we can find $k$ disjoint random walk paths, each of which occupies a distinct site in $O$
at time $t$. Reimer's inequality (\ref{Reimer}) then implies
$\P(A_k)\leq \P(A_1)^k$. Inequality (\ref{ocuineq}) then follows by the distributional convergence of coalescing
random walks to coalescing Brownian motions in the local uniform topology.
\qed

\bl[Moment bounds for occupation number]\label{L:negmom}
Let $\xi$ be a simple point process on $\R$ with a locally finite intensity measure $\mu$, which is absolutely continuous w.r.t.\
Lebesgue measure on $\R$. If $\xi$ is negatively correlated, then for any Lebesgue measurable $B\subset \R$ with
$\mu(B)<\infty$, and for any $k\in\N$, we have
\be{negmom}
  \mathbb{E}\Big[|\xi\cap B|^k\big] \leq  \sum_{m=1}^k {k\choose m} m^{k-m} \mu(B)^m.
\ee
\el
\noindent
{\bf Proof.}
Let $B$ be an open interval, say $(0,1)$. For $n\in\N$, let $D_n=\{i2^{-n} : 0\leq i\leq 2^n\}$, and
let $D=\bigcup_{n\in\N} D_n$. By our assumption that the intensity measure $\mu$ is absolutely continuous w.r.t.\
Lebesgue measure, $\xi\cap D=\emptyset$ almost surely. For $1\leq i\leq 2^n$, let
$I^{(n)}_i(\xi) = 1_{\{\xi\cap ((i-1)2^{-n}, i2^{-n})\neq\emptyset\}}$. By the assumption that
$\xi$ is a simple point process, and by monotone convergence,
\be{sunx}
\E\big[|\xi\cap (0,1)|^k\big] = \E\Big[\lim_{n\to\infty} \big(\!\!\!\sum_{1\leq i\leq 2^n}\!\!\! I^{(n)}_i\big)^k\Big]
=\lim_{n\to\infty} \E\Big[ \big(\!\!\!\sum_{1\leq i\leq 2^n}\!\!\! I^{(n)}_i\big)^k\Big]
=\lim_{n\to\infty} \!\!\!\sum_{1\leq i_1,\cdots, i_k\leq 2^n} \!\!\!\!\!\!\E\Big[\prod_{j=1}^{k}I^{(n)}_{i_j}\Big].
\ee
Note that $\mu(0,1)=\E[|\xi\cap (0,1)|]= \lim_{n\to\infty} \sum_{j=1}^{2^n}\E[I^{(n)}_j]$.
Given $\vec i:=(i_1,\cdots, i_k)\in \{1,\cdots, 2^n\}^k$, let $m(\vec i)=|\{i_1,\cdots, i_k\}|$, and denote the
$m(\vec i)$ distinct indices in $\{i_1,\cdots, i_k\}$ by $\sigma_1(\vec i)<\cdots< \sigma_m(\vec i)$. Then by
the negative correlation assumption,
\be{negmom2}\begin{array}{l}
\suml_{1\leq i_1,\cdots, i_k\leq 2^n} \mathbb{E}\Big[\prod_{j=1}^{k}I^{(n)}_{i_j}\Big]\\[2ex]
\hspace{1.5cm}
\leq  \sum_{1\leq i_1,\cdots, i_k\leq 2^n} \prod_{j=1}^{m(\vec i)}\E[I^{(n)}_{\sigma_j(\vec i)}]
= \sum_{m=1}^k \sum_{1\leq\sigma_1<\cdots <\sigma_m\leq 2^n} f(k,m) m!  \prod_{j=1}^m \E[I^{(n)}_j]
 \\[2ex]
\hspace{1.5cm}
\leq \sum_{m=1}^k f(k,m) \Big(\sum_{j=1}^{2^n} \E[I^{(n)}_j]\Big)^m \underset{n\to\infty}{\longrightarrow}
\sum_{m=1}^k f(k,m) \mu(0,1)^m,
\end{array}
\ee
where $f(k,m) m! = |\{(i_1,\cdots, i_k)\in \{1,\cdots, 2^n\}^k : \{i_1,\cdots, i_k\}=\{\sigma_1,\cdots, \sigma_m\}\}|$,
which is easily seen to be independent of the choice of $1\leq\sigma_1<\cdots<\sigma_m \leq 2^n$. By first picking
$m$ indices out of $\{i_1,\cdots, i_k\}$ and assign them values $\sigma_1<\cdots<\sigma_m$ respectively, we easily
verify that $f(k,m)\leq {k\choose m} m^{k-m}$, which proves (\ref{negmom}) for $B=(0,1)$. By the same argument,
(\ref{negmom}) holds for finite unions of disjoint open intervals, and by monotone convergence, for open sets as well.
Since any Lebesgue measurable set can be approximated from outside by open sets, again by monotone convergence and
the fact $\xi\cap E=\emptyset$ a.s.\ for a given $E\subset\R$ with zero Lebesgue measure, (\ref{negmom})
also holds for any Lebesgue measurable $B$.
\qed
\bigskip

We also need the following estimate on the constrained two-point correlation function for the Brownian web.

\bl[Constrained two point function for the Brownian web]\label{L:c2point}
Let $\xi^{\R\times\{0\}}_t$ be as in Lemma \ref{L:density1}. Let $t>0$. For $a<b$, let $I_{[a,b]}$ denote the event that
$\xi^{\R\times\{0\}}\cap [a,b]\neq\emptyset$. Then for any $x_1<x_2$ with $\Delta:=x_2-x_1$, we have
\be{detform}
K^{\rm c}_t(x_1,x_2) := K^{\rm c}_t(\Delta) = \lim_{\delta\downarrow 0} \frac{1}{\delta^2}\,
\P\big[ I_{[x_1,x_1+\delta]} \cap I^c_{[x_1+\delta, x_2]} \cap I_{[x_2, x_2+\delta]}\big]
= \frac{\Delta e^{-\frac{\Delta^2}{4t}}}{2\sqrt{\pi} t^\frac{3}{2}}. \qad
\ee
\el
\noindent
{\bf Proof.}
By translation invariance of $\xi^{\R\times\{0\}}$, we may assume $x_1=0$, and let $x_2-x_1=\Delta>0$.
Let $\hat w_{(0,t)}$, $\hat w_{(\delta,t)}$, $\hat w_{(\Delta,t)}$ and $\hat w_{(\Delta+\delta,t)}$ be the dual
coalescing Brownian motions in $\widehat\CW$ starting at respectively $(0,t)$, $(\delta, t)$, $(\Delta, t)$ and
$(\Delta+\delta, t)$. Let $\hat\tau_{a,b}$ be the time of coalescence between $\hat w_{(a,t)}$ and $\hat w_{(b,t)}$.
Then by the duality between $\CW$ and $\widehat\CW$, almost surely, the event $I_{[0,\delta]}\cap I_{[\delta,\Delta]}$
occurs if and only if $\hat\tau_{0,\delta}< 0$ and $\hat \tau_{\delta, \Delta}< 0$, and the event
$I_{[0,\delta]}\cap I_{[\delta,\Delta+\delta]}$ occurs if and only if $\hat\tau_{0,\delta}< 0$ and
$\hat \tau_{\delta, \Delta+\delta}< 0$. Since $I_{[\delta,\Delta]}\subset I_{[\delta,\Delta+\delta]}$ and
$I_{[\delta,\Delta+\delta]}\backslash I_{[\delta,\Delta]}= I_{[\delta, \Delta]}^c \cap I_{(\Delta, \Delta+\delta]}$,
we have
\begin{equation}\label{sw1}
\begin{aligned}
& \P\big[ I_{[x_1,x_1+\delta]} \cap I^c_{[x_1+\delta, x_2]} \cap I_{[x_2, x_2+\delta]}\big]  \\
=\ & \P^{B_1, B_2, B_3}_{0,\delta, \Delta}(\tau_{0,\delta}>t, \tau_{\delta,\Delta+\delta}>t) -
\P^{B_1, B_2, B_3}_{0,\delta, \Delta+\delta}(\tau_{0,\delta}>t, \tau_{\delta,\Delta}>t),
\end{aligned}
\end{equation}
where we have reversed time and replaced $\hat w_{(a,t)}$ for $a=0,\delta, \Delta, \Delta+\delta$ by independent
standard Brownian motions $B_1, B_2$ and $B_3$ starting from $0$, $\delta$, and $\Delta$ (or $\Delta+\delta$), and
$\tau_{a,b}$ is the first hitting time between the two Brownian motions starting from $a$ and $b$.

By the Karlin-McGregor formula \cite{KM59}, the transition density for three one-dimensional Brownian motions
$B_1,B_2$ and $B_3$ starting from $x_1<x_2<x_3$ to end at locations $y_1<y_2<y_3$ at time $t$ without ever intersecting
along the way is given by the determinant ${\rm Det}(p_t(x_i,y_j)_{1\leq i,j\leq 3})$, where
$p_t(x,y)=\frac{e^{-\frac{(y-x)^2}{2t}}}{\sqrt{2\pi t}}$. Therefore we define

\be{rs10b1}
D:= \P^{B_1, B_2, B_3}_{0,\delta, \Delta}(\tau_{0,\delta}>t, \tau_{\delta,\Delta+\delta}>t) -
\P^{B_1, B_2, B_3}_{0,\delta, \Delta+\delta}(\tau_{0,\delta}>t, \tau_{\delta,\Delta}>t).
\ee
The r.h.s. can be written as:
\be{rs10b2}
=
\iiint\limits_{y_1<y_2<y_3}   \nonumber
\begin{vmatrix}
p_t(0,y_1) & p_t(0,y_2) & p_t(0,y_3) \\
p_t(\delta, y_1) & p_t(\delta, y_2) & p_t(\delta, y_3) \\
p_t(\Delta+\delta, y_1) & p_t(\Delta+\delta, y_2) & p_t(\Delta+\delta, y_3)
\end{vmatrix}
-
\begin{vmatrix}
p_t(0,y_1) & p_t(0,y_2) & p_t(0,y_3) \\
p_t(\delta, y_1) & p_t(\delta, y_2) & p_t(\delta, y_3) \\
p_t(\Delta, y_1) & p_t(\Delta, y_2) & p_t(\Delta+, y_3)
\end{vmatrix}
d\vec y.
\ee
After some elementary manipulation, we obtain that $D$ is given by the expression:
\begin{equation}
D =
\!\!\! \iiint\limits_{y_1<y_2<y_3}
\frac{e^{-\frac{y_1^2+y_2^2+y_3^2+\delta^2+\Delta^2}{2t}}}{(2\pi  t)^{\frac{3}{2}}}
\begin{vmatrix}
1 & 1 & 1 \\
e^{\frac{\delta y_1}{t}}-1 & e^{\frac{\delta y_2}{t}}-1  &   e^{\frac{\delta y_3}{t}}- 1 \\
e^{\frac{\Delta y_1}{t}} f(t,y_1,\delta,\Delta) &
e^{\frac{\Delta y_2}{t}} f(t,y_2,\delta,\Delta) &
e^{\frac{\Delta y_3}{t}} f(t,y_3,\delta,\Delta)
\end{vmatrix}
d\vec y,
\end{equation}
where $f(t,y_i,\delta, \Delta) = e^{\frac{2(y_i-\Delta)\delta-\delta^2}{2t}}-1$.
We then Taylor expand in $\delta$, and note that the factor $e^{-\frac{\Vert y\Vert^2}{2t}}$ allows us to take $D\delta^{-2}$
and pass the limit $\delta\downarrow 0$ inside the integral to obtain
\be{detform2}
K^{\rm c}_t(\Delta)
= \liml_{\delta\downarrow 0} \frac{D}{\delta^2} = \!\!
\iiint\limits_{y_1<y_2<y_3}  \frac{e^{-\frac{y_1^2+y_2^2+y_3^2+\Delta^2}{2t}}}{(2\pi)^{\frac{3}{2}} t^{\frac{7}{2}}}
\begin{vmatrix}
1        &        1       &         1  \\
y_1      &       y_2      &        y_3 \\
e^{\frac{\Delta y_1}{t}}(y_1-\Delta) & e^{\frac{\Delta y_2}{t}}(y_2-\Delta) & e^{\frac{\Delta y_3}{t}}(y_3-\Delta)
\end{vmatrix}
d\vec y.
\ee
Expanding the determinant and performing the change of variable $y_i=x_i\sqrt{t}$ and $\Delta=\bar\Delta\sqrt{t}$
gives $K^{\rm c}_t = \big((2\pi)^{\frac{3}{2}}t\big)^{-1} (I_1 + I_2 + I_3)$, where
\be{rs5}
\begin{aligned}
I_1 & = \iiint\limits_{x_1<x_2<x_3} e^{-\frac{(x_1-\bar\Delta)^2+x_2^2+x_3^2}{2}}(x_1-\bar\Delta)(x_3-x_2) d\vec x, \\
I_2 & = \iiint\limits_{x_1<x_2<x_3} e^{-\frac{x_1^2+(x_2-\bar\Delta)^2+x_3^2}{2}}(x_2-\bar\Delta)(x_1-x_3) d\vec x, \\
I_3 & = \iiint\limits_{x_1<x_2<x_3} e^{-\frac{x_1^2+x_2^2+(x_3-\bar\Delta)^2}{2}}(x_3-\bar\Delta)(x_2-x_1) d\vec x. \\
\end{aligned}
\ee
We then have
\be{rs6}
\begin{aligned}
I_1+I_3 & = -\iint\limits_{x_2<x_3} e^{-\frac{x_2^2+(x_2-\bar\Delta)^2+x_3^2}{2}}(x_3-x_2)dx_2dx_3
+ \iint\limits_{x_1<x_2} e^{-\frac{x_1^2+x_2^2+(x_2-\bar\Delta)^2}{2}}(x_2-x_1)dx_1dx_2  \\
& = \iint\limits_{\R^2} e^{-\frac{x_1^2+x_2^2+(x_2-\bar\Delta)^2}{2}}(x_2-x_1)dx_1dx_2 = \frac{\pi}{\sqrt2}\bar\Delta
e^{-\frac{\bar\Delta^2}{4}},
\end{aligned}
\ee
and
\be{rs7}
\begin{aligned}
I_2 & = \iint\limits_{x_1<x_3} e^{-\frac{x_1^2+x_3^2}{2}}
\big(e^{-\frac{(x_1-\bar\Delta)^2}{2}}-e^{-\frac{(x_3-\bar\Delta)^2}{2}}\big) (x_1-x_3) dx_1 dx_3 \\
& = \iint\limits_{\R^2} e^{-\frac{x_1^2+(x_1-\bar\Delta)^2+x_3^2}{2}}(x_1-x_3)dx_1dx_3 = \frac{\pi}{\sqrt2}\bar\Delta
e^{-\frac{\bar\Delta^2}{4}}.
\end{aligned}
\ee
Together, they give $K^{\rm c}_t = \frac{\Delta e^{-\frac{\Delta^2}{4t}}}{2\sqrt{\pi}t^{\frac{3}{2}}}$.
\qed

\noindent
{\bf Acknowledgements.} R.~Sun is supported by AcRF Tier 1 grant R-146-000-185-112.
A. Greven was supported by the Deutsche Forschungsgemeinschaft DFG, Grant GR 876/15-1 und 15-2.
{A.~Winter is supported by DFG SPP 1590.  }

\bi

\bibliography{gsw}
\bibliographystyle{alpha}

\end{document}